\documentclass[a4paper,12pt]{article}
\usepackage{amsfonts,amsmath,amssymb}
\newtheorem{theorem}{Theorem}[section]
\newtheorem{lemma}[theorem]{Lemma}
\newtheorem{proposition}[theorem]{Proposition}

\newtheorem{remark}[theorem]{Remark}

\numberwithin{equation}{section}

\newcommand{\QED}{\hspace*{\fill}\rule{2.5mm}{2.5mm}}

\newcommand{\bea}{\begin{eqnarray}}
\newcommand{\eea}{\end{eqnarray}}
\def\be#1{\begin{equation} \label{#1}}
\def \eeq{\end{equation}}

\newcommand{\ds}{\displaystyle}
\providecommand{\norm}[1]{\lVert#1\rVert}
\providecommand{\normm}[1]{\left\lVert#1\right\rVert}
\newcommand{\nn}{\nonumber}

\newcommand{\R}{\mathbb{R}}
\newcommand{\N}{\mathbb{N}}

\newcommand{\la}{\lambda}
\newcommand{\f}{\mathcal{F}^{-1}}
\newcommand{\trt}{\textrm{tr}\theta}
\newcommand{\dd}{{\bf D}}
\renewcommand{\o}{\omega}
\renewcommand{\c}{\cdot}
\renewcommand{\th}{\theta}
\newcommand{\ep}{\varepsilon}

\newcommand{\s}{\Sigma}
\renewcommand{\S}{\mathbb{S}^2}
\newcommand{\xo}{x\cdot\omega}
\newcommand{\muu}{\mu_u}
\newcommand{\p}{P_u}
\newcommand{\gn}{g(N,N')}
\newcommand{\h}[1]{H^{#1}(\R^3)}
\renewcommand{\l}[2]{L^{#1}_uL^{#2}(P_u)}
\newcommand{\lprime}[2]{L^{#1}_{u'}L^{#2}(P_{u'})}
\renewcommand{\ll}[1]{L^{#1}(\Sigma)}
\renewcommand{\le}[1]{L^{#1}(\R^3)}

\newcommand{\nabb}{\mbox{$\nabla \mkern-13mu /$\,}}

\newcommand{\nabn}{\nabla_N}
\renewcommand{\lg}{\log(a)}
\newcommand{\po}{\partial_{\omega}}
\renewcommand{\a}{\alpha}

\newcommand{\ga}{\gamma}
\newcommand{\de}{\delta}

\begin{document}

\begin{center}
\Large{\bf Parametrix for wave equations on a rough background II: construction and control at initial time}
\end{center}

\vspace{0.5cm}

\begin{center}
\large{J\'er\'emie Szeftel}
\end{center}

\vspace{0.2cm}

\begin{center}
\large{DMA, Ecole Normale Sup\'erieure,\\
45 rue d'Ulm, 75005 Paris,\\
jeremie.szeftel@ens.fr}
\end{center}

\vspace{0.5cm}

{\bf Abstract.} This is the second of a sequence of four papers  \cite{param1}, \cite{param2}, \cite{param3}, \cite{param4} dedicated to the construction and the control of a parametrix to the homogeneous wave equation $\square_{\bf g} \phi=0$, where ${\bf g}$ is a rough metric satisfying the Einstein vacuum equations. Controlling such a parametrix as well as its error term when one only assumes $L^2$ bounds on the curvature tensor ${\bf R}$ of ${\bf g}$ is a major step of the proof of the bounded $L^2$ curvature conjecture proposed in \cite{Kl:2000}, and solved by S. Klainerman, I. Rodnianski and the author in \cite{boundedl2}. On a more general level, this sequence of papers deals with the control of the eikonal equation on a rough background, and with the  derivation of $L^2$ bounds for Fourier integral operators on manifolds with rough phases and symbols, and as such is also of independent interest.

\vspace{0.2cm}

\section{Introduction}

We consider the Einstein vacuum equations,
\begin{equation}\label{eq:I1}
{\bf R}_{\alpha\beta}=0
\end{equation}
where ${\bf R}_{\alpha\beta}$
denotes the  Ricci curvature tensor  of  a four dimensional Lorentzian space time  $(\mathcal{M},\,  {\bf g})$. The Cauchy problem consists in finding a metric ${\bf g}$ satisfying \eqref{eq:I1} such that the metric induced by ${\bf g}$ on a given space-like hypersurface $\Sigma_0$ and the second fundamental form of $\Sigma_0$ are prescribed. The initial data then consists of a Riemannian three dimensional metric $g_{ij}$ and a symmetric tensor $k_{ij}$ on the space-like hypersurface $\Sigma_0=\{t=0\}$. Now, \eqref{eq:I1} is an overdetermined system and the initial data set $(\Sigma_0,g,k)$ must satisfy the constraint equations
\begin{equation}\label{const}
\left\{\begin{array}{l}
\nabla^j k_{ij}-\nabla_i \textrm{Tr}k=0,\\
 R-|k|^2+(\textrm{Tr}k)^2=0, 
\end{array}\right.
\end{equation}
where the covariant derivative $\nabla$ is defined with respect to the metric $g$, $R$ is the scalar curvature of $g$, and $\textrm{Tr}k$ is the trace of $k$ with respect to the metric $g$.

The fundamental problem in general relativity is to 
 study the long term regularity and asymptotic 
properties
 of the Cauchy developments of general, asymptotically flat,  
initial data sets $(\Sigma_0, g, k)$. As far as local regularity is concerned it 
is natural to ask what are the minimal regularity properties of
the initial data which guarantee the existence and 
uniqueness of local developments. In \cite{boundedl2}, we obtain the 
following result which solves bounded $L^2$ curvature conjecture proposed 
in \cite{Kl:2000}:

\begin{theorem}[Theorem 1.10 in \cite{boundedl2}]\label{th:mainbl2}
Let $(\mathcal{M}, {\bf g})$ an asymptotically flat solution to the Einstein vacuum equations \eqref{eq:I1} together with a maximal foliation by space-like hypersurfaces $\Sigma_t$ defined as level hypersurfaces of a time function $t$. Let $r_{vol}(\Sigma_t,1)$ the volume radius on scales $\leq 1$ of $\Sigma_t$\footnote{See Remark \ref{rem:volrad} below for a  definition}. Assume that the initial slice $(\Sigma_0,g,k)$ is such that:
$$\norm{R}_{L^2(\Sigma_0)}\leq \ep,\,\norm{k}_{L^2(\Sigma_0)}+\norm{\nabla k}_{L^2(\Sigma_0)}\leq \ep\textrm{ and }r_{vol}(\Sigma_0,1)\geq \frac{1}{2}.$$
Then, there exists a small universal constant $\ep_0>0$ such that if $0<\ep<\ep_0$, then the following control holds on $0\leq t\leq 1$:
$$\norm{\R}_{L^\infty_{[0,1]}L^2(\Sigma_t)}\lesssim \ep,\,\norm{k}_{L^\infty_{[0,1]}L^2(\Sigma_t)}+\norm{\nabla k}_{L^\infty_{[0,1]}L^2(\Sigma_t)}\lesssim \ep\textrm{ and }\inf_{0\leq t\leq 1}r_{vol}(\Sigma_t,1)\geq \frac{1}{4}.$$
\end{theorem}

\begin{remark}
While  the first   nontrivial improvements  for well posedness for quasilinear  hyperbolic systems (in spacetime dimensions greater than $1+1$), based on Strichartz estimates,  were obtained in   \cite{Ba-Ch1}, \cite{Ba-Ch2}, \cite{Ta1}, \cite{Ta}, \cite{KR:Duke}, \cite{KR:Annals}, \cite{SmTa}, Theorem \ref{th:mainbl2}, is the first  result in which the  full nonlinear structure of the quasilinear system, not just its principal part,  plays  a  crucial  role.  We note that  though  the result is not  optimal with respect to the  standard  scaling  of the Einstein equations, it is  nevertheless critical   with respect to its causal geometry,  i.e. $L^2 $ bounds on the curvature is the minimum requirement necessary to obtain lower bounds on the radius of injectivity of null hypersurfaces. We refer the reader to section 1 in \cite{boundedl2} for more motivations and historical perspectives concerning Theorem \ref{th:mainbl2}. 
\end{remark}

\begin{remark}
The regularity assumptions on $\Sigma_0$ in Theorem \ref{th:mainbl2} - i.e. $R$ and $\nabla k$ bounded in $L^2(\Sigma_0)$ - correspond to an initial data set $(g,\, k )\in H^2_{loc}(\Sigma_0)\times H^1_{loc}(\Sigma_0)$.
\end{remark}

\begin{remark}\label{rem:reducsmallisok}
In \cite{boundedl2}, our main result is stated for corresponding large data. We then reduce the proof to the  small data statement of Theorem \ref{th:mainbl2} relying on a truncation and rescaling procedure, the control of the harmonic radius of $\Sigma_0$ based on Cheeger-Gromov convergence of Riemannian manifolds together with the assumption on the lower bound of the volume radius of $\Sigma_0$, and the gluing procedure in \cite{Co}, \cite{CoSc}. We refer the reader to section 2.3 in \cite{boundedl2} for the details.
\end{remark}

\begin{remark}\label{rem:volrad}
We recall for the convenience of the reader the definition of the volume radius of the Riemannian manifold $\Sigma_t$. Let $B_r(p)$ denote the geodesic ball of center $p$ and radius $r$. The volume radius $r_{vol}(p,r)$ at a point $p\in \Sigma_t$ and scales $\leq r$ is defined by
$$r_{vol}(p,r)=\inf_{r'\leq r}\frac{|B_{r'}(p)|}{r^3},$$
with $|B_r|$ the volume of $B_r$ relative to the metric $g_t$ on $\Sigma_t$. The volume radius $r_{vol}(\Sigma_t,r)$ of $\Sigma_t$ on scales $\leq r$ is the infimum of $r_{vol}(p,r)$ over all points $p\in \Sigma_t$.
\end{remark}

The proof of Theorem \ref{th:mainbl2}, obtained in the sequence of papers \cite{boundedl2}, \cite{param1}, \cite{param2}, \cite{param3}, \cite{param4}, \cite{bil2}, relies on the following ingredients\footnote{We also need trilinear estimates and an $L^4(\mathcal{M})$ Strichartz estimate (see the introduction in \cite{boundedl2})}: 
{\em\begin{enumerate}
\item[{\bf A}] Provide  a system of coordinates relative to which \eqref{eq:I1} exhibits a null structure.

\item[{\bf B}] Prove  appropriate bilinear estimates for solutions to $\square_{\bf g} \phi=0$, on
 a fixed Einstein vacuum  background\footnote{Note that the first bilinear estimate of this type was obtained in \cite{BIL}}.

\item[{\bf C}] Construct a parametrix for solutions to the homogeneous wave equations $\square_{\bf g} \phi=0$ on a fixed Einstein vacuum  background, and obtain control of the parametrix and of its error term only using the fact that the curvature tensor is bounded in $L^2$. 
\end{enumerate}
}

Steps {\bf A} and {\bf B} are carried out in \cite{boundedl2}. In particular, the proof of the bilinear estimates rests on a representation formula for the solutions of the wave equation using the following plane wave parametrix:
\begin{equation}\label{param}
Sf(t,x)=\int_{\S}\int_{0}^{+\infty}e^{i\lambda u(t,x,\o)}f(\lambda\o)\lambda^2 d\lambda d\o,\,(t,x)\in\mathcal{M} 
\end{equation}
where $u(.,.,\o)$ is a solution to the eikonal equation ${\bf g}^{\alpha\beta}\partial_\alpha u\partial_\beta u=0$ on $\mathcal{M}$ such that $u(0,x,\o)\sim x.\o$ when $|x|\rightarrow +\infty$ on $\Sigma_0$.

\begin{remark}
Actually, \eqref{param} only corresponds to a half wave parametrix. The full parametrix will be derived in section \ref{sec:fullparam}. 
\end{remark}

\begin{remark}
The asymptotic behavior for $u(0,x,\o)$ when $|x|\rightarrow +\infty$ will be important to generate arbitrary initial data for the wave equation (see \eqref{choicef}).
\end{remark}

\begin{remark}
Note that the parametrix \eqref{param} is invariantly defined\footnote{Our choice is reminiscent of the one used in \cite{SmTa} in the context of $H^{2+\epsilon}$ solutions of quasilinear wave equations. Note however that the construction in that paper is coordinate dependent}, i.e. without reference to any coordinate system. This is crucial since coordinate systems consistent with $L^2$ bounds on the curvature would not be regular enough to control a parametrix. 
\end{remark}

In order to complete the proof of the bounded $L^2$ curvature conjecture, we need to carry out step {\bf C} with the parametrix defined in \eqref{param}.  

\begin{remark}
In addition to their relevance to the resolution of the bounded $L^2$ curvature conjecture, the methods and results  of step {\bf C} are  also of independent interest. Indeed, they deal on the one hand with the control of the eikonal equation ${\bf g}^{\alpha\beta}\partial_\alpha u\partial_\beta u=0$ at a critical level\footnote{We need at least $L^2$ bounds on the curvature to obtain a lower bound on the radius of injectivity of the null level hypersurfaces of the solution $u$ of the eikonal equation, which in turn is necessary to control the local regularity of $u$ (see \cite{param3})}, and on the other hand with the derivation of $L^2$ bounds for Fourier integral operators with significantly lower differentiability  assumptions both for the corresponding phase and symbol compared to classical methods (see for example \cite{stein} and references therein). 
\end{remark}

In view of the energy estimates for the wave equation, it suffices to control the parametrix at $t=0$ (i.e. restricted to $\Sigma_0$)
\be{parami}
Sf(0,x)=\int_{\S}\int_{0}^{+\infty}e^{i\lambda u(0,x,\o)}f(\lambda\o)\lambda^2 d\lambda d\o,\,x\in\Sigma_0 
\eeq
 and the error term
\be{err} 
Ef(t,x)=\square_{\bf g}Sf(t,x)=\int_{\S}\int_{0}^{+\infty}e^{i\lambda u(t,x,\o)}\square_{\bf g}u(t,x,\o)f(\lambda\o)\lambda^3 d\lambda d\o,\,(t,x)\in\mathcal{M}. 
\eeq
This requires the following ingredients, the two first being related to the control of the parametrix restricted to $\Sigma_0$ \eqref{parami}, and the two others being related to the control of the error term \eqref{err}:
{\em\begin{enumerate}
\item[{\bf C1}] Make an appropriate choice for the equation satisfied by $u(0,x,\o)$ on $\Sigma_0$, and control the geometry of the foliation generated by the level surfaces of $u(0,x,\o)$ on $\Sigma_0$.

\item[{\bf C2}] Prove that the parametrix at $t=0$ given by \eqref{parami} is bounded in $\mathcal{L}(L^2(\mathbb{R}^3),\ll{2})$ using the estimates for $u(0,x,\o)$ obtained in {\bf C1}.

\item[{\bf C3}] Control the geometry of the foliation generated by the level hypersurfaces of $u$ on $\mathcal{M}$.

\item[{\bf C4}] Prove that the error term \eqref{err} satisfies the estimate $\norm{Ef}_{L^2(\mathcal{M})}\leq C\norm{\lambda f}_{L^2(\mathbb{R}^3)}$ using the estimates for $u$ and $\square_{{\bf g}}u$ proved in {\bf C3}.
\end{enumerate}
}

Concerning step {\bf C1}, let us note that the typical choice $u(0,x,\o)=x\c\o$ in a given coordinate system would not work for us, since we don't have enough control on the regularity of a given coordinate system within our framework\footnote{This issue appears because we are working at the level of $H^2$ solutions for Einstein equations. In particular, the choice $u(0,x,\o)=x\c\o$ in a given coordinate system is used in \cite{SmTa} in the context of $H^{2+\epsilon}$ solutions for quasilinear wave equations}. Instead, in \cite{param1}, we rely on a geometric definition for $u(0,x,\o)$ to achieve step {\bf C1}. In the present paper, we focus on step {\bf C2}. 

Note that the parametrix at $t=0$ given by \eqref{parami} is a Fourier integral operator (FIO) with phase $u(0,x,\o)$. Now, we only assume $R\in L^2(\Sigma_0)$ and $\nabla k\in L^2(\Sigma_0)$ in order to be consistent with the statement of Theorem \ref{th:mainbl2}. This severely limits the regularity we are able to obtain in step {\bf C1} for $u(0,x,\o)$ (see \cite{param1} and section \ref{sec:paulailey}). Although $R$ and $k$ do not depend on the parameter $\o$, the regularity in $\o$ we are able to obtain in step {\bf C1} for $u(0,x,\o)$ is very limited as well\footnote{This is due to the fact that our estimates are better in directions tangent to the $u$-foliation on $\Sigma_0$. Now, after differentiation with respect to $\o$, derivatives in tangential directions pick up a nonzero component along the normal direction to the $u$-foliation on $\Sigma_0$ (see \cite{param1} for details)}. In particular, we obtain for the phase $u(0,x,\o)$ of $Sf(x,0)$ in \eqref{parami}\footnote{Actually, we have weaker bounds for the estimates where all the spatial derivatives are taken in the direction normal to the $u$-foliation on $\Sigma_0$ (see section \ref{sec:paulailey})}:
\begin{equation}\label{cestpasbeaucoup}
\sup_\o\left(\norm{\nabla^3u}_{L^2(\Sigma_0)}+\norm{\nabla\po u}_{L^\infty(\Sigma_0)}+\norm{\nabla^2\po u}_{L^2(\Sigma_0)}+\norm{\po^3u}_{L^\infty_{\textrm{loc}}(\Sigma_0)}\right)\lesssim\ep.
\end{equation}
Let us note that the classical arguments for proving $L^2$ bounds for FIO are based either on a $T T^*$ argument, or a $T^* T$ argument, which requires in our setting\footnote{Since $\Sigma_0$ is 3-dimensional} taking at least 4 derivatives of the phase in $L^\infty(\Sigma_0\times\S)$ either with respect to $x$ for $T^*T$, or with respect to $(\la, \o)$ for $TT^*$ (see for example \cite{stein}). Both methods would fail by a large margin, in particular in view of the regularity \eqref{cestpasbeaucoup} obtained for the phase of the parametrix at initial time $Sf(x,0)$. In order to obtain the control required in step {\bf C2} with the regularity of the phase of the FIO $Sf(x,0)$ given by \eqref{cestpasbeaucoup}, we are forced to design a method which allows us to take advantage both of the regularity in $x$ and $\o$. This is achieved using in particular the following ingredients:
\begin{itemize}
\item geometric integrations by parts taking full advantage of the better regularity properties in directions tangent to the level surfaces of $u(0,x,\o)$\footnote{Let us repeat that we actually obtain a weaker bound than \eqref{cestpasbeaucoup} for the estimates where all the spatial derivatives are taken in the direction normal to the $u$-foliation on $\Sigma_0$ (see section \ref{sec:paulailey})},

\item the standard first and second dyadic decomposition in frequency and angle (see \cite{stein}), as well as another decomposition involving frequency and angle,

\item after localization in frequency and angle, an estimate for the diagonal term using the $TT^*$ argument and a  change of variable tied to $u(0,x,\o)$. 
\end{itemize}

\vspace{0.2cm}

The rest of the paper is as follows. In section 2, we present the full parametrix for solutions to the homogeneous wave equation $\square_{\bf g}\phi=0$, we recall the regularity for the phase $u(0,x,\o)$ obtained in \cite{param1}, and we state our main results. In section 3, we prove the boundedness on $L^2$ of a pseudodifferential operator acting on $\R^3$ with a rough symbol introducing the main ideas in a simple setting. In section 4, we prove the boundedness on $L^2$ of a Fourier integral operator acting on $\Sigma_0$ with phase 
$u(0,x,\o)$ and a symbol having limited regularity consistent with the one given by our parametrix. Finally, we use the results of section 4 to show the existence and to control our parametrix in section 5.\\

\noindent{\bf Acknowledgments.} The author wishes to express his deepest gratitude to Sergiu Klainerman and Igor Rodnianski for stimulating discussions and constant encouragements during the long years where this work has matured. He also would like to stress that  the basic strategy of the construction of the parametrix and how it fits  into the whole proof of the bounded $L^2$ curvature conjecture has been done in collaboration with them. Finally, he  would like to mention the influential work \cite{SmTa} providing construction and control of parametrices for $H^{2+\epsilon}$ solutions of quasilinear wave equations. The author is supported by ANR jeunes chercheurs SWAP.\\

\section{Main results}

From now on, there will be no further reference to $\Sigma_t$ for $t>0$. Since there is no confusion, we will denote $\Sigma_0$ simply by $\Sigma$ in the rest of the paper. 

\subsection{Presentation of the parametrix}\label{sec:fullparam}

In this section, we construct a parametrix for the following homogeneous wave equation:
\begin{equation}\label{waveeq}
\left\{\begin{array}{l}
\ds\square_{\bf g}\phi=0\textrm{ on }\mathcal{M},\\
\ds\phi_{|_\s}=\phi_0,\, T(\phi)_{|_\s}=\phi_1,
\end{array}\right.
\end{equation}
where $\phi_0$ and $\phi_1$ are two given functions on $\s$ and $T$ is the future oriented unit normal to $\s$ in $\mathcal{M}$.

We recall the plane wave representation of the solution of the flat wave equation. This corresponds to the case where ${\bf g}$ is the Minkowski metric. \eqref{waveeq} becomes:
\begin{equation}\label{flatwaveeq}
\left\{\begin{array}{l}
\ds\square\phi=0\textrm{ on }\R^{1+3},\\
\ds\phi(0,.)=\phi_0,\, \partial_t\phi(0,.)=\phi_1\textrm{ on }\R^3.
\end{array}\right.
\end{equation}
The plane wave representation of the solution $\phi$ of \eqref{flatwaveeq} is given by:
\begin{equation}\label{flatparam}
\begin{array}{l}
\ds\int_{\S}\int_0^{+\infty}e^{i(-t+\xo)\la}\frac{1}{2}\left(\mathcal{F}\phi_0(\la\o)+i\frac{\mathcal{F}\phi_1(\la\o)}{\la}\right)d\la d\o\\
\ds +\int_{\S}\int_0^{+\infty}e^{i(t+\xo)\la}\frac{1}{2}\left(\mathcal{F}\phi_0(\la\o)-i\frac{\mathcal{F}\phi_1(\la\o)}{\la}\right)d\la d\o,
\end{array}
\end{equation}
where $\mathcal{F}$ denotes the Fourier transform on $\R^3$.

We would like to construct a parametrix in the curved case similar to \eqref{flatparam}. We introduce two solutions $u_\pm$ of the eikonal equation
\begin{equation}\label{eikonal}
{\bf g}^{\alpha\beta}\partial_\alpha u_\pm\partial_\beta u_\pm=0\textrm{ on }\mathcal{M},
\end{equation}
such that:
\begin{equation}\label{eikonal1}
T(u_\pm)=\mp |\nabla u_\pm| =\mp a^{-1}_\pm\textrm{ on }\s,
\end{equation}
where $T$ is the future oriented unit normal to $\s$ in the space-time $\mathcal{M}$, $\nabla$ is the gradient on $\s$ associated to the metric $g$, $|\c|$ is the length associated to $g$ for vectorfields on $\s$, and $a_\pm$ is the lapse of $u_\pm$ on $\s$. We look for a parametrix for \eqref{waveeq} of the form:
\be{param1}
\begin{array}{ll}
\ds S_+f_+(t,x)+S_-f_-(t,x)= & \ds\int_{\S}\int_{0}^{+\infty}e^{i\lambda u_+(t,x,\o)}f_+(\lambda\o)\lambda^2 d\lambda d\o\\
& \ds +\int_{\S}\int_{0}^{+\infty}e^{i\lambda u_-(t,x,\o)}f_-(\lambda\o)\lambda^2 d\lambda d\o,\,(t,x)\in\mathcal{M}.
\end{array}
\eeq
Thanks to \eqref{eikonal}, this parametrix generates the following error term:
\be{error1}
\begin{array}{r}
\ds E_+f_+(t,x)+E_-f_-(t,x)=\int_{\S}\int_{0}^{+\infty}e^{i\lambda u_+(t,x,\o)}\square_gu_+(t,x,\o)f_+(\lambda\o)\lambda^3 d\lambda d\o\\
\ds +\int_{\S}\int_{0}^{+\infty}e^{i\lambda u_-(t,x,\o)}\square_gu_-(t,x,\o)f_-(\lambda\o)\lambda^3 d\lambda d\o,\,(t,x)\in\mathcal{M}.
\end{array}
\eeq
In the next two sections, we precise the parametrix \eqref{param1} by prescribing $u_\pm$ on $\s$ and by making our choice for $f_\pm$ explicit.

\subsubsection{Prescription of $u_+$ and $u_-$ on $\s$}

\eqref{eikonal} and \eqref{eikonal1} are not enough to define $u_\pm$ in a unique manner. Indeed, we 
still need to prescribe $u_\pm$ on $\s$. To motivate our choice, we need to introduce some geometric 
objects connected to $u_\pm$. Let $N_\pm$ the vectorfield on $\s$ defined by:
\begin{equation}\label{eikonal2}
N_\pm=\frac{\nabla u_\pm}{|\nabla u_\pm|}=a_\pm\nabla u_\pm,
\end{equation}
and $L_\pm$ the vectorfield on $\mathcal{M}$ which is given on $\s$ by: 
\begin{equation}\label{eikonal3}
L_\pm=a_\pm{\bf g}^{\alpha\beta}\partial_\alpha u_\pm\partial_\beta=a_\pm(-T(u_\pm)T+\nabla u_\pm)=\pm T+N_\pm.
\end{equation}
Let $P_{u_\pm}=\{x\in\s/\,u_\pm(x)=u_\pm\}$ denote the level surfaces of $u_\pm$ in $\s$. Since $N_\pm$ is the unit normal to $P_{u_\pm}$, the second fundamental form of $P_{u_\pm}$ in $\s$ is given by:
\begin{equation}\label{eikonal4}
\th_\pm(e^\pm_A,e^\pm_{B})=g(D_{e^\pm_A}N_\pm,e^\pm_{B}),\,A, B=1,2,
\end{equation}
where $(e^\pm_1,e^\pm_2)$ is an arbitrary orthonormal frame of $TP_{u_\pm}$. Let$$\mathcal{H}_{u_\pm}=\{(t,x)\in\mathcal{M}/\,u_\pm(t,x)=u_\pm\}$$ 
denote the null level hypersurfaces of $u_\pm$ in $\mathcal{M}$. Since $L_\pm$ is null and orthogonal to $P_{u_\pm}$ in $\mathcal{H}_{u_\pm}$, the null second fundamental form $\chi_\pm$ is given on $P_{u_\pm}$ by:
\begin{equation}\label{eikonal5}
\chi_\pm(e^\pm_A,e^\pm_{B})=g(\dd_{e^\pm_A}L_\pm,e^\pm_{B}),\,A, B=1,2.
\end{equation}
Taking the trace in \eqref{eikonal4} and \eqref{eikonal5}, and using \eqref{eikonal3} and the fact that $k$ is the second fundamental form of $\s$, we obtain:
\begin{equation}\label{eikonal6}
\textrm{tr}\chi_\pm=\pm\textrm{tr}k+\textrm{tr}\th_\pm.
\end{equation}
Note that $\textrm{Tr}k=\textrm{tr}k+k_{NN}$, where Tr denotes the trace for 2-tensors on $\s$. In addition to the constraint equations \eqref{const}, we choose a maximal foliation to be consistent with the statement of Theorem \ref{th:mainbl2}. This corresponds to $\textrm{Tr}k=0$. Together with \eqref{eikonal6}, this yields:
\begin{equation}\label{eikonal7}
\textrm{tr}\chi_\pm=\mp k_{N_\pm N_\pm}+\textrm{tr}\th_\pm.
\end{equation}
Now, an easy computation yields:
\begin{equation}\label{eikonal8}
\square_{{\bf g}}u_\pm=a^{-1}_\pm\textrm{tr}\chi_\pm,
\end{equation}
so that the error term \eqref{error1} may be rewritten:
\be{error2}
\begin{array}{ll}
& \ds E_+f_+(t,x)+E_-f_-(t,x)\\
\ds = & \ds\int_{\S}\int_{0}^{+\infty}e^{i\lambda u_+(t,x,\o)}a_+(t,x,\o)^{-1}\textrm{tr}\chi_+(t,x,\o)f_+(\lambda\o)\lambda^3 d\lambda d\o\\
& \ds +\int_{\S}\int_{0}^{+\infty}e^{i\lambda u_-(t,x,\o)}a_-(t,x,\o)^{-1}\textrm{tr}\chi_-(t,x,\o)f_-(\lambda\o)\lambda^3 d\lambda d\o,\,(t,x)\in\mathcal{M}.
\end{array}
\eeq
In view of \eqref{error2}, one has to show in particular that $\textrm{tr}\chi_\pm$ belongs to $L^\infty(\mathcal{M})$ as part of step {\bf C3} in order to complete step {\bf C4}. This estimate is obtained in \cite{FLUX} using a transport equation (the Raychadhouri equation). Thus, one needs the corresponding estimate on $\s$ (i.e. at $t=0$):
\begin{equation}\label{eikonal9}
\textrm{tr}\chi_\pm\in\ll{\infty},
\end{equation}
which in view of \eqref{eikonal7} is equivalent to:
\begin{equation}\label{eikonal10}
\mp k_{N_\pm N_\pm}+\textrm{tr}\th_\pm\in\ll{\infty}.
\end{equation}
Now, we construct in \cite{param1} a function $u(x,\o)$ on $\s\times\S$ such that 
\begin{equation}\label{eikonal11}
- k_{NN}+\textrm{tr}\th\in\ll{\infty}.
\end{equation}
Note that $-u(x,-\o)$ satisfies:
\begin{equation}\label{eikonal12}
k_{NN}+\textrm{tr}\th\in\ll{\infty}.
\end{equation}
Thus, in view of \eqref{eikonal10}, \eqref{eikonal11} and \eqref{eikonal12}, we initialize $u_\pm$ on $\s$ by:
\begin{equation}\label{eikonal13}
u_+(0,x,\o)=u(x,\o)\textrm{ and }u_-(0,x,\o)=-u(x,-\o)\textrm{ for }(x,\o)\in\s\times\S.
\end{equation}

\begin{remark}
Note that in the particular case where $k\equiv 0$ - the so-called time symmetric case-, we may take 
$$u_+(0,x,\o)=u_-(0,x,\o)=u(x,\o)\textrm{ for }(x,\o)\in\s\times\S.$$
In particular, we have $u_+(0,x,\o)=u_-(0,x,\o)=\xo$ in the flat case.
\end{remark}

\subsubsection{The choice of $f_+$ and $f_-$}\label{sec:choicef}

Having defined $u_\pm$, we still need to define $f_\pm$ in the parametrix \eqref{param1}. According to \eqref{waveeq}, the half wave parametrix $S_+$ and $S_-$ should satisfy on $\s$:
\begin{equation}\label{choicef}
\left\{\begin{array}{l}
\ds S_+f_+(0,x)+S_-f_-(0,x)=\phi_0(x),\\
\ds T(S_+f_+)(0,x)+T(S_-f_-)(0,x)=\phi_1(x).
\end{array}\right.
\end{equation}
Let us introduce the following operators acting on functions of $\R^3$:
\begin{equation}\label{choicef1}
M_\pm f(x)=\int_{\S}\int_{0}^{+\infty}e^{\pm i\lambda u(x,\pm\o)}f(\lambda\o)\lambda^2 d\lambda d\o
\end{equation}
and 
\begin{equation}\label{choicef2}
Q_\pm f(x)=\int_{\S}\int_{0}^{+\infty}e^{\pm i\lambda u(x,\pm\o)}a(x,\pm\o)^{-1} f(\lambda\o)\lambda^2 d\lambda d\o,
\end{equation}
where $a(x,\o)=|\nabla u(x,\o)|^{-1}$ is the lapse of $u$. 
Using \eqref{eikonal1}, the definition of $S_\pm$ in \eqref{param1}, \eqref{eikonal13}, the definition \eqref{choicef1} of $M_\pm$ and the definition \eqref{choicef2} of $Q_\pm$, we may rewrite \eqref{choicef} as: 
\begin{equation}\label{choicef3}
\left\{\begin{array}{l}
M_+f_++M_-f_-=\phi_0,\\
Q_+(\la f_+)-Q_-(\la f_-)=i\phi_1.
\end{array}\right.
\end{equation}
The goal of this paper will be to show that there exist a unique $(f_+,f_-)$ satisfying \eqref{choicef3}, 
and that $(f_+,f_-)$ satisfies the following estimate:
\begin{equation}\label{choicef4}
\norm{\la f_+}_{L^2(\R^3)}+\norm{\la f_-}_{L^2(\R^3)}\lesssim \norm{\nabla\phi_0}_{\ll{2}}+\norm{\phi_1}_{\ll{2}}.
\end{equation}

\begin{remark}
In the case of the flat wave equation \eqref{flatwaveeq}, we have $(\s,g)=(\R^3,\de)$, $u_\pm(t,x,\o)=\mp t+\xo$, $u(x,\o)=\xo$ and $a(x,\o)= 1$. In particular, the operators $M_\pm$ and $Q_\pm$ defined respectively by \eqref{choicef1} and \eqref{choicef2} all coincide with the inverse Fourier transform. Then, the system \eqref{choicef3} admits the following solutions:
$$f_\pm (\la\o) =\frac{1}{2}\left(\mathcal{F}\phi_0(\la\o)\pm i\frac{\mathcal{F}\phi_1(\la\o)}{\la} \right),$$
which clearly satisfy the estimate \eqref{choicef4}.
\end{remark}

Before stating precisely the main results of this paper, we will first recall the regularity obtained for the 
phase $u(x,\o)$ constructed in \cite{param1}. 

\subsection{Regularity assumptions on the phase $u(x,\o)$}\label{sec:paulailey}

The operators $M_\pm$ and $Q_\pm$ defined respectively in \eqref{choicef1} and \eqref{choicef2} are 
Fourier integral operators with phase $\pm u(x,\pm\o)$. The regularity assumptions on $u(x,\o)$ will be crucial to show the existence of $(f_+,f_-)$ satisfying \eqref{choicef3} and the estimate \eqref{choicef4}. 
In this section, we state our assumptions on $u(x,\o)$. 

We define the lapse $a(x,\o)=|\nabla u(x,\o)|^{-1}$, and the unit vector $N$ such that $\nabla u(x,\o)=a(x,\o)^{-1}N(x,\o)$. We also define the level surfaces $\p=\{x\,/\,u(x,\o) =u\}$ so that $N$ is the normal to $\p$. The second fundamental form $\th$ of $\p$ is defined by 
\be{eq:I6}
\th(X,Y)=g(\nabla_{X}N,Y)
\end{equation}
with $X,Y$ arbitrary vectorfields tangent to the  $u$-foliation $P_{u}$ of $\s$ and where $\nabla$ denotes the covariant differentiation with respect to $g$. We denote by $\trt$ the trace of $\th$,
i.e. $\trt=\delta^{AB}\th_{AB}$ where $\th_{AB}$ are the components of $\th$ relative to an
orthonormal frame $(e_A)_{A=1,2}$ on $P_u$. 

Let $\muu$ denote the area element of $\p$. Then, for all integrable function $f$ on $\s$, the coarea formula implies:
\begin{equation}\label{coarea}
\int_{\s} fd\s =\int_{u}\int_{\p} fad\muu du.
\end{equation}
It is also well-known that for a scalar function $f$:
\begin{equation}\label{du}
\frac{d}{du}\left(\int_{\p}fd\muu\right) =\int_{\p}\left(\frac{df}{du}+\trt f\right)d\muu.
\end{equation}
For $1\leq p,q\leq +\infty$, we define the spaces $\l{p}{q}$ using the norm 
$$\norm{F}_{\l{p}{q}}=\left(\int_{u}\norm{F}^p_{L^q(\p)}du\right)^{1/p}.$$
We assume that $1/2\leq a(x)\leq 2$ for all $x\in\s$ (see {\bf Assumption 1} below) so that $\l{p}{p}$ coincides with $L^p(\s)$ for all $1\leq p\leq +\infty$. We denote by $\gamma$ the metric induced by $g$ on $\p$, and by $\nabb$ the induced covariant derivative. 


We now state our assumptions for the phase $u(x,\o)$ of our Fourier integral operators. These assumptions are compatible with the regularity obtained for the function $u(x,\o)$ constructed in \cite{param1} (this construction corresponds to step {\bf C1}). The constant $\ep>0$ below satisfies $0<\ep<1$ and will be chosen later to be sufficiently small.

\vspace{0.5cm}

{\em \noindent{\bf Assumption 1} (regularity with respect to $x$):
\begin{equation}\label{thregx1}
\ds\norm{\nabla a}_{\l{\infty}{2}}+\norm{a-1}_{\ll{\infty}}+\norm{\nabb\nabla a}_{\ll{2}}+\norm{\th}_{\l{\infty}{2}}+\norm{\nabla\th}_{\ll{2}}\lesssim\ep.
\end{equation}

\noindent {\bf Assumption 2} (regularity with respect to $\o$): 

\begin{equation}\label{threomega1}
\ds\norm{\po a}_{\ll{2}}+\norm{\nabla\po a}_{\ll{2}}+\norm{\po\th}_{\ll{2}}+\norm{\nabla\po\th}_{\ll{2}}\lesssim\ep,
\end{equation}
\begin{equation}\label{threomega3}
\norm{\po^{\a}a}_{\ll{\infty}}\lesssim 1\textrm{ for some }0<\a<1.
\end{equation}
\begin{equation}\label{threomega1bis}
\norm{\po N}_{\ll{\infty}}\lesssim 1,
\end{equation}
\begin{equation}\label{threomega1ter}
||N(x,\o)-N(x,\o')|-|\o-\o'||\lesssim (\ep+|\o-\o'|)|\o-\o'|,\,\forall x\in\s, \o,\o'\in\S,
\end{equation}
\begin{equation}\label{threomega2}
\norm{\nabla\po^2N}_{\ll{2}}\lesssim\ep.
\end{equation}
and 
\begin{equation}\label{threomega3bis}
\norm{\po^3u}_{L^\infty_{\textrm{loc}}(\Sigma)}\lesssim 1.
\end{equation}

\noindent {\bf Assumption 3} (additional regularity with respect to $x$): 

\noindent For all $j\geq 0$, there are scalar functions $a^j_1$ and $a^j_2$ such that:
\begin{equation}\label{cordecfr1}
\begin{array}{l}
\ds\nabn a=a^j_1+a^j_2\textrm{ where }\norm{a^j_1}_{\ll{2}}\lesssim 2^{- j/2}\ep,\,\norm{a^j_2}_{\l{\infty}{2}}\lesssim \ep\\
\ds\textrm{and }\norm{\nabn a^j_2}_{\ll{2}}+\norm{a^j_2}_{\l{2}{\infty}}\lesssim 2^{j/2}\ep.
\end{array}
\end{equation}

\noindent {\bf Assumption 4} (global change of variable on $\s$):

\noindent Let $\o\in\S$. Let $\phi_\o:\s\rightarrow\R^3$ defined by:
\begin{equation}\label{gl21}
\phi_\o(x):=u(x,\o)\o+\po u(x,\o).
\end{equation}
Then $\phi_\o$ is a bijection, and the determinant of its Jacobian satisfies the following estimate:
\begin{equation}\label{gl22}
\norm{|\det(\textrm{Jac}\phi_\o)|-1}_{\ll{\infty}}\lesssim\ep.
\end{equation}

\noindent {\bf Assumption 5} (comparison of $u(x,\o)$ with a phase linear in $\o$):

\noindent Let $\nu\in\S$ and $\phi_\nu$ the map defined in \eqref{gl21}. Then, we have:
\begin{equation}\label{gl23}
\begin{array}{l}
\ds u(x,\o)-\phi_\nu(x)\c\o=O(\ep|\o-\nu|^2),\\
\ds \po u(x,\o)-\po(\phi_\nu(x)\c\o)=O(\ep|\o-\nu|),\\
\ds \po^2 u(x,\o)-\po^2(\phi_\nu(x)\c\o)=O(\ep).
\end{array}
\end{equation}

\noindent {\bf Assumption 6} (comparison of $N(x,-\o)$ with $N(x,\o)$):

\noindent For all $x\in\s$ and $\o\in\S$, we have:
 \begin{equation}\label{gl24}
 |N(x,\o)+N(x,-\o)|\lesssim\ep.
 \end{equation}}

\begin{remark}
In {\bf Assumptions 1-6}, all inequalities hold for any $\o\in\S$ with the constant in the right-hand side being independent of $\o$. Thus, one may take the supremum in $\o$ everywhere. To ease the notations, we do not explicitly write down this supremum. 
\end{remark}

\begin{remark}\label{smallreduc}
The fact that we may take a small constant $\ep>0$ in  {\bf Assumptions 1-6} is directly related to the assumptions on $\Sigma$ for $R$ and $k$ in Theorem \ref{th:mainbl2}.
\end{remark}

\begin{remark}
In the case of the flat wave equation \eqref{flatwaveeq}, we have $(\s,g)=(\R^3,\de)$, $u(x,\o)=\xo$, $a= 1$, $N=\o$ and $\phi_\o=I\! d_{\R^3}$. Thus, {\bf Assumptions 1-6} are clearly satisfied with $\ep=0$.
\end{remark}

\begin{remark}\label{reduccompact}
In \cite{param1}, the phase $u(x,\o)$ is actually exactly equal to $\xo$ on $|x|\geq 2$. This is made possible by exploiting the finite speed of propagation of Einstein vacuum equations (see \cite{param1}).
\end{remark}

\begin{remark}
Recall that the lapse $a$ is at the level of one derivative of $u$ with respect to $x$. Thus, we obtain from \eqref{thregx1} that some components of $\nabla^3u$ are in $\ll{2}$. Note that this is not true for all components since \eqref{cordecfr1} does not allow us to control $\nabn^2a$ in $\ll{2}$. In fact, \eqref{cordecfr1} is only at the level of $3/2$ derivatives of $a$ with respect to $N$ in $L^2$. 
\end{remark}



\subsection{Main results}

We first state a result of boundedness on $L^2$ for Fourier integral operators with phase $u(x,\o)$. 
\begin{theorem}\label{th1}
Let $u$ be a function on $\s\times\S$ satisfying {\bf Assumption 1}, {\bf Assumption 2} and {\bf Assumption 4}. Let $U$ the Fourier integral operator with phase $u(x,\o)$ and symbol $b(x,\o)$:
\begin{equation}\label{fio}
Uf(x)=\int_{\S}\int_{0}^{+\infty}e^{i\lambda u(x,\o)}b(x,\o)f(\lambda\o)\lambda^2 d\lambda d\o.
\end{equation}
Let $D>0$. We assume furthermore that $b(x,\o)$ satisfies:
\begin{equation}\label{thregx1s}
\norm{b}_{\ll{\infty}}+\norm{\nabla b}_{\l{\infty}{2}}+\norm{\nabb\nabla b}_{\ll{2}}\lesssim D,
\end{equation}
\begin{equation}\label{threomega1s}
\norm{\po b}_{\ll{2}}+\norm{\nabla\po b}_{\ll{2}}\lesssim D,
\end{equation}
and
\begin{equation}\label{cordecfr1s}
\begin{array}{l}
\ds\nabn b=b^j_1+b^j_2\textrm{ where }\norm{b^j_1}_{\ll{2}}\lesssim 2^{-\frac{j}{2}}D, \norm{b^j_2}_{\l{\infty}{2}}\lesssim D,\\
\ds\textrm{ and }\norm{\nabn b^j_2}_{\ll{2}}+\norm{b^j_2}_{\l{2}{\infty}}\lesssim 2^{\frac{j}{2}}D.
\end{array}
\end{equation}
Then, $U$ is bounded on $L^2$ and satisfies the estimate:
\begin{equation}\label{l2}
\norm{Uf}_{\ll{2}}\lesssim D\norm{f}_{L^2(\R^3)}.
\end{equation}
\end{theorem}

\begin{remark}
We intend to apply Theorem \ref{th1} to the Fourier integral operators $M_\pm$ and $Q_\pm$ introduced in section \ref{sec:choicef} whose symbol are respectively 1 and $a^{-1}$. Thus, our assumptions on the regularity of the symbol $b(x,\o)$ are consistent with the assumptions on the regularity of $a(x,\o)$ given by {\bf Assumptions 1-3}.  
\end{remark}

Recall the definition of the Fourier integral operators $M_\pm$ and $Q_\pm$ introduced in section \ref{sec:choicef}: 
\begin{equation}\label{choicef1bis}
M_\pm f(x)=\int_{\S}\int_{0}^{+\infty}e^{\pm i\lambda u(x,\pm\o)}f(\lambda\o)\lambda^2 d\lambda d\o,
\end{equation}
and 
\begin{equation}\label{choicef2bis}
Q_\pm f(x)=\int_{\S}\int_{0}^{+\infty}e^{\pm i\lambda u(x,\pm\o)}a(x,\pm\o)^{-1} f(\lambda\o)\lambda^2 d\lambda d\o.
\end{equation}
The following theorem is the main result of this paper and achieves step {\bf C2}.

\begin{theorem}\label{th2}
Let $u$ be a function on $\s\times\S$ satisfying {\bf Assumptions 1-6}. Then, there exist a unique $(f_+,f_-)$ satisfying:
\begin{equation}\label{choicef3bis}
\left\{\begin{array}{l}
M_+f_++M_-f_-=\phi_0,\\
Q_+(\la f_+)-Q_-(\la f_-)=i\phi_1.
\end{array}\right.
\end{equation}
Furthermore, $(f_+,f_-)$ satisfies the following estimate:
\begin{equation}\label{choicef4bis}
\norm{\la f_+}_{L^2(\R^3)}+\norm{\la f_-}_{L^2(\R^3)}\lesssim \norm{\nabla\phi_0}_{\ll{2}}+\norm{\phi_1}_{\ll{2}}.
\end{equation}
\end{theorem}

\begin{remark}
In view of the definition of $U$, $M_\pm$ and $Q_\pm$, the estimates \eqref{l2} and \eqref{choicef4bis} 
correspond to the obtention of $L^2$ bounds for Fourier integral operators. Let us repeat that the classical arguments for proving $L^2$ bounds for Fourier operators are based either on a $T T^*$ argument, or a $T^* T$ argument, 
which requires in our setting\footnote{Since $\Sigma_0$ is 3-dimensional} taking at least 4 derivatives of the phase in $L^\infty(\Sigma_0\times\S)$ either with respect to $x$ for $T^*T$, or with respect to $(\la, \o)$ for $TT^*$ (see for example \cite{stein}). Both methods would fail by far within the regularity for the phase $u(x,\o)$ given by {\bf Assumptions 1-4} and for the symbol $b(x,\o)$ given by \eqref{thregx1s} \eqref{threomega1s} \eqref{cordecfr1s}. 
\end{remark}

\subsection{Boundnessness on $L^2$ for pseudodifferential operators acting on $\R^3$ with rough symbols}

Theorem \ref{th1} yields the following result on the $L^2$ boundedness of pseudodifferential operators acting on $\R^3$ which corresponds to the case $\s=\R^3$, $g=\de$ and $u(x,\o)=\xo$. 
\begin{theorem}\label{th3}
Let $B$ the pseudodifferential operator with symbol $b(x,\o)$:
\begin{equation}\label{opd}
Bf(x)=\int_{\S}\int_{0}^{+\infty}e^{i\lambda \xo}b(x,\o)\mathcal{F}f(\lambda\o)\lambda^2 d\lambda d\o.
\end{equation}
We assume furthermore that $b(x,\o)$ satisfies:
\begin{equation}\label{hyp1}
\norm{b}_{H^{3/2}(\R^3)}+\norm{\nabla b}_{\l{\infty}{2}}+\norm{\nabb\nabla b}_{\le{2}}+\norm{\po b}_{H^{1/2+\a}(\R^3)}\leq D,
\end{equation}
for some constant $D>0$ and $\a>0$. Then, $B$ is bounded on $L^2$ and satisfies the estimate:
\begin{equation}\label{l2e}
\norm{Bf}_{L^2(\R^3)}\lesssim D\norm{f}_{\le{2}}.
\end{equation}
\end{theorem}

\begin{remark}
We do not claim that Theorem \ref{th3} is an improvement compared to the vast literature on boundedness on $L^2$ for pseudodifferential operators. Its purpose is to give a warm up for the proof of Theorem \ref{th1}, i.e. boundedness on $L^2$  for Fourier integral operators on a 3 dimensional Riemannian manifold $\s$. 
\end{remark}

\begin{remark}
The assumptions \eqref{hyp1} hold for any $\o\in\S$ with the constant in the right-hand side being independent of $\o$. Thus, one may take the supremum in $\o$ everywhere. To ease the notations, we do not explicitly write down this supremum. 
\end{remark}

\begin{remark}
In the Euclidean setting, the derivative $\nabb$ simply refers to derivatives in directions orthogonal to $\o$. Also, the space $\l{\infty}{2}$ is defined with respect to $u(x,\o)=\xo$ and the level surfaces of $u$ are now planes $\p=\{x/\,\xo=u\}$.
\end{remark}

\begin{remark}\label{rm:hyp1}
The assumptions on the symbol $b(x,\o)$ in Theorem \ref{th3} are slightly different from the ones in Theorem \ref{th1}. In particular, we do not assume that $b\in\le{\infty}$ since this is a consequence of  the assumption \eqref{hyp1} and Sobolev embeddings in dimension 3. Also, the assumption \eqref{cordecfr1s} follows from assumption \eqref{hyp1}. Indeed, let $\Delta_j$ denote the usual Littlewood Paley projections in $\R^3$ which localizes at frequencies of size $2^j$. We may decompose $\nabla b=b^j_1+b^j_2$ with $b^j_1=\Delta_{> j}\nabla b$ and $b^2_j=\Delta_{\leq j}\nabla b$ and we obtain \eqref{cordecfr1s} by using $\norm{b}_{\h{3/2}}\leq D$. Finally, \eqref{hyp1} only assumes $\norm{\po b}_{H^{1/2+\a}(\R^3)}\leq D$ while \eqref{threomega1s} assumes essentially that $\norm{\po b}_{H^{1}(\R^3)}$ is bounded. We may actually relax \eqref{threomega1s} by replacing it with the analog of $\norm{\po b}_{H^{1/2+\a}(\R^3)}\leq D$. However, this would require to discuss fractional Sobolev spaces on $\s$ and would complicate the exposition.
\end{remark}

The rest of the paper is as follows. In section \ref{sec:proofth3}, we prove Theorem \ref{th3}. In section \ref{sec:th1}, we prove Theorem \ref{th1}. Finally, we prove Theorem \ref{th2} in section \ref{sec:th2}. 

\section{Proof of Theorem \ref{th3}}\label{sec:proofth3}

While the conclusion of Theorem \ref{th3} follows from Theorem \ref{th1} in the case $(\s,g)=(\R^3,\de)$ where $\de$ is the euclidean metric, and 
$u(x,\o)=\xo$, it will be instructive to perform the proof first in this simple case of a pseudodifferential operator on $\R^3$. This will clarify the main ideas, before turning to the proof of Theorem \ref{th1} for Fourier integral operators on a 3 dimensional Riemannian manifold $\s$ in section \ref{sec:th1}.

\subsection{The basic computation}

Since the Fourier transform is an isomorphism of $\le{2}$, we may remove the Fourier transform in the definition \eqref{opd} of $B$ in order to ease the notations: 
\begin{equation}\label{opd1}
Bf(x)=\int_{\S}\int_{0}^{+\infty}e^{i\lambda \xo}b(x,\o)f(\lambda\o)\lambda^2 d\lambda d\o.
\end{equation}
We start the proof of Theorem \ref{th3} with the following instructive computation:
\begin{equation}\label{b5}
\begin{array}{ll}
\ds\norm{Bf}_{\le{2}}& \ds\leq\int_{\S}\normm{b(x,\o)\int_{0}^{+\infty}e^{i\lambda \xo}f(\lambda\o)\lambda^2 d\lambda}_{\le{2}}d\o\\
& \ds\leq\int_{\S}\norm{b(x,\o)}_{\l{\infty}{2}}\normm{\int_{0}^{+\infty}e^{i\lambda \xo}f(\lambda\o)\lambda^2 d\lambda}_{L^2_{\xo}}d\o\\
& \ds\leq D\norm{\lambda f}_{L^2(\R^3)},
\end{array}
\end{equation}
where we have used Plancherel with respect to $\lambda$, Cauchy-Schwarz with respect to $\o$ 
and \eqref{hyp1} to bound $\norm{b}_{\l{\infty}{2}}$ (note that the space $\l{\infty}{2}$ is defined with respect to $u(x,\o)=\xo$ and the level surfaces of $u$ are now planes $\p=\{x/\,\xo=u\}$). \eqref{b5} misses the conclusion \eqref{l2e} of Theorem \ref{th3} by a power of $\lambda$. Now, assume for a moment that we may replace a power of $\lambda$ by a derivative on $b(x,\o)$. Then, the same computation yields:
\begin{equation}\label{b6}
\begin{array}{ll}
& \ds\normm{\int_{\S}\int_{0}^{+\infty}\nabla b(x,\o) e^{i\lambda \xo}f(\lambda\o)\lambda d\lambda d\o}_{\le{2}}\\
\ds\leq & \ds\int_{\S}\norm{\nabla b(x,\o)}_{\l{\infty}{2}}\normm{\int_{0}^{+\infty}e^{i\lambda\xo}f(\lambda\o)\lambda^2 d\lambda}_{L^2_{\xo}}d\o\\
\ds\leq & \ds D\norm{f}_{L^2(\R^3)},
\end{array}
\end{equation}
which is \eqref{l2e}. This suggests a strategy which consists in making integrations by parts to trade powers of $\lambda$ against derivatives of the symbol $b(x,\o)$. 

\subsection{Structure of the proof of Theorem \ref{th3}}

The proof of Theorem \ref{th3} proceeds in three steps. We first localize in frequencies of size $\la\sim 2^j$. We then localize the angle $\o$ in patches on the sphere $\S$ of diameter $2^{-j/2}$. Finally, we estimate the diagonal terms.

\subsubsection{Step 1: decomposition in frequency}
 
For the first step, we introduce $\varphi$  and $\psi$ two smooth compactly supported functions on $\R$ such that: 
\begin{equation}\label{b7}
\varphi(\lambda)+\sum_{j\geq 0}\psi(2^{-j}\lambda)=1\textrm{ for all }\lambda\in\R.
\end{equation}
We use \eqref{b7} to decompose $Bf$ as follows:
\begin{equation}\label{b8}
Bf(x)=\sum_{j\geq -1}B_jf(x),
\end{equation}
where for $j\geq 0$:
\begin{equation}\label{b9}
B_jf(x)=\int_{\S}\int_{0}^{+\infty}e^{i\lambda\xo}b(x,\o)\psi(2^{-j}\lambda)f(\lambda\o)\lambda^2 d\lambda d\o,
\end{equation}
and 
\begin{equation}\label{b10}
B_{-1}f(x)=\int_{\S}\int_{0}^{+\infty}e^{i\lambda\xo}b(x,\o)\varphi(\lambda)f(\lambda\o)\lambda^2 d\lambda d\o.
\end{equation}
This decomposition is classical and is known as the first dyadic decomposition (see \cite{stein}). The goal of this first step is to prove the following proposition:
\begin{proposition}\label{orthofreq}
The decomposition \eqref{b8} satisfies an almost orthogonality property:
\begin{equation}\label{orthofreq1}
\norm{Bf}_{\le{2}}^2\lesssim\sum_{j\geq -1}\norm{B_jf}_{\le{2}}^2+D^2\norm{f}^2_{\le{2}}.
\end{equation}
\end{proposition}
The proof of Proposition \ref{orthofreq} is postponed to section \ref{sec:orthofreq}. 

\subsubsection{Step 2: decomposition in angle}

Proposition \ref{orthofreq} allows us to  estimate $\norm{B_jf}_{\le{2}}$ instead of 
$\norm{Bf}_{\le{2}}$. The analog of computation \eqref{b5} for $\norm{B_jf}_{\le{2}}$ yields:
\begin{equation}\label{b5bis}
\begin{array}{l}
\ds\norm{B_jf}_{\le{2}}\leq D\norm{\lambda\psi(2^{-j}\la) f}_{L^2(\R^3)}\lesssim D2^j\norm{\psi(2^{-j}\la) f}_{L^2(\R^3)},
\end{array}
\end{equation}
which misses the wanted estimate by a power of $2^j$. We thus need to perform a second dyadic decomposition (see \cite{stein}). We introduce a smooth partition of unity on the sphere $\S$:
\begin{equation}\label{b14}
\sum_{\nu\in\Gamma}\eta^\nu_j(\o)=1\textrm{ for all }\o\in\S, 
\end{equation}
where the support of $\eta^\nu_j$ is a patch on $\S$ of diameter $\sim 2^{-j/2}$. We use \eqref{b14} to decompose $B_jf$ as follows:
\begin{equation}\label{b15}
B_jf(x)=\sum_{\nu\in\Gamma}B^\nu_jf(x),
\end{equation}
where:
\begin{equation}\label{b16}
B^\nu_jf(x)=\int_{\S}\int_{0}^{+\infty}e^{i\lambda\xo}b(x,\o)\psi(2^{-j}\lambda)\eta^\nu_j(\o)f(\lambda\o)\lambda^2 d\lambda d\o.
\end{equation}
We also define:
\begin{equation}\label{decf}
\begin{array}{l}
\ds\ga_{-1}=\norm{\varphi(\la)f}_{\le{2}},\, \ga_j=\norm{\psi(2^{-j}\la)f}_{\le{2}},\,j\geq 0, \\
\ds\ga^\nu_j=\norm{\psi(2^{-j}\la)\eta^\nu_j(\o)f}_{\le{2}},\,j\geq 0,\,\nu\in\Gamma, 
\end{array}
\end{equation}
which satisfy:
\begin{equation}\label{decf1}
\norm{f}_{\le{2}}^2=\sum_{j\geq -1}\ga_j^2=\sum_{j\geq -1}\sum_{\nu\in\Gamma}(\ga^\nu_j)^2.
\end{equation}
The goal of this second step is to prove the following proposition:
\begin{proposition}\label{orthoangle}
The decomposition \eqref{b15} satisfies an almost orthogonality property:
\begin{equation}\label{orthoangle1}
\norm{B_jf}_{\le{2}}^2\lesssim\sum_{\nu\in\Gamma}\norm{B^\nu_jf}_{\le{2}}^2+D^2\ga_j^2.
\end{equation}
\end{proposition}
The proof of Proposition \ref{orthoangle} is postponed to section \ref{sec:orthoangle}. 

\subsubsection{Step 3: control of the diagonal term}

Proposition \ref{orthoangle} allows us to  estimate $\norm{B^\nu_jf}_{\le{2}}$ instead of $\norm{B_jf}_{\le{2}}$. The analog of computation \eqref{b5} for $\norm{B^\nu_jf}_{\le{2}}$ yields:
\begin{equation}\label{b5ter}
\hspace{-0.15cm}\begin{array}{ll}
\ds\norm{B_j^\nu f}_{\le{2}} & \ds\leq\int_{\S}\norm{b(x,\o)}_{\l{\infty}{2}}\normm{\int_{0}^{+\infty}e^{i\lambda \xo}\psi(2^{-j}\la)\eta^\nu_j(\o)f(\lambda\o)\lambda^2 d\lambda}_{L^2_{\xo}}d\o\\
& \ds\leq D\sqrt{\textrm{vol}(\textrm{supp}(\eta^\nu_j))}\norm{\lambda\psi(2^{-j}\la)\eta^\nu_j(\o)f}_{L^2(\R^3)}\\
& \ds\lesssim D2^{j/2}\gamma^\nu_j,
\end{array}
\end{equation}
where the term $\sqrt{\textrm{vol}(\textrm{supp}(\eta^\nu_j))}$ comes from the fact that we apply Cauchy-Schwarz in $\o$. Note that we have used in \eqref{b5ter} the fact that the support of $\eta^\nu_j$ is 2 dimensional and has diameter $2^{-j/2}$ so that:
\begin{equation}\label{suppangle}
\sqrt{\textrm{vol}(\textrm{supp}(\eta^\nu_j))}\lesssim 2^{-j/2}.
\end{equation}
Now, \eqref{b5ter} still misses the wanted estimate by a power of $2^{j/2}$. Nevertheless, taking advantage of the regularity of $\po b$ given by \eqref{hyp1}, we are able to estimate the diagonal term:
\begin{proposition}\label{diagonal}
The diagonal term $B^\nu_jf$ satisfies the following estimate:
\begin{equation}\label{diagonal1}
\norm{B^\nu_jf}_{\le{2}}\lesssim D\ga^\nu_j.
\end{equation}
\end{proposition}
The proof of Proposition \ref{diagonal} is postponed to section \ref{sec:diagonal}. 

\subsubsection{Proof of Theorem \ref{th3}}

Proposition \ref{orthofreq}, \ref{orthoangle} and \ref{diagonal} 
immediately yield the proof of Theorem \ref{th3}. Indeed, \eqref{orthofreq1}, \eqref{decf1}, \eqref{orthoangle1} and \eqref{diagonal1} imply:
\begin{equation}\label{cclth3}
\begin{array}{ll}
\ds\norm{Bf}_{\le{2}}^2 & \ds\lesssim\sum_{j\geq -1}\norm{B_jf}_{\le{2}}^2+D^2\norm{f}^2_{\le{2}}\\
& \ds\lesssim\sum_{j\geq -1}\sum_{\nu\in\Gamma}\norm{B^\nu_jf}_{\le{2}}^2+D^2\sum_{j\geq -1}\gamma_j^2+D^2\norm{f}^2_{\le{2}}\\
& \ds\lesssim D^2\sum_{j\geq -1}\sum_{\nu\in\Gamma}(\gamma_j^\nu)^2+D^2\sum_{j\geq -1}\gamma_j^2+D^2\norm{f}^2_{\le{2}}\\
& \ds\lesssim D^2\norm{f}^2_{\le{2}},
\end{array}
\end{equation}
which is the conclusion of Theorem \ref{th3}. \QED

The remainder of section \ref{sec:proofth3} is dedicated to the proof of Proposition \ref{orthofreq}, \ref{orthoangle} and \ref{diagonal}.

\subsection{Proof of Proposition \ref{orthofreq} (almost orthogonality in frequency)}\label{sec:orthofreq}

We have to prove \eqref{orthofreq1}:
\begin{equation}\label{of1}
\norm{Bf}_{\le{2}}^2\lesssim\sum_{j\geq -1}\norm{B_jf}_{\le{2}}^2+D^2\norm{f}^2_{\le{2}}.
\end{equation}
This will result from the following inequality using Shur's Lemma:
\begin{equation}\label{of2}
\left|\int_{\R^3}B_jf(x)\overline{B_kf(x)}dx\right| \lesssim D^22^{-\frac{|j-k|}{2}}\gamma_j\ga_k\textrm{ for }|j-k|> 2.
\end{equation}

\subsubsection{A first integration by parts} 

From now on, we focus on proving \eqref{of2}. We may assume $j\geq k+3$. We have:
\begin{equation}\label{of3}
\begin{array}{ll}
\ds\int_{\R^3}B_jf(x)\overline{B_kf(x)}dx= & \ds\int_{\S}\int_{0}^{+\infty}\int_{\S}\int_{0}^{+\infty}\left(\int_{\R^3}e^{i\lambda\xo-i\la'\xo'}b(x,\o)\overline{b(x,\o')}dx\right)\\
& \ds\times\psi(2^{-j}\lambda)f(\lambda\o)\lambda^2 \psi(2^{-k}\lambda')\overline{f(\lambda'\o')}(\lambda')^2 d\lambda d\o d\lambda' d\o'.
\end{array}
\end{equation}
We integrate by parts with respect to $\partial_{\xo}$ in $\int_{\R^3}e^{i\lambda\xo-i\la'\xo'}b(x,\o)\overline{b(x,\o')}dx$ using the fact that:
\begin{equation}\label{of4}
e^{i\lambda\xo-i\la'\xo'}=-\frac{i}{\la-\la'\o\cdot\o'}\partial_{\xo}(e^{i\lambda\xo-i\la'\xo'}).
\end{equation}
We obtain:
\begin{equation}\label{of5}
\begin{array}{ll}
\ds\int_{\R^3}e^{i\lambda\xo-i\la'\xo'}b(x,\o)\overline{b(x,\o')}dx= & \ds i\int_{\R^3}e^{i\lambda\xo-i\la'\xo'}\frac{\partial_{\xo}b(x,\o)\overline{b(x,\o')}}{\la-\la'\o\cdot\o'}dx\\
& \ds +i\int_{\R^3}e^{i\lambda\xo-i\la'\xo'}\frac{b(x,\o)\partial_{\xo}\overline{b(x,\o')}}{\la-\la'\o\cdot\o'}dx.
\end{array}
\end{equation}
Since $|\la'\o\cdot\o'|<\la$, we may expand the fractions in \eqref{of5}:
\begin{equation}\label{of6}
\frac{1}{\la-\la'\o\cdot\o'}=\sum_{p\geq 0}\frac{1}{\la}\left(\frac{\la'\o\cdot\o'}{\la}\right)^p.
\end{equation}
For $p\in\mathbb{Z}$, We introduce the notation $F_{j,p}(\xo)$:
\begin{equation}\label{of7}
F_{j,p}(\xo)=\int_{0}^{+\infty}e^{i\lambda\xo}\psi(2^{-j}\lambda)f(\lambda\o)(2^{-j}\la)^{p}\lambda^2 d\lambda.
\end{equation}
Together with \eqref{of3}, \eqref{of5} and \eqref{of6}, this implies:
\begin{equation}\label{of8}
\int_{\R^3}B_jf(x)\overline{B_kf(x)}dx=\sum_{p\geq 0}A^1_p+\sum_{p\geq 0}A^2_p,
\end{equation}
where $A^1_p$ and $A^2_p$ are given by:
\begin{equation}\label{of9}
\begin{array}{l}
\ds A^1_p=\ds 2^{-j-p(j-k)}\\
\ds\times\int_{\R^3}\left(\int_{\S}\partial_{\xo}b(x,\o)\o^pF_{j,-p-1}(\xo)d\o\right)\cdot\overline{\left(\int_{\S}b(x,\o'){\o'}^pF_{k,p}(\xo')d\o'\right)}dx,
\end{array}
\end{equation}
and
\begin{equation}\label{of10}
\begin{array}{l}
\ds A^2_p=\ds 2^{-j-p(j-k)}\\
\ds\times\int_{\R^3}\left(\int_{\S}b(x,\o)\o^{p+1}F_{j,-p-1}(\xo)d\o\right)\cdot\overline{\left(\int_{\S}\nabla b(x,\o'){\o'}^pF_{k,p}(\xo')d\o'\right)}dx.
\end{array}
\end{equation}

\begin{remark}\label{rmksep}
The expansion \eqref{of6} allows us to rewrite $\int_{\R^3}B_jf(x)\overline{B_kf(x)}dx$ in the form \eqref{of8}, i.e. as a sum of terms $A^1_p$, $A^2_p$. The key point is that in each of these terms - according to \eqref{of9} and \eqref{of10} - one may separate the terms depending of $(\la,\o)$ from the terms depending on $(\la',\o')$. 
\end{remark}

\subsubsection{Estimates for $A^1_p$ and $A^2_p$} The term containing one derivative of $b$ in \eqref{of9} may be estimated using the basic computation \eqref{b5}:
\begin{equation}\label{of11}
\begin{array}{ll}
& \ds\normm{\int_{\S}\partial_{\xo}b(x,\o)\o^pF_{j,-p-1}(\xo)d\o}_{\le{2}}\\
\ds\leq &\ds\int_{\S}\norm{\partial_{\xo}b(x,\o)\o^p}_{\l{\infty}{2}}\norm{F_{j,-p-1}(\xo)}_{L^2_{\xo}}d\o\\
\ds\leq &\ds\norm{\nabla b}_{\l{\infty}{2}}\norm{\psi(2^{-j}\lambda)f(\lambda\o)(2^{-j}\la)^{-p-1}\la}_{\le{2}}\\
\leq &\ds D2^{p+1+j}\ga_j,
\end{array}
\end{equation}
where we have used the assumption \eqref{hyp1} on $b$ and the fact that $(2^{-j}\la)^{-1}\leq 2$ on the support of $\psi(2^{-j}\la)$. In the same way, the term containing one derivative of $b$ in \eqref{of10} may be estimated by:
\begin{equation}\label{of12}
\begin{array}{ll}
&\ds\normm{\int_{\S}\nabla b(x,\o'){\o'}^pF_{k,p}(\xo')d\o'}_{\le{2}}\\
\ds\leq &\ds\int_{\S}\norm{\nabla b(x,\o'){\o'}^p}_{\l{\infty}{2}}\norm{F_{k,p}(\xo')}_{L^2_{\xo'}}d\o'\\
\ds\leq &\ds\norm{\nabla b}_{\l{\infty}{2}}\norm{\psi(2^{-k}\lambda')f(\lambda'\o')(2^{-k}\la')^{p}\la'}_{\le{2}}\\
\ds \leq &\ds D2^{p+k}\ga_k,
\end{array}
\end{equation}
where we have used the assumption \eqref{hyp1} on $b$ and the fact that $(2^{-k}\la')\leq 2$ on the support of $\psi(2^{-k}\la')$.

Note that Proposition \ref{orthoangle} together with Proposition \ref{diagonal} yields the estimate:
\begin{equation}\label{of13}
\norm{B_jf}_{\le{2}}\lesssim D\ga_j,
\end{equation}
for any symbol $b$ satisfying the assumptions \eqref{hyp1}. Now, the term containing no derivative of $b$ in \eqref{of9} has a symbol given by $b(x,\o'){\o'}^p$ which satisfies the assumptions \eqref{hyp1} since $b$ does. Applying \eqref{of13}, we obtain: 
\begin{equation}\label{of14}
\begin{array}{ll}
& \ds\normm{\int_{\S}b(x,\o'){\o'}^pF_{k,p}(\xo')d\o'}_{\le{2}}\\
\ds\lesssim &\ds D\norm{\psi(2^{-k}\lambda')f(\lambda'\o')(2^{-k}\la')^{p}}_{\le{2}}\\
\leq & \ds D2^{p}\ga_k.
\end{array}
\end{equation}
In the same way, the term containing no derivative of $b$ in \eqref{of10} has a symbol given by $b(x,\o){\o}^{p+1}$ which satisfies the assumptions \eqref{hyp1} since $b$ does. Applying again \eqref{of13}, we obtain: 
\begin{equation}\label{of15}
\begin{array}{ll}
&\ds\normm{\int_{\S}b(x,\o){\o}^{p+1}F_{j,-p-1}(\xo)d\o}_{\le{2}}\\
\ds\lesssim & \ds D\norm{\psi(2^{-j}\lambda)f(\lambda\o)(2^{-j}\la)^{-p-1}}_{\le{2}}\\
\lesssim & \ds D2^{p+1}\ga_j.
\end{array}
\end{equation}

Finally, the definition of $A_p^1$ \eqref{of9} and the estimates \eqref{of11} and \eqref{of14} yield:
\begin{equation}\label{of16}
|A^1_p|\lesssim D2^{2p-p(j-k)}\ga_j\ga_k,\,\forall p\geq 0.
\end{equation}
Similarly, the definition of $A_p^2$ \eqref{of10} and the estimates \eqref{of12} and \eqref{of15} yield:
\begin{equation}\label{of17}
|A^2_p|\lesssim D2^{2p-(p+1)(j-k)}\ga_j\ga_k,\,\forall p\geq 0.
\end{equation}
\eqref{of16} and \eqref{of17} imply:
\begin{equation}\label{of18}
\sum_{p\geq 1}|A^1_p|+\sum_{p\geq 0}|A^2_p|\lesssim D2^{-(j-k)}\left(\sum_{p\geq 0}2^{-p(j-k-2)}\right)\ga_j\ga_k\lesssim D2^{-(j-k)}\ga_j\ga_k,
\end{equation}
where we have used the assumption $j-k-2>0$. \eqref{of8} and \eqref{of18} will yield \eqref{of2} provided we obtain a similar estimate for $A^1_0$. Now, the estimate of $A^1_0$ provided by \eqref{of16} is not sufficient since it does not contain any decay in $j-k$. We will need to perform a second integration by parts for this term.

\subsubsection{A more precise estimate for $A^1_0$} From \eqref{of9} with $p=0$, we have:
\begin{equation}\label{of19}
\ds A^1_0=\ds 2^{-j}\int_{\R^3}\left(\int_{\S}\partial_{\xo}b(x,\o)F_{j,0}(\xo)d\o\right)\overline{B_k(x)}dx.
\end{equation}
Since $b(x,\o)$ is assumed to be in $\h{3/2}$, we may only make one half integration by parts. To this end, we decompose $\partial_{\xo}b$ as in Remark \ref{rm:hyp1}. Let $\Delta_j$ denote the usual Littlewood Paley projections in $\R^3$ which localizes at frequencies of size $2^j$. We decompose $\partial_{\xo}b=b^j_1+b^j_2$ with $b^j_1=\Delta_{> j}\partial_{\xo}b$ and $b^j_2=\Delta_{\leq j}\partial_{\xo}b$ and we obtain 
\begin{equation}\label{of20}
\norm{b^1_j}_{\le{2}}\lesssim D2^{-\frac{j}{2}}\textrm{ and }\norm{\nabla b^j_2}_{\le{2}}\lesssim D2^{\frac{j}{2}}
\end{equation}
by using $\norm{b}_{\h{3/2}}\leq D$. In turn, this yields a decomposition for $A^1_0$:
\begin{equation}\label{of21}
A^1_0=A^1_{0,1}+A^1_{0,2}
\end{equation}
where:
\begin{equation}\label{of22}
\begin{array}{l}
\ds A^1_{0,1}=2^{-j}\int_{\R^3}\left(\int_{\S}b^j_1(x,\o)F_{j,0}(\xo)d\o\right)\overline{B_k(x)}dx,\\[3mm]
\ds A^1_{0,2}=2^{-j}\int_{\R^3}\left(\int_{\S}b^j_2(x,\o)F_{j,0}(\xo)d\o\right)\overline{B_k(x)}dx.
\end{array}
\end{equation}

We first estimate $A^1_{0,1}$. We have: 
\begin{equation}\label{of23}
\begin{array}{ll}
\ds |A^1_{0,1}| & \ds\leq 2^{-j}\int_{\S}\left|\int_{\R^3}b^j_1(x,\o)F_{j,0}(\xo)\overline{B_k(x)}dx\right|d\o\\
& \ds\leq 2^{-j}\int_{\S}\norm{b^j_1(.,\o)}_{\le{2}}\norm{F_{j,0}}_{L^2_{\xo}}\norm{B_k}_{\l{\infty}{2}}d\o\\
& \ds\lesssim D2^{-\frac{3j}{2}}\int_{\S}\norm{F_{j,0}}_{L^2_{\xo}}\norm{B_k}_{\l{\infty}{2}}d\o,
\end{array}
\end{equation}
where we have used \eqref{of20} in the last inequality. Plancherel yields:
\begin{equation}\label{of24}
\norm{F_{j,0}}_{L^2_{\o,\xo}}\leq\norm{\psi(2^{-j}\la)f(\la\o)\la}_{\le{2}}\lesssim 2^j\ga_j.
\end{equation}
In view of \eqref{of23}, we also need to estimate $\norm{B_k}_{\l{\infty}{2}}$. We have:
\begin{equation}\label{of25}
\norm{B_k}_{\l{\infty}{2}}\lesssim\norm{B_k}^{\frac{1}{2}}_{\h{1}}\norm{B_k}^{\frac{1}{2}}_{\le{2}}\lesssim D^{\frac{1}{2}}\ga_k^{\frac{1}{2}}\norm{B_k}^{\frac{1}{2}}_{\h{1}},
\end{equation}
where we have used a standard trace Theorem for the first inequality, and \eqref{of13} for the second inequality. We still need to estimate $\norm{\nabla B_k}_{\le{2}}$. We have:
\begin{equation}\label{of26}
\begin{array}{ll}
\ds\nabla B_k(x)= & \ds\int_{\S}\int_0^{+\infty} e^{i\la\xo}\nabla b(x,\o)\psi(2^{-k}\la)f(\la\o)\la^2d\la d\o\\
& \ds +i2^k\int_{\S}\int_0^{+\infty} e^{i\la\xo}\o b(x,\o)\psi(2^{-k}\la)(2^{-k}\la)f(\la\o)\la^2d\la d\o.
\end{array}
\end{equation}
Using the basic computation \eqref{b5} for the first term together with the fact that $\nabla b\in\l{\infty}{2}$, and \eqref{of13} for the second term together with the fact that $\o b(x,\o)$ satisfies the assumption \eqref{hyp1}, we obtain:
\begin{equation}\label{of27}
\norm{\nabla B_k}_{\le{2}}\lesssim D2^k\ga_k.
\end{equation}
Finally, \eqref{of23}, \eqref{of24}, \eqref{of25} and \eqref{of27} yield:
\begin{equation}\label{of28}
|A^1_{0,1}|\lesssim D2^{-\frac{j-k}{2}}\ga_j\ga_k.
\end{equation}

\subsubsection{A second integration by parts} We now estimate the term $A^1_{0,2}$ defined in \eqref{of22}. We perform a second integration by parts relying again on \eqref{of4}. We obtain: 
\begin{equation}\label{of29}
\hspace{-0.1cm}\begin{array}{l}
\ds A^1_{0,2}=2^{-2j}\int_{\R^3}\left(\int_{\S}\partial_{\xo}b^j_2(x,\o)F_{j,0}(\xo)d\o\right)\overline{B_k(x)}dx\\
\ds +2^{-2j}\int_{\R^3}\left(\int_{\S}b^j_2(x,\o)\o F_{j,0}(\xo)d\o\right)\cdot\overline{\left(\int_{\S}\nabla b(x,\o')F_{k,0}(\xo')d\o'\right)}dx+\cdots,
\end{array}
\end{equation}
where we only mention the first term generated by the expansion \eqref{of6}. In fact, the other terms are estimated in the same way and generate more decay in $j-k$ similarly to the estimates \eqref{of16} \eqref{of17}.

The first term in the right-hand side of \eqref{of29} has the same form than $A^1_{0,1}$ defined in \eqref{of22} where $b^j_1$ is replaced by $2^{-j}\partial_{\xo}b^j_2$. By \eqref{of20}, $2^{-j}\partial_{\xo}b^j_2$ satisfies:
$$\norm{2^{-j}\partial_{\xo}b^j_2}_{\le{2}}\lesssim D2^{-\frac{j}{2}}.$$
Since $b^1_j$ and $2^{-j}\partial_{\xo}b^j_2$ satisfy the same estimate, we obtain the analog of \eqref{of28} for the first term in the right-hand side of \eqref{of29}:
\begin{equation}\label{of30}
\ds \left|2^{-2j}\int_{\R^3}\left(\int_{\S}\partial_{\xo}b^j_2(x,\o)F_{j,0}(\xo)d\o\right)\overline{B_k(x)}dx\right|\lesssim D2^{-\frac{j-k}{2}}\ga_j\ga_k.
\end{equation}

We now estimate the second term in the right-hand side of \eqref{of29}. Recall that $b^j_2=\Delta_{\leq j}\partial_{\xo}b$ so that together with the assumption \eqref{hyp1}, we have:
\begin{equation}\label{of31}
\norm{b^j_2}_{\l{\infty}{2}}\lesssim D.
\end{equation}
We estimate the scorn term in the right-hand side of \eqref{of29} using the assumption \eqref{hyp1}, the basic computation \eqref{b5} and \eqref{of31}: 
\begin{equation}\label{of32}
\begin{array}{ll}
&\ds\left|2^{-2j}\int_{\R^3}\left(\int_{\S}b^j_2(x,\o)\o F_{j,0}(\xo)d\o\right)\cdot\overline{\left(\int_{\S}\nabla b(x,\o')F_{k,0}(\xo')d\o'\right)}dx\right|\\
\ds\leq &\ds 2^{-2j}\normm{\int_{\S}b^j_2(x,\o)\o F_{j,0}(\xo)d\o}_{\le{2}}\normm{\int_{\S}\nabla b(x,\o')F_{k,0}(\xo')d\o'}_{\le{2}}\\
\ds\leq &\ds 2^{-2j}\left(\int_{\S}\norm{b^j_2(.,\o)\o}_{\l{\infty}{2}}\norm{F_{j,0}}_{L^2_{\xo}}d\o\right)\\
&\ds\times\left(\int_{\S}\norm{\nabla b(.,\o)}_{\l{\infty}{2}}\norm{F_{k,0}}_{L^2_{\xo}}d\o\right)\\
\lesssim &\ds D^22^{-(j-k)}\ga_j\ga_k.
\end{array}
\end{equation}
Finally, \eqref{of29}, \eqref{of30} and \eqref{of32} imply:
\begin{equation}\label{of33}
|A^1_{0,2}|\lesssim D^22^{-\frac{j-k}{2}}\ga_j\ga_k.
\end{equation}

\subsubsection{End of the proof of Proposition \ref{orthofreq}} Since $A^1_0=A^1_1+A^1_2$, the estimate \eqref{of28} of $A^1_{0,1}$ and the estimate \eqref{of33} of $A^1_{0,2}$ yield:
\begin{equation}\label{of34}
|A^1_{0}|\lesssim D^22^{-\frac{j-k}{2}}\ga_j\ga_k.
\end{equation}
Together with \eqref{of8} and \eqref{of18}, this implies:
\begin{equation}\label{of35}
\left|\int_{\R^3}B_jf(x)\overline{B_kf(x)}dx\right|\lesssim D^22^{-\frac{|j-k|}{2}}\ga_j\ga_k\textrm{ for }|j-k|>2.
\end{equation}
Finally, \eqref{of35} together with Shur's Lemma yields:
\begin{equation}\label{of36}
\norm{Bf}_{\le{2}}^2\lesssim\sum_{j\geq -1}\norm{B_jf}_{\le{2}}^2+D^2\norm{f}^2_{\le{2}}.
\end{equation}
This concludes the proof of Proposition \ref{orthofreq}. \QED

\subsection{Proof of Proposition \ref{orthoangle} (almost orthogonality in angle)}\label{sec:orthoangle}

We have to prove \eqref{orthoangle1}:
\begin{equation}\label{oa1}
\norm{B_jf}_{\le{2}}^2\lesssim\sum_{\nu\in\Gamma}\norm{B^\nu_jf}_{\le{2}}^2+D^2\ga_j^2.
\end{equation}
This will result from the following inequality:
\begin{equation}\label{oa2}
\left|\int_{\R^3}B^\nu_jf(x)\overline{B^{\nu'}f(x)}dx\right|\lesssim \frac{D^2\ga_j^\nu\ga_j^{\nu'}}{2^{j\a/2}(2^{j/2}|\nu-\nu'|)^{2-\a}},\,|\nu-\nu'|\neq 0,
\end{equation}
where $\a>0$. Indeed, since $\S$ is 2 dimensional and $1\leq 2^{j/2}|\nu-\nu'|\leq 2^{j/2}$ for $\nu, \nu'\in\Gamma$ and $\nu\neq\nu'$, we have:
\begin{equation}\label{oa3}
\sup_{\nu}\sum_{\nu'} \frac{1}{2^{j\a/2}(2^{j/2}|\nu-\nu'|)^{2-\a}}\leq C_\a<+\infty\,\forall\a>0.
\end{equation}
Thus, \eqref{oa2} and \eqref{oa3} together with Shur's Lemma imply \eqref{oa1}.

\subsubsection{A second decomposition in frequency} From now on, we focus on proving \eqref{oa2}. Integrating by parts twice in $\int_{\R^3}B^\nu_jf(x)\overline{B^{\nu'}_jf(x)}dx$ would ultimately yield:
\begin{equation}\label{oa4}
\left|\int_{\R^3}B^\nu_jf(x)\overline{B^{\nu'}_jf(x)}dx\right|\lesssim \frac{D^2\ga_j^\nu\ga_j^{\nu'}}{(2^{j/2}|\nu-\nu'|)^{2}},\,|\nu-\nu'|\neq 0.
\end{equation}
This corresponds to the case $\a=0$ in \eqref{oa3} and yields to a log-loss since we have:
\begin{equation}\label{oa5}
\sup_{\nu}\sum_{\nu'} \frac{1}{(2^{j/2}|\nu-\nu'|)^{2}}\sim j.
\end{equation}
To avoid this log-loss, we do a second decomposition in frequency. $\la$ belongs to the interval $[2^{j-1},2^{j+1}]$ which we decompose in intervals $I_k$:
\begin{equation}\label{oa6}
[2^{j-1},2^{j+1}]=\bigcup_{1\leq k\leq |\nu-\nu'|^{-\a}}I_k\textrm{ where }\textrm{diam}(I_k)\sim 2^j|\nu-\nu'|^\a.
\end{equation}
Let $\phi_k$ a partition of unity of the interval $[2^{j-1},2^{j+1}]$ associated to the $I_k$'s. We decompose $B^\nu_jf$ as follows:
\begin{equation}\label{oa7}
B^\nu_jf(x)=\sum_{1\leq k\leq |\nu-\nu'|^{-\a}}B^{\nu,k}_jf(x),
\end{equation}
where:
\begin{equation}\label{oa8}
B^{\nu,k}_jf(x)=\int_{\S}\int_{0}^{+\infty}e^{i\lambda\xo}b(x,\o)\psi(2^{-j}\lambda)\eta^\nu_j(\o)\phi_k(\la)f(\lambda\o)\lambda^2 d\lambda d\o.
\end{equation}
We also define:
\begin{equation}\label{oa9}
\begin{array}{l}
\ds\ga^{\nu,k}_j=\norm{\psi(2^{-j}\la)\eta^\nu_j(\o)\phi_k(\la)f}_{\le{2}},\,j\geq 0,\,\nu\in\Gamma,1\leq k\leq |\nu-\nu'|^{-\a}, 
\end{array}
\end{equation}
which satisfy:
\begin{equation}\label{oa10}
(\gamma^{\nu}_j)^2=\sum_{1\leq k\leq |\nu-\nu'|^{-\a}}(\ga^{\nu,k}_j)^2.
\end{equation}

\subsubsection{The two key estimates}\label{sec:keyest} 

We will prove the following two estimates:
\begin{equation}\label{oa11}
\left|\int_{\R^3}B^{\nu,k}_jf(x)\overline{B^{\nu',k}_jf(x)}dx\right|\lesssim \frac{D^2\ga_j^{\nu,k}\ga_j^{\nu',k}}{2^{j\a/2}(2^{j/2}|\nu-\nu'|)^{2-\a}},\,|\nu-\nu'|\neq 0,\,1\leq k\leq |\nu-\nu'|^{-\a},
\end{equation}
and
\begin{equation}\label{oa12}
\begin{array}{r}
\ds\left|\int_{\R^3}B^{\nu,k}_jf(x)\overline{B^{\nu',k'}_jf(x)}dx\right|\lesssim \frac{D^2\ga_j^{\nu,k}\ga_j^{\nu',k'}}{|k-k'|2^{j(1-\a/2)/2}(2^{j/2}|\nu-\nu'|)^{1+\a/2}},\\
\ds\textrm{ for }|\nu-\nu'|\neq 0,\,1\leq k,k'\leq |\nu-\nu'|^{-\a},\,k\neq k'.
\end{array}
\end{equation}
\eqref{oa11} and \eqref{oa12} imply:
\begin{equation}\label{oa13}
\begin{array}{lll}
\ds\left|\int_{\R^3}B^{\nu}_jf(x)\overline{B^{\nu'}_jf(x)}dx\right| & \ds\leq & \ds\sum_{1\leq k\leq |\nu-\nu'|^{-\a}}\left|\int_{\R^3}B^{\nu,k}_jf(x)\overline{B^{\nu',k}_jf(x)}dx\right|\\
& & \ds +\sum_{1\leq k\neq k'\leq |\nu-\nu'|^{-\a}}\left|\int_{\R^3}B^{\nu,k}_jf(x)\overline{B^{\nu',k'}_jf(x)}dx\right|\\
& \ds\lesssim & \ds\sum_{1\leq k\leq |\nu-\nu'|^{-\a}}\frac{D^2\ga_j^{\nu,k}\ga_j^{\nu',k}}{2^{j\a/2}(2^{j/2}|\nu-\nu'|)^{2-\a}}\\
& &\ds +\sum_{1\leq k\neq k'\leq |\nu-\nu'|^{-\a}}\frac{D^2\ga_j^{\nu,k}\ga_j^{\nu',k'}}{|k-k'|2^{\frac{j}{2}(1-\frac{\a}{2})}(2^{j/2}|\nu-\nu'|)^{1+\a/2}}\\
& \lesssim & \ds\frac{D^2\ga_j^{\nu}\ga_j^{\nu'}}{2^{j\a/2}(2^{j/2}|\nu-\nu'|)^{2-\a}},
\end{array}
\end{equation}
where we have used \eqref{oa10} in the last inequality and the fact that:
\begin{equation}\label{oa14}
\sup_{1\leq k\leq |\nu-\nu'|^{-\a}}\sum_{1\leq k'\leq |\nu-\nu'|^{-\a},\,k'\neq k}\frac{1}{|k-k'|}\lesssim \a |\log(|\nu-\nu'|)|.
\end{equation}
Since \eqref{oa13} yields the wanted estimate \eqref{oa2}, we are left with proving \eqref{oa11} and \eqref{oa12}.

\subsubsection{Proof of \eqref{oa11}} The estimate \eqref{oa11} will result of two integrations by parts with respect to tangential derivatives. By definition of $\nabb$, we have $\nabb h=\nabla h-(\nabla_\o h)\o$ for any function $h$ on $\R^3$. In particular, we have $\nabb(\xo)=0$ and $\nabb(\xo')=\o'-(\o'\cdot\o) \o$. Now, since $|\o'-(\o'\cdot\o)\o|^2=1-(\o'\cdot\o)^2$, this yields:
\begin{equation}\label{oa15} 
e^{i\la\xo-i\la'\xo'}=\frac{i}{\la'\sqrt{1-(\o'\cdot\o)^2}}\nabb_e(e^{i\la\xo-i\la'\xo'}),
\end{equation}
where
\begin{equation}\label{oa16} 
e=\frac{\o'-(\o'\cdot\o) \o}{\sqrt{1-(\o'\cdot\o)^2}}
\end{equation}
is a tangent vector with respect of the level surfaces of $\xo$. Similarly, we have:
\begin{equation}\label{oa17} 
e^{i\la\xo-i\la'\xo'}=-\frac{i}{\la\sqrt{1-(\o'\cdot\o)^2}}\nabb'_{e'}(e^{i\la\xo-i\la'\xo'}),
\end{equation}
where
\begin{equation}\label{oa18} 
e'=\frac{\o-(\o\cdot\o') \o'}{\sqrt{1-(\o'\cdot\o)^2}}
\end{equation}
is a tangent vector with respect of the level surfaces of $\xo'$. For $p\in\mathbb{Z}$, We introduce the notation $F_{j,k,p}(\xo)$:
\begin{equation}\label{oa19}
F_{j,k,p}(\xo)=\int_{0}^{+\infty}e^{i\lambda\xo}\psi(2^{-j}\lambda)\phi_k(\la)f(\lambda\o)(2^{-j}\la)^{p}\lambda^2 d\lambda.
\end{equation}
We integrate once by parts using \eqref{oa15} in $\int_{\R^3}B^{\nu,k}_jf(x)\overline{B^{\nu',k}_jf(x)}dx$ and we obtain:
\bea\label{oa20}
&&\int_{\R^3}B^{\nu,k}_jf(x)\overline{B^{\nu',k}_jf(x)}dx\\
\nn &=& 2^{-j}\int_{\R^3\times\S\times\S}\frac{i\nabb_eb(x,\o)\overline{b(x',\o')}}{\sqrt{1-(\o'\cdot\o)^2}}F_{j,k,0}(\xo)\overline{F_{j,k,-1}(\xo')}\eta^\nu_j(\o)\eta^{\nu'}_j(\o')d\o d\o'dx\\
\nn&& +2^{-j}\int_{\R^3\times\S\times\S}\frac{ib(x,\o)\overline{\nabb_eb(x',\o')}}{\sqrt{1-(\o'\cdot\o)^2}}F_{j,k,0}(\xo)\overline{F_{j,k,-1}(\xo')}\eta^\nu_j(\o)\eta^{\nu'}_j(\o')d\o d\o'dx.
\eea
We then integrate a second time by parts using \eqref{oa17} (so that there is at least one tangential derivative on $b(x,\o')$):
\bea\label{oa21}
&&\int_{\R^3}B^{\nu,k}_jf(x)\overline{B^{\nu',k}_jf(x)}dx\\
\nn& = &  2^{-2j}\int_{\R^3\times\S\times\S}\frac{\nabb_e\nabla_{e'}b(x,\o)\overline{b(x',\o')}}{1-(\o'\cdot\o)^2}F_{j,k,0}(\xo)\overline{F_{j,k,-2}(\xo')}\eta^\nu_j(\o)\eta^{\nu'}_j(\o')d\o d\o'dx\\
\nn&& +2^{-2j}\int_{\R^3\times\S\times\S}\frac{\nabb_eb(x,\o)\overline{\nabb_{e'}b(x',\o')}}{1-(\o'\cdot\o)^2}F_{j,k,0}(\xo)\overline{F_{j,k,-2}(\xo')}\eta^\nu_j(\o)\eta^{\nu'}_j(\o')d\o d\o'dx\\
\nn&& +2^{-2j}\int_{\R^3\times\S\times\S}\frac{\nabla_{e'}b(x,\o)\overline{\nabla_eb(x',\o')}}{1-(\o'\cdot\o)^2}F_{j,k,-1}(\xo)\overline{F_{j,k,-1}(\xo')}\eta^\nu_j(\o)\eta^{\nu'}_j(\o')d\o d\o'dx\\
\nn&& +2^{-2j}\int_{\R^3\times\S\times\S}\frac{b(x,\o)\overline{\nabb'_{e'}\nabla_eb(x',\o')}}{1-(\o'\cdot\o)^2}F_{j,k,-1}(\xo)\overline{F_{j,k,-1}(\xo')}\eta^\nu_j(\o)\eta^{\nu'}_j(\o')d\o d\o'dx.
\eea

\vspace{0.2cm}

\noindent{\bf Control of the right-hand side of \eqref{oa21}.} We now estimate the four terms in \eqref{oa21}. Using the fact that:
\begin{equation}\label{oa22}
\o\cdot\o'=1-\frac{|\o-\o'|^2}{2},
\end{equation}
we obtain the following expansions:
\begin{equation}\label{oa23}
\frac{1}{1-(\o\cdot\o')^2}=\frac{1}{|\nu-\nu'|^2}\left(1+\sum_{p+q\geq 1}c^1_{p,q}\left(\frac{\o-\nu}{|\nu-\nu'|}\right)^p\left(\frac{\o'-\nu'}{|\nu-\nu'|}\right)^q\right),
\end{equation}
and
\begin{equation}\label{oa24}
e=\frac{\nu'-(\nu'\cdot\nu) \nu}{\sqrt{1-(\nu'\cdot\nu)^2}}+\sum_{p+q\geq 1}c^2_{p,q}\left(\frac{\o-\nu}{|\nu-\nu'|}\right)^p\left(\frac{\o'-\nu'}{|\nu-\nu'|}\right)^q,
\end{equation}
where $c^1{p,q}$ and $c^2_{p,q}$ are constants. The expansions \eqref{oa23} and \eqref{oa24} allow us to rewrite the four terms of $\int_{\R^3}B_j^{\nu,k}f(x)\overline{B_j^{\nu',k}f(x)}dx$ such that one may separate the terms depending of $(\la,\o)$ from the terms depending on $(\la',\o')$. For instance, the first term in the right-hand side of \eqref{oa21} becomes:
\begin{equation}\label{oa25}
\begin{array}{r}
\ds\frac{2^{-j}}{(2^{j/2}|\nu-\nu'|)^2}\sum_{p+q\geq 0}c_{p,q}\int_{\R^3}\left(\int_{\S}\nabb\nabla b(x,\o)F_{j,k,0}(\xo)\left(\frac{\o-\nu}{|\nu-\nu'|}\right)^p\eta^\nu_j(\o)d\o\right)\\ 
\ds\times\left(\int_{\S}b(x,\o')F_{j,k,-2}(\xo')\left(\frac{\o'-\nu'}{|\nu-\nu'|}\right)^q\eta^{\nu'}_j(\o')d\o'\right)dx,
\end{array}
\end{equation}
where $c_{p,q}$ are constants. Since we have:
\begin{equation}\label{oa26}
\frac{|\o-\nu|}{|\nu-\nu'|}\lesssim \frac{1}{2^{j/2}|\nu-\nu'|}\textrm{ and }\frac{|\o'-\nu'|}{|\nu-\nu'|}\lesssim \frac{1}{2^{j/2}|\nu-\nu'|},
\end{equation}
the terms in the expansion \eqref{oa25} have more and more decay, and it is enough to consider the first one. We have:
\begin{equation}\label{oa27}
\begin{array}{rr}
&\ds\frac{2^{-j}}{(2^{j/2}|\nu-\nu'|)^2}\bigg|\int_{\R^3}\left(\int_{\S}\nabb\nabla b(x,\o)F_{j,k,0}(\xo)\eta^\nu_j(\o)d\o\right)\\[3mm]
&\ds\times\left(\int_{\S}b(x,\o')F_{j,k,-2}(\xo')\eta^{\nu'}_j(\o')d\o'\right)dx\bigg|\\[3mm]
\leq &\ds\frac{2^{-j}}{(2^{j/2}|\nu-\nu'|)^2}\left(\int_{\S}\norm{\nabb\nabla b}_{\le{2}}\norm{F_{j,k,0}}_{L^\infty_{\xo}}\eta^\nu_j(\o)d\o\right)\\[3mm]
&\ds\times\normm{\int_{\S}b(x,\o')F_{j,k,-2}(\xo')\eta^{\nu'}_j(\o')d\o'}_{\le{2}}.
\end{array}
\end{equation}
Using the estimate for the diagonal term \eqref{diagonal1} yields:
\begin{equation}\label{oa28}
\normm{\int_{\S}b(x,\o')F_{j,k,-2}(\xo')\eta^{\nu'}_j(\o')d\o'}_{\le{2}}\lesssim D\ga_j^{\nu',k'}.
\end{equation}
Using Cauchy Schwartz in $\la$ together with the size of the support of $\phi_k$ yields:
\begin{equation}\label{oa29}
\norm{F_{j,k,0}}_{L^\infty_{\xo}}\lesssim 2^{3j/2}|\nu-\nu'|^{\frac{\a}{2}}\norm{\psi(2^{-j}\la)\phi_k(\la)f(\la\o)\la}_{L^2_\la}.
\end{equation}
Finally, the assumption \eqref{hyp1} on $b(x,\o)$, the size of the support in $\o$, \eqref{oa27}, \eqref{oa28} and \eqref{oa29} imply:
\begin{equation}\label{oa30}
\begin{array}{l}
\ds\frac{2^{-j}}{(2^{j/2}|\nu-\nu'|)^2}\bigg|\int_{\R^3}\left(\int_{\S}\nabb\nabla b(x,\o)F_{j,k,0}(\xo)\eta^\nu_j(\o)d\o\right)\\
\ds\times\left(\int_{\S}b(x,\o')F_{j,k,-2}(\xo')\eta^{\nu'}_j(\o')d\o'\right)dx\bigg|\lesssim \frac{D^2|\nu-\nu'|^{\frac{\a}{2}}}{(2^{j/2}|\nu-\nu'|)^2}\ga_j^{\nu,k}\ga_j^{\nu',k'},
\end{array}
\end{equation}
which satisfies the wanted estimate \eqref{oa11}. The last term in the right-hand side of \eqref{oa21} is estimated exactly in the same way. 

\vspace{0.2cm}

\noindent{\bf Control of the second term in the right-hand side of \eqref{oa21}.} We still need to estimate the second and the third term in the right-hand side of \eqref{oa21}. Estimating them directly would yield the estimate \eqref{oa4} and ultimately the log-loss \eqref{oa5}. Thus, we need to integrate by parts once more. We first consider the second term in the right-hand side of \eqref{oa21}. Integrating by parts using \eqref{oa17} yields:
\bea\label{oa31}
&& 2^{-2j}\int_{\R^3\times\S\times\S}\frac{\nabb_eb(x,\o)\overline{\nabb_{e'}b(x',\o')}}{1-(\o'\cdot\o)^2}F_{j,k,0}(\xo)\overline{F_{j,k,-2}(\xo')}\eta^\nu_j(\o)\eta^{\nu'}_j(\o')d\o d\o'dx\\
\nn& =& i2^{-3j}\int_{\R^3\times\S\times\S}\frac{\nabla_{e'}\nabb_eb(x,\o)\overline{\nabla_{e'}b(x',\o')}}{(1-(\o'\cdot\o)^2)^{3/2}}F_{j,k,-1}(\xo)\overline{F_{j,k,-2}(\xo')}\eta^\nu_j(\o)\eta^{\nu'}_j(\o')d\o d\o'dx\\
\nn&& + i2^{-3j}\int_{\R^3\times\S\times\S}\frac{\nabb_eb(x,\o)\overline{\nabb_{e'}\nabla_{e'}b(x',\o')}}{(1-(\o'\cdot\o)^2)^{3/2}}F_{j,k,-1}(\xo)\overline{F_{j,k,-2}(\xo')}\eta^\nu_j(\o)\eta^{\nu'}_j(\o')d\o d\o'dx.
\eea
The two terms in the right-hand side of \eqref{oa31} are estimated in the same way, so we only consider  the first one. It is estimated by:
\bea\label{oa32}
\nn&& 2^{-3j}\bigg|\int_{\R^3\times\S\times\S}\frac{\nabla_{e'}\nabb_eb(x,\o)\overline{\nabla_{e'}b(x',\o')}}{(1-(\o'\cdot\o)^2)^{3/2}}F_{j,k,-1}(\xo)\overline{F_{j,k,-2}(\xo')}\eta^\nu_j(\o)\eta^{\nu'}_j(\o')d\o d\o'dx\bigg|\\
\nn&\leq& 2^{-3j}\int_{\S\times\S}\frac{1}{(1-(\o'\cdot\o)^2)^{3/2}}\norm{\nabla\nabb b(x,\o)}_{\le{2}}\norm{F_{j,k,-1}}_{L^\infty_{\xo}}\\
\nn&& \times\norm{\nabla b(x',\o')}_{\lprime{\infty}{2}}\norm{F_{j,k,-2}}_{L^2_{\xo'}}\eta^\nu_j(\o)\eta^{\nu'}_j(\o')d\o d\o'\\
&\lesssim& \frac{D^2|\nu-\nu'|^{\frac{\a}{2}}}{(2^{j/2}|\nu-\nu'|)^3}\ga_j^{\nu,k}\ga_j^{\nu',k'},
\eea
where we have used Plancherel to estimate $\norm{F_{j,k,-2}}_{L^2_{\xo'}}$, Cauchy-Schwartz in $\o$ and $\o'$, the assumption \eqref{hyp1} on $b$, and the estimate \eqref{oa29}. \eqref{oa32} satisfies the wanted estimate \eqref{oa11}.

\vspace{0.2cm}

\noindent{\bf Control of the third term in the right-hand side of \eqref{oa21} and end of the proof of \eqref{oa11}.} Finally, we consider the third term in the right-hand side of \eqref{oa21}. Neither of the two terms $\nabla_{e'}b$ and $\nabla_eb$ contain tangential derivatives, so integrating by parts directly would require to control two normal derivatives of $b$, which is not part of the assumptions \eqref{hyp1}. 
We first remark using the definition \eqref{oa16} of $e$ and \eqref{oa18} of $e'$ that:
\begin{equation}\label{oa34} 
e+e'=\frac{(1-\o'\cdot\o)(\o+\o')}{\sqrt{1-(\o'\cdot\o)^2}},
\end{equation}
 which yields the estimate:
\begin{equation}\label{oa34bis} 
e+e'\lesssim |\nu-\nu'|.
\end{equation} 
This allows us to rewrite the third term in the right-hand side of \eqref{oa21} as:
\begin{equation}\label{oa35}
\begin{array}{l}
\ds 2^{-2j}\int_{\R^3\times\S\times\S}\frac{\nabla_{e'}b(x,\o)\overline{\nabla_eb(x',\o')}}{1-(\o'\cdot\o)^2}F_{j,k,-1}(\xo)\overline{F_{j,k,-1}(\xo')}\eta^\nu_j(\o)\eta^{\nu'}_j(\o')d\o d\o'dx\\
\ds = 2^{-2j}\int_{\R^3\times\S\times\S}\frac{\nabb_{e}b(x,\o)\overline{\nabla_eb(x',\o')}}{1-(\o'\cdot\o)^2}F_{j,k,-1}(\xo)\overline{F_{j,k,-1}(\xo')}\eta^\nu_j(\o)\eta^{\nu'}_j(\o')d\o d\o'dx\\
\ds +2^{-2j}\int_{\R^3\times\S\times\S}\frac{\nabla_{e+e'}b(x,\o)\overline{\nabla_eb(x',\o')}}{1-(\o'\cdot\o)^2}F_{j,k,-1}(\xo)\overline{F_{j,k,-1}(\xo')}\eta^\nu_j(\o)\eta^{\nu'}_j(\o')d\o d\o'dx.
\end{array}
\end{equation}
The first term in the right-hand side of \eqref{oa35} is estimated in exactly as we proceeded for the second term in the right-hand side of \eqref{oa21} (i.e. by performing an additional integration by parts with the help of \eqref{oa17}). 
The second term in the right-hand side of \eqref{oa35} is estimated directly by:
\bea\label{oa36}
\nn&& 2^{-2j}\bigg|\int_{\R^3\times\S\times\S}\frac{\nabla_{e+e'}b(x,\o)\overline{\nabla_eb(x',\o')}}{1-(\o'\cdot\o)^2}F_{j,k,-1}(\xo)\overline{F_{j,k,-1}(\xo')}\eta^\nu_j(\o)\eta^{\nu'}_j(\o')d\o d\o'dx\bigg|\\
\nn&\leq& 2^{-2j}\int_{\S\times\S}\frac{|e+e'|}{1-(\o'\cdot\o)^2}\norm{\nabla b(x,\o)}_{\l{\infty}{2}}\norm{F_{j,k,-1}}_{L^2_{\xo}}\norm{\nabla b(x',\o')}_{\lprime{\infty}{2}}\\
\nn&&\times\norm{F_{j,k,-1}}_{L^2_{\xo'}}\eta^\nu_j(\o)\eta^{\nu'}_j(\o')d\o d\o'\\
&\lesssim& \frac{D^2}{2^{j/2}(2^{j/2}|\nu-\nu'|)}\gamma_j^{\nu,k}\gamma_{j}^{\nu',k},
\eea
where we have used Plancherel to estimate $\norm{F_{j,k,-1}}_{L^2_{\xo}}$ and $\norm{F_{j,k,-1}}_{L^2_{\xo'}}$, Cauchy-Schwartz in $\o$ and $\o'$, the assumption \eqref{hyp1} on $b$, and the estimate 
\eqref{oa34bis}. \eqref{oa36} satisfies the wanted estimate \eqref{oa11} for $0<\a\leq 1$. We now control all the terms in the right-hand side of \eqref{oa21} which concludes the proof of \eqref{oa11}.

\subsubsection{Proof of \eqref{oa12}} The estimate \eqref{oa12} will result of two integrations by parts, one with respect to the normal derivative, and one with respect to tangential derivatives. We first integrate by parts with respect to $\partial_{\xo}$ in $\int_{\R^3}B^{\nu,k}_jf(x)\overline{B^{\nu',k'}_jf(x)}dx$ using \eqref{of4}. We obtain:
\begin{equation}\label{oa43}
\begin{array}{l}
\ds \int_{\R^3}B^{\nu,k}_jf(x)\overline{B^{\nu',k'}_jf(x)}dx = \int_{\R^3\times\S\times\S}\int_0^{+\infty}\int_0^{+\infty}\frac{i}{\la-\la'\o\cdot\o'}\\
\ds \times\partial_{\xo}b(x,\o)b(x,\o')\eta_j^\nu(\o)\eta_j^{\nu'}(\o')\psi(2^{-j}\la)\psi(2^{-j}\la')\phi_k(\la)\phi_{k'}(\la') \\
\ds \times f(\la\o)f(\la'\o')\la^2 {\la'}^2d\la d\la' d\o d\o'dx+\int_{\R^3\times\S\times\S}\int_0^{+\infty}\int_0^{+\infty}\frac{i}{\la-\la'\o\cdot\o'}\\
\ds \times b(x,\o)\partial_{\xo}b(x,\o')\eta_j^\nu(\o)\eta_j^{\nu'}(\o')\psi(2^{-j}\la)\psi(2^{-j}\la')\phi_k(\la)\phi_{k'}(\la') \\
\ds\times f(\la\o)f(\la'\o')\la^2 {\la'}^2d\la d\la' d\o d\o'dx.
\end{array}
\end{equation}
We then integrate a second time by parts using \eqref{oa15} for the first term in the right-hand side of \eqref{oa43}, and using \eqref{oa17} for the second term in the right-hand side of \eqref{oa43} (so that there is at least one tangential derivative on $b(x,\o')$). We obtain:
\begin{equation}\label{oa44}
\begin{array}{l}
\ds \int_{\R^3}B^{\nu,k}_jf(x)\overline{B^{\nu',k'}_jf(x)}dx = \int_{\R^3\times\S\times\S}\int_0^{+\infty}\int_0^{+\infty}\frac{1}{(\la-\la'\o\cdot\o')\la'\sqrt{1-(\o\cdot\o')^2}}\\
\ds \times\nabb_e\partial_{\xo}b(x,\o)b(x,\o')\eta_j^\nu(\o)\eta_j^{\nu'}(\o')\psi(2^{-j}\la)\psi(2^{-j}\la')\phi_k(\la)\phi_{k'}(\la') f(\la\o)f(\la'\o')\\
\ds\times \la^2 {\la'}^2 d\la d\la' d\o d\o'dx+\int_{\R^3\times\S\times\S}\int_0^{+\infty}\int_0^{+\infty}\frac{1}{(\la-\la'\o\cdot\o')\la'\sqrt{1-(\o\cdot\o')^2}}\\
\ds \times\partial_{\xo}b(x,\o)\nabla_eb(x,\o')\eta_j^\nu(\o)\eta_j^{\nu'}(\o')\psi(2^{-j}\la)\psi(2^{-j}\la')\phi_k(\la)\phi_{k'}(\la') f(\la\o)f(\la'\o')\\
\ds\times\la^2 {\la'}^2d\la d\la' d\o d\o'dx + \int_{\R^3\times\S\times\S}\int_0^{+\infty}\int_0^{+\infty}\frac{1}{(\la-\la'\o\cdot\o')\la\sqrt{1-(\o\cdot\o')^2}}\\
\ds \times\nabla_{e'}b(x,\o)\partial_{\xo}b(x,\o')\eta_j^\nu(\o)\eta_j^{\nu'}(\o')\psi(2^{-j}\la)\psi(2^{-j}\la')\phi_k(\la)\phi_{k'}(\la') f(\la\o)f(\la'\o')\\
\ds\times\la^2 {\la'}^2d\la d\la' d\o d\o'dx + \int_{\R^3\times\S\times\S}\int_0^{+\infty}\int_0^{+\infty}\frac{1}{(\la-\la'\o\cdot\o')\la\sqrt{1-(\o\cdot\o')^2}}\\
\ds \times b(x,\o)\nabb_{e'}\partial_{\xo}b(x,\o')\eta_j^\nu(\o)\eta_j^{\nu'}(\o')\psi(2^{-j}\la)\psi(2^{-j}\la')\phi_k(\la)\phi_{k'}(\la') f(\la\o)f(\la'\o')\\
\ds\times\la^2 {\la'}^2d\la d\la' d\o d\o'dx.
\end{array}
\end{equation}

Since $|\la-k2^j|\nu-\nu'|^\a|\leq 2^j|\nu-\nu'|^\a$ on the support of $\phi_k$ and $|\la'-k'2^j|\nu-\nu'|^\a|\leq 2^j|\nu-\nu'|^\a$ on the support of $\phi_{k'}$, we have the following expansion:
\begin{equation}\label{of6b}
\begin{array}{r}
\ds\frac{1}{\la-\la'\o\cdot\o'}=\frac{1}{(k-k')2^j|\nu-\nu'|^\a}\sum_{p,q,r\geq 0}c_{p,q,r}\left(\frac{\la-k2^j|\nu-\nu'|^\a}{(k-k')2^j|\nu-\nu'|^\a}\right)^p\\\ds\times\left(\frac{\la'-k'2^j|\nu-\nu'|^\a}{(k-k')2^j|\nu-\nu'|^\a}\right)^q\left(\frac{\la'|\o-\o'|^2}{(k-k')2^j|\nu-\nu'|^\a}\right)^r
\end{array}
\end{equation}
For $p\in\mathbb{Z}$ and $q\in\N$, we introduce the notation $F_{j,k,p,q}(\xo)$:
\begin{equation}\label{of7b}
F_{j,k,p,q}(\xo)=\int_{0}^{+\infty}e^{i\lambda\xo}\psi(2^{-j}\lambda)\phi_k(\la)f(\lambda\o)(2^{-j}\la)^{p}\left(\frac{\la-k2^j|\nu-\nu'|^\a}{2^j|\nu-\nu'|^\a}\right)^q\lambda^2 d\lambda.
\end{equation}
\eqref{oa44}, \eqref{of6b} and \eqref{of7b} yield:
\begin{equation}\label{oa37}
\int_{\R^3}B^{\nu,k}_jf(x)\overline{B^{\nu',k'}_jf(x)}dx=\sum_{p,q,r\geq 0}c_{p,q,r}(A^{1,1}_{p,q,r}+A^{1,2}_{p,q,r}+A^{2,1}_{p,q,r}+A^{2,2}_{p,q,r}),
\end{equation}
where $A^{1,1}_{p,q,r}$, $A^{1,2}_{p,q,r}$, $A^{2,1}_{p,q,r}$ and $A^{2,2}_{p,q,r}$ are given by:
\begin{equation}\label{oa38}
\begin{array}{l}
\ds A^{1,1}_{p,q,r}=\ds \frac{1}{(k-k')^{p+q+r+1}2^{2j}|\nu-\nu'|^{\a}}\int_{\R^3\times\S\times\S}\frac{1}{\sqrt{1-(\o\cdot\o')^2}}\\
\ds\times\left(\frac{|\o-\o'|^2}{|\nu-\nu'|^\a}\right)^r\nabb_e\partial_{\xo}b(x,\o)F_{j,k,0,p}(\xo)\overline{b(x,\o')F_{j,k',r-1,q}(\xo')}d\o d\o'dx,
\end{array}
\end{equation}
\begin{equation}\label{oa39}
\begin{array}{l}
\ds A^{1,2}_{p,q,r}=\ds \frac{1}{(k-k')^{p+q+r+1}2^{2j}|\nu-\nu'|^{\a}}\int_{\R^3\times\S\times\S}\frac{1}{\sqrt{1-(\o\cdot\o')^2}}\\
\ds\times\left(\frac{|\o-\o'|^2}{|\nu-\nu'|^\a}\right)^r\partial_{\xo}b(x,\o)F_{j,k,0,p}(\xo)\overline{\nabla_eb(x,\o')F_{j,k',r-1,q}(\xo')}d\o d\o'dx,
\end{array}
\end{equation}
\begin{equation}\label{oa40}
\begin{array}{l}
\ds A^{2,1}_{p,q,r}=\ds \frac{1}{(k-k')^{p+q+r+1}2^{2j}|\nu-\nu'|^{\a}}\int_{\R^3\times\S\times\S}\frac{1}{\sqrt{1-(\o\cdot\o')^2}}\\
\ds\times\left(\frac{|\o-\o'|^2}{|\nu-\nu'|^\a}\right)^r\nabla_{e'}b(x,\o)F_{j,k,-1,p}(\xo)\overline{\partial_{\xo}b(x,\o')F_{j,k',r,q}(\xo')}d\o d\o'dx,
\end{array}
\end{equation}
and
\begin{equation}\label{oa41}
\begin{array}{l}
\ds A^{2,2}_{p,q,r}=\ds \frac{1}{(k-k')^{p+q+r+1}2^{2j}|\nu-\nu'|^{\a}}\int_{\R^3\times\S\times\S}\frac{1}{\sqrt{1-(\o\cdot\o')^2}}\\
\ds\times\left(\frac{|\o-\o'|^2}{|\nu-\nu'|^\a}\right)^rb(x,\o)F_{j,k,-1,p}(\xo)\overline{\nabb_{e'}\partial_{\xo}b(x,\o')F_{j,k',r,q}(\xo')}d\o d\o'dx.
\end{array}
\end{equation}

\vspace{0.2cm}

\noindent{\bf Control of $A^{1,1}_{p,q,r}$, $A^{1,2}_{p,q,r}$, $A^{2,1}_{p,q,r}$ and $A^{2,2}_{p,q,r}$.} We start by evaluating $A^{1,2}_{p,q,r}$. We have:
\bea\label{oa42}
\nn  |A^{1,2}_{p,q,r}|&\leq &\!\!\!\!\ds\frac{1}{(k-k')^{p+q+r+1}2^{2j}|\nu-\nu'|^{\a}}\int_{\R^3\times\S\times\S}\frac{1}{\sqrt{1-(\o\cdot\o')^2}}\norm{\nabla b(x,\o)}_{\l{\infty}{2}}\\
\nn&&\times\norm{F_{j,k,0,p}}_{L^2_{\xo}}\norm{\nabla b(x,\o')}_{\lprime{\infty}{2}}\norm{F_{j,k',r-1,q}}_{L^2_{\xo'}}d\o d\o'\\
&\lesssim &\!\!\!\!\ds\frac{D^2\gamma_j^{\nu,k}\gamma_j^{\nu',k'}}{(k-k')^{p+q+r+1}2^{j/2(1-\a)}(2^{j/2}|\nu-\nu'|)^{1+\a}},
\eea
where we have used Plancherel to estimate $\norm{F_{j,k,0,p}}_{L^2_{\xo}}$ and $\norm{F_{j,k',r,q}}_{L^2_{\xo'}}$, Cauchy-Schwartz in $\o$ and $\o'$ and the assumption \eqref{hyp1} on $b$. We control $A^{2,1}_{p,q,r}$ in the same way.

It remains to estimate $A^{1,1}_{p,q,r}$ and $A^{2,2}_{p,q,r}$. They are controlled in the the same way, so we focus on estimating $A^{1,1}_{p,q,r}$. Using the expansions \eqref{oa23} and \eqref{oa24}, we obtain:
\begin{equation}\label{oa45}
A^{1,1}_{p,q,r}=\sum_{l,m\geq 0}c_{p,q,r,l,m}A^{1,1}_{p,q,r,l,m},
\end{equation}
where $A^{1,1}_{p,q,r,l,m}$ are given by:
\begin{equation}\label{oa46}
\begin{array}{ll}
\ds A^{1,1}_{p,q,r,l,m}= & \ds \frac{1}{(k-k')^{p+q+r+1}2^{j(3/2-\a/2)}(2^{j/2}|\nu-\nu'|)^{1+\a}}\\
& \ds\times\int_{\R^3}\left(\int_{\S}\left(\frac{\o-\nu}{|\nu-\nu'|}\right)^l\nabb\nabla b(x,\o)F_{j,k,0,p}(\xo)\eta_j^{\nu}(\o)d\o\right)\\
& \ds\times\overline{\left(\int_{\S}\left(\frac{\o'-\nu'}{|\nu-\nu'|}\right)^mb(x,\o')F_{j,k',r-1,q}(\xo')\eta_j^{\nu'}(\o')d\o'\right)}dx,
\end{array}
\end{equation}
The terms in the expansion \eqref{oa45} have more and more decay, and it is enough to consider the first one. We have:
\begin{equation}\label{oa47}
\begin{array}{ll}
&\ds\frac{1}{(k-k')^{p+q+r+1}2^{j(3/2-\a/2)}(2^{j/2}|\nu-\nu'|)^{1+\a}}\\
&\ds\times\bigg|\int_{\R^3}\left(\int_{\S}\nabb\nabla b(x,\o)F_{j,k,0,p}(\xo)\eta_j^{\nu}(\o)d\o\right)\\
&\ds\times\overline{\left(\int_{\S}b(x,\o')F_{j,k',r-1,q}(\xo')\eta_j^{\nu'}(\o')d\o'\right)}dx\bigg|\\
\ds \leq &\ds\frac{1}{(k-k')^{p+q+r+1}2^{j(3/2-\a/2)}(2^{j/2}|\nu-\nu'|)^{1+\a}}\\
&\ds\times\left(\int_{\S}\norm{\nabb\nabla b(x,\o)}_{\le{2}}\norm{F_{j,k,0,p}}_{L^\infty_{\xo}}\eta_j^{\nu}(\o)d\o\right)\\
&\ds\times\normm{\int_{\S}b(x,\o')F_{j,k',r-1,q}(\xo')\eta_j^{\nu'}(\o')d\o'}_{\le{2}}
\end{array}
\end{equation}
Using the estimate for the diagonal term \eqref{diagonal1} yields:
\begin{equation}\label{oa48}
\normm{\int_{\S}b(x,\o')F_{j,k',r-1,q}(\xo')\eta_j^{\nu'}(\o')d\o'}_{\le{2}}\lesssim D\ga_j^{\nu',k'}.
\end{equation}
Finally, the assumption \eqref{hyp1} on $b(x,\o)$, the size of the support in $\o$, the bound \eqref{oa29} on $\norm{F_{j,k,0,p}}_{L^\infty_{\xo}}$, \eqref{oa46}, \eqref{oa47} and \eqref{oa48} imply:
\begin{equation}\label{oa49}
\begin{array}{l}
\ds |A^{1,1}_{p,q,r}| \lesssim\frac{D^2|\nu-\nu'|^\a\ga_j^{\nu,k}\ga_j^{\nu',k'}}{(k-k')^{p+q+r+1}2^{j/2(1-\a)}(2^{j/2}|\nu-\nu'|)^{1+\a}}.
\end{array}
\end{equation}
Summing in $p, q, r$ the estimate \eqref{oa42} and its analog for $A^{2,1}_{p,q,r}$ together with \eqref{oa49} and its analog for $A^{2,2}_{p,q,r}$, and using \eqref{oa37}, we obtain the wanted estimate \eqref{oa12}. 

\subsubsection{End of the proof of Proposition \ref{orthoangle}}

We have proved the estimates \eqref{oa11} and \eqref{oa12} in the two previous sections. Since \eqref{oa11} and \eqref{oa12} yield \eqref{oa2} (see section \ref{sec:keyest}), this concludes the proof of Proposition \ref{orthoangle}. \QED

\subsection{Proof of Proposition \ref{diagonal} (control of the diagonal term)}\label{sec:diagonal}

We have to prove \eqref{diagonal1}:
\begin{equation}\label{di1}
\norm{B^\nu_jf}_{\le{2}}\lesssim D\ga^\nu_j.
\end{equation}
Recall that $B^\nu_j$ is given by:
\begin{equation}\label{di2}
B^\nu_jf(x)=\int_{\S}b(x,\o)F_j(\xo)\eta_j^\nu(\o)d\o,
\end{equation}
where $F_j(\xo)$ is defined by:
\begin{equation}\label{di3}
F_j(\xo)=\int_0^{+\infty}e^{i\la\xo}\psi(2^{-j}\la)f(\la\o)\la^2d\la.
\end{equation}
We decompose $B^\nu_j$ in the sum of two terms:
\begin{equation}\label{di4}
B^\nu_jf(x)=b(x,\nu)\int_{\S}F_j(\xo)\eta_j^\nu(\o)d\o+\int_{\S}(b(x,\o)-b(x,\nu))F_j(\xo)\eta_j^\nu(\o)d\o.
\end{equation}

Notice that the first term in the right-hand side of \eqref{di4} is equal to 
\begin{equation}\label{di5}
b(x,\nu)\int_{\S}F_j(\xo)\eta_j^\nu(\o)d\o=b(x,\nu)\mathcal{F}^{-1}(\psi(2^{-j}\la)\eta_j^\nu(\o)f(\la\o))(x),
\end{equation}
where $\mathcal{F}$ denotes the Fourier transform on $\R^3$. Now, the assumption \eqref{hyp1} on $b$ imply that $\norm{b}_{\le{\infty}}\lesssim D$. Together with \eqref{di5}, this yields:
\begin{equation}\label{di6}
\normm{b(x,\nu)\int_{\S}F_j(\xo)\eta_j^\nu(\o)d\o}_{\le{2}}\lesssim D\gamma^\nu_j.
\end{equation}

We turn to the second term in the right-hand side of \eqref{di4}. We have:
\begin{equation}\label{di7}
\begin{array}{ll}
&\ds\normm{\int_{\S}(b(x,\o)-b(x,\nu))F_j(\xo)\eta_j^\nu(\o)d\o}_{\le{2}}\\
\ds\leq &\ds\int_{\S}\norm{b(x,\o)-b(x,\nu)}_{\l{\infty}{2}}\norm{F_j}_{L^2_{\xo}}\eta_j^\nu(\o)d\o.
\end{array}
\end{equation}
Now, $\h{1/2+\a}$ embeds in $\l{\infty}{2}$ for any $\a>0$, thus:
\begin{equation}\label{di8}
\norm{b(x,\o)-b(x,\nu)}_{\l{\infty}{2}}\lesssim \norm{b(x,\o)-b(x,\nu)}_{\h{1/2+\a}}\lesssim |\o-\nu|\norm{\po b}_{H^{1/2+\a}}.
\end{equation}
Together with \eqref{di7}, this yields:
\begin{equation}\label{di9}
\begin{array}{ll}
&\ds\normm{\int_{\S}(b(x,\o)-b(x,\nu))F_j(\xo)\eta_j^\nu(\o)d\o}_{\le{2}}\\
\ds\leq &\ds\int_{\S}|\o-\nu|\norm{\po b}_{H^{1/2+\a}(\R^3)}\norm{F_j}_{L^2_{\xo}}\eta_j^\nu(\o)d\o\\
\ds\lesssim &\ds D\gamma_j^\nu,
\end{array}
\end{equation}
where we have used Plancherel to estimate $\norm{F_j}_{L^2_{\xo}}$, Cauchy-Schwartz in $\o$, the assumption \eqref{hyp1} on $b$, and the fact that $|\o-\nu|\lesssim 2^{-j/2}$ on the support of $\eta_j^\nu$.

Finally, \eqref{di4}, \eqref{di6} and \eqref{di9} yield the wanted estimate \eqref{di1} which concludes the proof of Proposition \ref{diagonal}. \QED.

\section{Proof of Theorem \ref{th1} ($L^2$ boundedness for Fourier integral operator)}\label{sec:th1}

\subsection{The basic computation}

We start the proof of Theorem \ref{th1} with the following instructive computation:
\begin{equation}\label{bisb5}
\begin{array}{ll}
\ds\norm{Uf}_{\ll{2}}\leq &\ds\int_{\S}\normm{b(x,\o)\int_{0}^{+\infty}e^{i\lambda u}f(\lambda\o)\lambda^2 d\lambda}_{\ll{2}}d\o\\
&\ds\leq\int_{\S}\norm{b(x,\o)}_{\l{\infty}{2}}\normm{\int_{0}^{+\infty}e^{i\lambda u}f(\lambda\o)\lambda^2 d\lambda}_{L^2_{u}}d\o\\
&\ds\leq D\norm{\lambda f}_{L^2(\s)},
\end{array}
\end{equation}
where we have used Plancherel with respect to $\lambda$, Cauchy-Schwarz with respect to $\o$ 
and \eqref{hyp1} to bound $\norm{b}_{\l{\infty}{2}}$. \eqref{bisb5} misses the conclusion \eqref{l2} of Theorem \ref{th1} by a power of $\lambda$. Now, assume for a moment that we may replace a power of $\lambda$ by a derivative on $b(x,\o)$. Then, the same computation yields:
\begin{equation}\label{bisb6}
\begin{array}{ll}
&\ds\normm{\int_{\S}\int_{0}^{+\infty}\nabla b(x,\o) e^{i\lambda u}f(\lambda\o)\lambda d\lambda d\o}_{\ll{2}}\\
\ds\leq &\ds \int_{\S}\norm{\nabla b(x,\o)}_{\l{\infty}{2}}\normm{\int_{0}^{+\infty}e^{i\lambda u}f(\lambda\o)\lambda^2 d\lambda}_{L^2_{ u}}d\o\\
\ds\leq &\ds D\norm{f}_{L^2(\s)},
\end{array}
\end{equation}
which is \eqref{l2}. This suggests a strategy which consists in making integrations by parts to trade powers of $\lambda$ against derivatives of the symbol $b(x,\o)$. 

\subsection{Structure of the proof of Theorem \ref{th1}}

The proof of Theorem \ref{th1} proceeds in three steps. We first localize in frequencies of size $\la\sim 2^j$. We then localize the angle $\o$ in patches on the sphere $\S$ of diameter $2^{-j/2}$. Finally, we estimate the diagonal terms.

\subsubsection{Step 1: decomposition in frequency}
 
For the first step, we introduce $\varphi$  and $\psi$ two smooth compactly supported functions on $\R$ such that: 
\begin{equation}\label{bisb7}
\varphi(\lambda)+\sum_{j\geq 0}\psi(2^{-j}\lambda)=1\textrm{ for all }\lambda\in\R.
\end{equation}
We use \eqref{bisb7} to decompose $Uf$ as follows:
\begin{equation}\label{bisb8}
Uf(x)=\sum_{j\geq -1}U_jf(x),
\end{equation}
where for $j\geq 0$:
\begin{equation}\label{bisb9}
U_jf(x)=\int_{\S}\int_{0}^{+\infty}e^{i\lambda u}b(x,\o)\psi(2^{-j}\lambda)f(\lambda\o)\lambda^2 d\lambda d\o,
\end{equation}
and 
\begin{equation}\label{bisb10}
U_{-1}f(x)=\int_{\S}\int_{0}^{+\infty}e^{i\lambda u}b(x,\o)\varphi(\lambda)f(\lambda\o)\lambda^2 d\lambda d\o.
\end{equation}
This decomposition is classical and is known as the first dyadic decomposition (see \cite{stein}). The goal of this first step is to prove the following proposition:
\begin{proposition}\label{bisorthofreq}
The decomposition \eqref{bisb8} satisfies an almost orthogonality property:
\begin{equation}\label{bisorthofreq1}
\norm{Uf}_{\ll{2}}^2\lesssim\sum_{j\geq -1}\norm{U_jf}_{\ll{2}}^2+D^2\norm{f}^2_{\le{2}}.
\end{equation}
\end{proposition}
The proof of Proposition \ref{bisorthofreq} is postponed to section \ref{bissec:orthofreq}. 

\subsubsection{Step 2: decomposition in angle}

Proposition \ref{bisorthofreq} allows us to  estimate $\norm{U_jf}_{\ll{2}}$ instead of 
$\norm{Uf}_{\ll{2}}$. The analog of computation \eqref{bisb5} for $\norm{U_jf}_{\ll{2}}$ yields:
\begin{equation}\label{bisb5bis}
\begin{array}{l}
\ds\norm{U_jf}_{\ll{2}}\leq D\norm{\lambda\psi(2^{-j}\la) f}_{L^2(\s)}\lesssim D2^j\norm{\psi(2^{-j}\la) f}_{\le{2}},
\end{array}
\end{equation}
which misses the wanted estimate by a power of $2^j$. We thus need to perform a second dyadic decomposition (see \cite{stein}). We introduce a smooth partition of unity on the sphere $\S$:
\begin{equation}\label{bisb14}
\sum_{\nu\in\Gamma}\eta^\nu_j(\o)=1\textrm{ for all }\o\in\S, 
\end{equation}
where the support of $\eta^\nu_j$ is a patch on $\S$ of diameter $\sim 2^{-j/2}$. We use \eqref{bisb14} to decompose $U_jf$ as follows:
\begin{equation}\label{bisb15}
U_jf(x)=\sum_{\nu\in\Gamma}U^\nu_jf(x),
\end{equation}
where:
\begin{equation}\label{bisb16}
U^\nu_jf(x)=\int_{\S}\int_{0}^{+\infty}e^{i\lambda u}b(x,\o)\psi(2^{-j}\lambda)\eta^\nu_j(\o)f(\lambda\o)\lambda^2 d\lambda d\o.
\end{equation}
We also define:
\begin{equation}\label{bisdecf}
\begin{array}{l}
\ds\ga_{-1}=\norm{\varphi(\la)f}_{\le{2}},\, \ga_j=\norm{\psi(2^{-j}\la)f}_{\le{2}},\,j\geq 0, \\
\ds\ga^\nu_j=\norm{\psi(2^{-j}\la)\eta^\nu_j(\o)f}_{\le{2}},\,j\geq 0,\,\nu\in\Gamma, 
\end{array}
\end{equation}
which satisfy:
\begin{equation}\label{bisdecf1}
\norm{f}_{\le{2}}^2=\sum_{j\geq -1}\ga_j^2=\sum_{j\geq -1}\sum_{\nu\in\Gamma}(\ga^\nu_j)^2.
\end{equation}
The goal of this second step is to prove the following proposition:
\begin{proposition}\label{bisorthoangle}
The decomposition \eqref{bisb15} satisfies an almost orthogonality property:
\begin{equation}\label{bisorthoangle1}
\norm{U_jf}_{\ll{2}}^2\lesssim\sum_{\nu\in\Gamma}\norm{U^\nu_jf}_{\ll{2}}^2+D^2\ga_j^2.
\end{equation}
\end{proposition}
The proof of Proposition \ref{bisorthoangle} is postponed to section \ref{bissec:orthoangle}. 

\subsubsection{Step 3: control of the diagonal term}

Proposition \ref{bisorthoangle} allows us to  estimate $\norm{U^\nu_jf}_{\ll{2}}$ instead of $\norm{U_jf}_{\ll{2}}$. The analog of computation \eqref{bisb5} for $\norm{U^\nu_jf}_{\ll{2}}$ yields:
\begin{equation}\label{bisb5ter}
\begin{array}{rl}
\ds\norm{U_j^\nu f}_{\ll{2}}\leq &\ds\int_{\S}\norm{b(x,\o)}_{\l{\infty}{2}}\normm{\int_{0}^{+\infty}e^{i\lambda  u}\psi(2^{-j}\la)\eta^\nu_j(\o)f(\lambda\o)\lambda^2 d\lambda}_{L^2_{ u}}d\o\\
\ds\leq &\ds D\sqrt{\textrm{vol}(\textrm{supp}(\eta^\nu_j))}\norm{\lambda\psi(2^{-j}\la)\eta^\nu_j(\o)f}_{\le{2}}\\
\ds\lesssim &\ds D2^{j/2}\gamma^\nu_j,
\end{array}
\end{equation}
where the term $\sqrt{\textrm{vol}(\textrm{supp}(\eta^\nu_j))}$ comes from the fact that we apply Cauchy-Schwarz in $\o$. Note that we have used in \eqref{bisb5ter} the fact that the support of $\eta^\nu_j$ is 2 dimensional and has diameter $2^{-j/2}$ so that:
\begin{equation}\label{bissuppangle}
\sqrt{\textrm{vol}(\textrm{supp}(\eta^\nu_j))}\lesssim 2^{-j/2}.
\end{equation}
Now, \eqref{bisb5ter} still misses the wanted estimate by a power of $2^{j/2}$. Nevertheless, we are able to estimate the diagonal term:
\begin{proposition}\label{bisdiagonal}
The diagonal term $U^\nu_jf$ satisfies the following estimate:
\begin{equation}\label{bisdiagonal1}
\norm{U^\nu_jf}_{\ll{2}}\lesssim D\ga^\nu_j.
\end{equation}
\end{proposition}
The proof of Proposition \ref{bisdiagonal} is postponed to section \ref{bissec:diagonal}. 

\subsubsection{Proof of Theorem \ref{th1}}

Proposition \ref{bisorthofreq}, \ref{bisorthoangle} and \ref{bisdiagonal} 
immediately yield the proof of Theorem \ref{th1}. Indeed, \eqref{bisorthofreq1}, \eqref{bisdecf1}, \eqref{bisorthoangle1} and \eqref{bisdiagonal1} imply:
\begin{equation}\label{biscclth3}
\begin{array}{ll}
\ds\norm{Uf}_{\ll{2}}^2 & \ds\lesssim\sum_{j\geq -1}\norm{U_jf}_{\ll{2}}^2+D^2\norm{f}^2_{\le{2}}\\
& \ds\lesssim\sum_{j\geq -1}\sum_{\nu\in\Gamma}\norm{U^\nu_jf}_{\ll{2}}^2+D^2\sum_{j\geq -1}\gamma_j^2+D^2\norm{f}^2_{\le{2}}\\
& \ds\lesssim D^2\sum_{j\geq -1}\sum_{\nu\in\Gamma}(\gamma_j^\nu)^2+D^2\sum_{j\geq -1}\gamma_j^2+D^2\norm{f}^2_{\le{2}}\\
& \ds\lesssim D^2\norm{f}^2_{\le{2}},
\end{array}
\end{equation}
which is the conclusion of Theorem \ref{th1}. \QED

The remainder of section \ref{sec:th1} is dedicated to the proof of Proposition \ref{bisorthofreq}, \ref{bisorthoangle} and \ref{bisdiagonal}.

\subsection{Proof of Proposition \ref{bisorthofreq} (almost orthogonality in frequency)}\label{bissec:orthofreq}

We have to prove \eqref{bisorthofreq1}:
\begin{equation}\label{bisof1}
\norm{Uf}_{\ll{2}}^2\lesssim\sum_{j\geq -1}\norm{U_jf}_{\ll{2}}^2+D^2\norm{f}^2_{\le{2}}.
\end{equation}
This will result from the following inequality using Shur's Lemma:
\begin{equation}\label{bisof2}
\left|\int_{\s}U_jf(x)\overline{U_kf(x)}d\s \right| \lesssim D^22^{-\frac{|j-k|}{2}}\gamma_j\ga_k\textrm{ for }|j-k|> 2.
\end{equation}

\subsubsection{A first integration by parts} 

From now on, we focus on proving \eqref{bisof2}. We may assume $j\geq k+3$. We have:
\begin{equation}\label{bisof3}
\begin{array}{ll}
\ds\int_{\s}U_jf(x)\overline{U_kf(x)}d\s = & \ds\int_{\S}\int_{0}^{+\infty}\int_{\S}\int_{0}^{+\infty}\left(\int_{\s}e^{i\lambda u-i\la' u'}b(x,\o)\overline{b(x,\o')}d\s \right)\\
& \ds\times\psi(2^{-j}\lambda)f(\lambda\o)\lambda^2 \psi(2^{-k}\lambda')\overline{f(\lambda'\o')}(\lambda')^2 d\lambda d\o d\lambda' d\o'.
\end{array}
\end{equation}
We integrate by parts with respect to $\partial_{ u}$ in $\int_{\s}e^{i\lambda u-i\la' u'}b(x,\o)\overline{b(x,\o')}d\s $ using the coarea formula \eqref{coarea} and the fact that:
\begin{equation}\label{bisof4}
e^{i\lambda u-i\la' u'}=-\frac{i}{\la-\la'\frac{a}{a'}g(N,N')}\partial_{ u}(e^{i\lambda u-i\la' u'}),
\end{equation}
where we use the notation $u$ for $u(x,\o)$, $a$ for $a(x,\o)$, $N$ for $N(x,\o)$, $u'$ for $u(x,\o')$, $a'$ for $a(x,\o')$ and $N'$ for $N(x,\o')$. We will also use the notation $b$ for $b(x,\o)$, $b'$ for $b(x,\o')$, $\th$ for $\th(x,\o)$ and $\th'$ for $\th(x,\o')$. Using \eqref{bisof4}, we obtain:
\begin{equation}\label{bisof5}
\begin{array}{ll}
\ds\int_{\s}e^{i\lambda u-i\la' u'}b\overline{b'}d\s = & \ds i\int_{\s}e^{i\lambda u-i\la' u'}\frac{\partial_{ u}b\overline{b'}}{\la-\la'\frac{a}{a'}g(N,N')}d\s\\
&\ds +i\int_{\s}e^{i\lambda u-i\la' u'}\frac{b\partial_{ u}\overline{b'}}{\la-\la'\frac{a}{a'}g(N,N')}d\s \\
&\ds +i\int_{\s}e^{i\lambda u-i\la' u'}\frac{b\overline{b'}\trt}{\la-\la'\frac{a}{a'}g(N,N')}d\s\\
&\ds +i\la'\int_{\s}e^{i\lambda u-i\la' u'}\frac{b\overline{b'}(\frac{\nabn a}{a'}g(N,N')-\frac{a\nabn a'}{a'^2}g(N,N'))}{(\la-\la'\frac{a}{a'}g(N,N'))^2}d\s\\
&\ds +i\la'\int_{\s}e^{i\lambda u-i\la' u'}\frac{b\overline{b'}\frac{a}{a'}(g(\nabn N,N')+g(N,\nabn N'))}{(\la-\la'\frac{a}{a'}g(N,N'))^2}d\s,
\end{array}
\end{equation}
where we have used \eqref{du} to obtain the third term in the right-hand side of \eqref{bisof5}. Since $|\la'\frac{a}{a'}g(N,N')|<\la$, we may expand the fractions in \eqref{bisof5}:
\begin{equation}\label{bisof6}
\frac{1}{\la-\la'\frac{a}{a'}g(N,N')}=\sum_{p\geq 0}\frac{1}{\la}\left(\frac{\la'\frac{a}{a'}g(N,N')}{\la}\right)^p.
\end{equation}
and
\begin{equation}\label{bisof6bis}
\frac{1}{(\la-\la'\frac{a}{a'}g(N,N'))^2}=\sum_{p\geq 0}\frac{p+1}{\la^2}\left(\frac{\la'\frac{a}{a'}g(N,N')}{\la}\right)^p.
\end{equation}
For $p\in\mathbb{Z}$, We introduce the notation $F_{j,p}( u)$:
\begin{equation}\label{bisof7}
F_{j,p}( u)=\int_{0}^{+\infty}e^{i\lambda u}\psi(2^{-j}\lambda)f(\lambda\o)(2^{-j}\la)^{p}\lambda^2 d\lambda.
\end{equation}
Together with \eqref{bisof3}, \eqref{bisof5} and \eqref{bisof6}, this implies:
\begin{equation}\label{bisof8}
\int_{\s}U_jf(x)\overline{U_kf(x)}d\s =\sum_{p\geq 0}A^1_p+\sum_{p\geq 0}A^2_p+\sum_{p\geq 0}A^3_p+\sum_{p\geq 0}A^4_p,
\end{equation}
where $A^1_p$, $A^2_p$, $A^3_p$ and $A^4_p$ are given by:
\begin{equation}\label{bisof9}
\begin{array}{ll}
\ds A^1_p= & \ds 2^{-j-p(j-k)}\int_{\s}\left(\int_{\S}(\nabn b+ b\trt)a^{p+1}N^pF_{j,-p-1}( u)d\o\right)\\
& \ds\cdot\overline{\left(\int_{\S}b'{a'}^{-p}{N'}^pF_{k,p}( u')d\o'\right)}d\s ,
\end{array}
\end{equation}
\begin{equation}\label{bisof10}
\begin{array}{ll}
\ds A^2_p= & \ds 2^{-j-p(j-k)}\int_{\s}\left(\int_{\S}ba^{p+1}N^{p+1}F_{j,-p-1}( u)d\o\right)\\
& \ds\cdot\overline{\left(\int_{\S}\nabla b'{a'}^{-p}{N'}^pF_{k,p}( u')d\o'\right)}d\s .
\end{array}
\end{equation}
\begin{equation}\label{bisof9b}
\begin{array}{ll}
\ds A^3_p= & \ds (p+1)2^{-j-(p+1)(j-k)}\int_{\s}\left(\int_{\S} b(\nabn a N+a\nabn N)a^{p}N^pF_{j,-p-2}( u)d\o\right)\\
& \ds\cdot\overline{\left(\int_{\S}b'{a'}^{-p-1}{N'}^{p+1}F_{k,p+1}( u')d\o'\right)}d\s ,
\end{array}
\end{equation}
and
\begin{equation}\label{bisof10b}
\begin{array}{ll}
\ds A^4_p= & \ds (p+1)2^{-j-(p+1)(j-k)}\int_{\s}\left(\int_{\S}ba^{p+1}N^{p+2}F_{j,-p-2}( u)d\o\right)\\
& \ds\cdot\overline{\left(\int_{\S} b'(\nabla\log(a')N'+\nabla N'){a'}^{-p-1}{N'}^pF_{k,p+1}( u')d\o'\right)}d\s .
\end{array}
\end{equation}

\begin{remark}\label{bisrmksep}
The expansion \eqref{bisof6} allows us to rewrite $\int_{\s}U_jf(x)\overline{U_kf(x)}d\s $ in the form \eqref{bisof8}, i.e. as a sum of terms $A^1_p$, $A^2_p$, $A^3_p$, $A^4_p$. The key point is that in each of these terms - according to \eqref{bisof9}-\eqref{bisof10b} - one may separate the terms depending of $(\la,\o)$ from the terms depending on $(\la',\o')$. 
\end{remark}

\subsubsection{Estimates for $A^1_p$ and $A^2_p$} 

Let $H(x,\o)$ a tensor such that $\norm{H}_{\l{\infty}{2}}\lesssim D$. Then proceeding as in the basic computation \eqref{bisb5}, we have for any $p\in\mathbb{Z}$:
\begin{equation}\label{bisof11}
\begin{array}{ll}
\ds\normm{\int_{\S}H(x,\o)F_{j,p}( u)d\o}_{\ll{2}} &\ds\leq\int_{\S}\norm{H}_{\l{\infty}{2}}\norm{F_{j,p}( u)}_{L^2_{ u}}d\o\\
&\ds\leq\norm{H}_{\l{\infty}{2}}\norm{\psi(2^{-j}\lambda)f(\lambda\o)(2^{-j}\la)^{p}\la}_{\le{2}}\\
&\ds\lesssim D2^{|p|+j}\ga_j.
\end{array}
\end{equation}
where we have used the fact that $1/2\leq 2^{-j}\la\leq 2$ on the support of $\psi(2^{-j}\la)$. Now, {\bf Assumption 1} on the regularity of $a$, $N$, $\th$ and assumption \eqref{thregx1s} on the regularity of $b$ yield:
\begin{equation}\label{bisof12}
\begin{array}{rr}
\ds\norm{(\nabn b+ b\trt)a^{p+1}N^p}_{\l{\infty}{2}}+\norm{\nabla b'{a'}^{-p}{N'}^p}_{\lprime{\infty}{2}}&\\
\ds +\norm{b(\nabn a N+a\nabn N)a^{p}N^p}_{\l{\infty}{2}}&\\
\ds +\norm{b'(\nabla\log(a')N'+\nabla N'){a'}^{-p-1}{N'}^p}_{\lprime{\infty}{2}} & \ds\lesssim D,
\end{array}
\end{equation}
which together with \eqref{bisof11} implies:
\begin{equation}\label{bisof13b}
\begin{array}{ll}
&\ds\normm{\int_{\S}(\nabn b+ b\trt)a^{p+1}N^pF_{j,-p-1}( u)d\o}_{\ll{2}}\\
&\ds +\normm{\int_{\S}\nabla b'{a'}^{-p}{N'}^pF_{k,p}( u')d\o'}_{\ll{2}}\\
&\ds+\normm{\int_{\S} b(\nabn a N+a\nabn N)a^{p}N^pF_{j,-p-2}( u)d\o}_{\ll{2}}\\
&\ds +\normm{\int_{\S} b'(\nabla\log(a')N'+\nabla N'){a'}^{-p-1}{N'}^pF_{k,p+1}( u')d\o'}_{\ll{2}}\\ 
\ds\lesssim & \ds D2^{p+j}\ga_j.
\end{array}
\end{equation}

Note that Proposition \ref{bisorthoangle} together with Proposition \ref{bisdiagonal} yields the estimate:
\begin{equation}\label{bisof13}
\norm{U_jf}_{\ll{2}}\lesssim D\ga_j,
\end{equation}
for any symbol $b$ satisfying the assumptions \eqref{thregx1s} and \eqref{threomega1s}. Now, the terms containing no derivative in \eqref{bisof9}-\eqref{bisof10b} have a symbol given respectively by $b'{a'}^{-p}{N'}^p$, $ba^{p+1}N^{p+1}$, $b'{a'}^{-p-1}{N'}^{p+1}$ and $ba^{p+1}N^{p+2}$. These symbols satisfies the assumptions \eqref{thregx1s} and \eqref{threomega1s} since $b$ does, and since $a$, $N$, $\th$ satisfy {\bf Assumption 1} and {\bf Assumption 2}. Applying \eqref{bisof13}, we obtain: 
\begin{equation}\label{bisof14}
\ds\normm{\int_{\S}b'{a'}^{-p}{N'}^pF_{k,p}( u')d\o'}_{\ll{2}}+\normm{\int_{\S}b'{a'}^{-p-1}{N'}^{p+1}F_{k,p+1}( u')d\o'}_{\ll{2}}\lesssim D2^{p}\ga_k,
\end{equation}
and 
\begin{equation}\label{bisof15}
\ds\normm{\int_{\S}ba^{p+1}N^{p+1}F_{j,-p-1}( u)d\o}_{\ll{2}}+\normm{\int_{\S}ba^{p+1}N^{p+2}F_{j,-p-2}( u)d\o}_{\ll{2}}\lesssim D2^{p}\ga_j,
\end{equation}
where we have used the fact that $1/2\leq 2^{-j}\la\leq 2$ on the support of $\psi(2^{-j}\la)$.

Finally, the definition of $A_p^1-A_p^4$ given by \eqref{bisof9}-\eqref{bisof10b} and the estimates \eqref{bisof13b}, \eqref{bisof14} and \eqref{bisof15} yield:
\begin{equation}\label{bisof16}
|A^1_p|\lesssim D2^{2p-p(j-k)}\ga_j\ga_k,\,\forall p\geq 0,
\end{equation}
and
\begin{equation}\label{bisof17}
|A^2_p|+|A^3_p|+|A^4_p|\lesssim D2^{2p-(p+1)(j-k)}\ga_j\ga_k,\,\forall p\geq 0.
\end{equation}
\eqref{bisof16} and \eqref{bisof17} imply:
\bea\label{bisof18}
\nn\sum_{p\geq 1}|A^1_p|+\sum_{p\geq 0}(|A^2_p|+|A^3_p|+|A^4_p|)&\lesssim& D2^{-(j-k)}\left(\sum_{p\geq 0}2^{-p(j-k-2)}\right)\ga_j\ga_k\\
&\lesssim& D2^{-(j-k)}\ga_j\ga_k,
\eea
where we have used the assumption $j-k-2>0$. \eqref{bisof8} and \eqref{bisof18} will yield \eqref{bisof2} provided we obtain a similar estimate for $A^1_0$. Now, the estimate of $A^1_0$ provided by \eqref{bisof16} is not sufficient since it does not contain any decay in $j-k$. We will need to perform a second integration by parts for this term.

\subsubsection{A more precise estimate for $A^1_0$} From \eqref{bisof9} with $p=0$, we have:
\begin{equation}\label{bisof19}
\ds A^1_0=2^{-j}\int_{\s}\left(\int_{\S}(a\nabn b+ b\trt)F_{j,-1}( u)d\o\right)\overline{U_k(x)}.
\end{equation}
We decompose $\nabn b=b^j_1+b^j_2$ using the assumption \eqref{cordecfr1s}.
In turn, this yields a decomposition for $A^1_0$:
\begin{equation}\label{bisof21}
A^1_0=A^1_{0,1}+A^1_{0,2}
\end{equation}
where:
\begin{equation}\label{bisof22}
\begin{array}{l}
\ds A^1_{0,1}=2^{-j}\int_{\s}\left(\int_{\S}ab^j_1F_{j,0}( u)d\o\right)\overline{U_k(x)}d\s ,\\
\ds A^1_{0,2}=2^{-j}\int_{\s}\left(\int_{\S}(ab^j_2+ b\trt)F_{j,0}( u)d\o\right)\overline{U_k(x)}d\s .
\end{array}
\end{equation}

We first estimate $A^1_{0,1}$. We have: 
\begin{equation}\label{bisof23}
\begin{array}{ll}
\ds |A^1_{0,1}| & \ds\leq 2^{-j}\int_{\S}\left|\int_{\s}ab^j_1F_{j,0}( u)\overline{U_k(x)}d\s \right|d\o\\
& \ds\leq 2^{-j}\int_{\S}\norm{b^j_1}_{\ll{2}}\norm{a}_{\ll{\infty}}\norm{F_{j,0}}_{L^2_{ u}}\norm{U_k}_{\l{\infty}{2}}d\o\\
& \ds\lesssim D2^{-\frac{3j}{2}}\int_{\S}\norm{F_{j,0}}_{L^2_{ u}}\norm{U_k}_{\l{\infty}{2}}d\o,
\end{array}
\end{equation}
where we have used {\bf Assumption 1} on $a$ and the assumption \eqref{cordecfr1s} on $b^j_1$ in the last inequality. Plancherel yields:
\begin{equation}\label{bisof24}
\norm{F_{j,0}}_{L^2_{\o, u}}\leq\norm{\psi(2^{-j}\la)f(\la\o)\la}_{\le{2}}\lesssim 2^j\ga_j.
\end{equation}
In view of \eqref{bisof23}, we also need to estimate $\norm{U_k}_{\l{\infty}{2}}$. We have:
\begin{equation}\label{bisof25}
\norm{U_k}_{\l{\infty}{2}}\lesssim(\norm{\nabla U_k}_{\ll{2}}+\norm{U_k}_{\ll{2}})^{\frac{1}{2}}\norm{U_k}^{\frac{1}{2}}_{\ll{2}}\lesssim D\gamma_k+D^{\frac{1}{2}}\ga_k^{\frac{1}{2}}\norm{\nabla U_k}^{\frac{1}{2}}_{\ll{2}},
\end{equation}
where we have used the fact that $H^1(\s)$ embeds in $\l{\infty}{2}$ for the first inequality (see \cite{param1} Corollary 3.6 for a proof only using the regularity given by {\bf Assumption 1}), and \eqref{bisof13} for the second inequality. We still need to estimate $\norm{\nabla U_k}_{\ll{2}}$. We have:
\begin{equation}\label{bisof26}
\begin{array}{ll}
\ds\nabla U_k(x)= & \ds\int_{\S}\int_0^{+\infty} e^{i\la u}\nabla b\psi(2^{-k}\la)f(\la\o)\la^2d\la d\o\\
& \ds +i2^k\int_{\S}\int_0^{+\infty} e^{i\la u} ba^{-1}N\psi(2^{-k}\la)(2^{-k}\la)f(\la\o)\la^2d\la d\o.
\end{array}
\end{equation}
Using the basic computation \eqref{bisb5} for the first term together with the fact that $\nabla b\in\l{\infty}{2}$, and \eqref{bisof13} for the second term together with the fact that $ba^{-1}N$ satisfies the assumptions \eqref{thregx1s} and \eqref{threomega1s}, we obtain:
\begin{equation}\label{bisof27}
\norm{\nabla U_k}_{\ll{2}}\lesssim D2^k\ga_k.
\end{equation}
Finally, \eqref{bisof23}, \eqref{bisof24}, \eqref{bisof25} and \eqref{bisof27} yield:
\begin{equation}\label{bisof28}
|A^1_{0,1}|\lesssim D2^{-\frac{j-k}{2}}\ga_j\ga_k.
\end{equation}

\subsubsection{A second integration by parts} We now estimate the term $A^1_{0,2}$ defined in \eqref{bisof22}. We perform a second integration by parts relying again on \eqref{bisof4}. We obtain: 
\begin{equation}\label{bisof29}
\begin{array}{ll}
\ds A^1_{0,2}= &\ds 2^{-2j}\int_{\s}\left(\int_{\S}(\nabn b^j_2a+b^j_2\nabn a+b^j_2a\trt)F_{j,0}( u)d\o\right)\overline{U_k(x)}d\s \\
&\ds +2^{-2j}\int_{\s}\left(\int_{\S}b^j_2a^2N F_{j,0}( u)d\o\right)\cdot\overline{\nabla U_k(x)}d\s +\cdots,
\end{array}
\end{equation}
where we only mention the first term generated by the expansion \eqref{bisof6}. In fact, the other terms generated by \eqref{bisof6} and the ones generated by \eqref{bisof6bis} are estimated in the same way and generate more decay in $j-k$ similarly to the estimates \eqref{bisof16} \eqref{bisof17}.

The first term in the right-hand side of \eqref{bisof29} has the same form than $A^1_{0,1}$ defined in \eqref{bisof22} where $ab^j_1$ is replaced by $2^{-j}(\nabn b^j_2a+b^j_2\nabn a+ab^j_2\trt)$. By {\bf Assumption 1} and \eqref{cordecfr1s}, $2^{-j}(\nabn b^j_2a+b^j_2\nabn a+ab^j_2\trt)$ satisfies:
\begin{equation}
\begin{array}{ll}
&\ds\norm{2^{-j}(\nabn b^j_2a+b^j_2\nabn a+ab^j_2\trt)}_{\ll{2}}\\
\leq &\ds 2^{-j}\norm{\nabn b^j_2}_{\ll{2}}\norm{a}_{\ll{\infty}}+2^{-j}\norm{b^j_2}_{\l{2}{\infty}}\norm{\nabn a}_{\l{\infty}{2}}\\
&\ds +2^{-j}\norm{a}_{\ll{\infty}}\norm{b^j_2}_{\l{2}{\infty}}\norm{\trt}_{\l{\infty}{2}}\\
\lesssim &\ds D2^{-\frac{j}{2}}.
\end{array}
\end{equation}
Since $b^1_j$ and $2^{-j}(\nabn b^j_2a+b^j_2\nabn a+ab^j_2\trt)$ satisfy the same estimate, we obtain the analog of \eqref{bisof28} for the first term in the right-hand side of \eqref{bisof29}:
\begin{equation}\label{bisof30}
\ds \left|2^{-2j}\int_{\s}\left(\int_{\S}(\nabn b^j_2a+b^j_2\nabn a+ab^j_2\trt)F_{j,0}( u)d\o\right)\overline{U_k(x)}d\s \right|\lesssim D2^{-\frac{j-k}{2}}\ga_j\ga_k.
\end{equation}

We now estimate the second term in the right-hand side of \eqref{bisof29}. Using {\bf Assumption 1} on $a$ together with \eqref{cordecfr1s}, we have:
\begin{equation}\label{bisof31}
\norm{b^j_2a^2N}_{\l{\infty}{2}}\lesssim D.
\end{equation}
We estimate the second term in the right-hand side of \eqref{bisof29} using the assumption \eqref{thregx1s} on $b$, the basic computation \eqref{bisb5} and \eqref{bisof31}: 
\begin{equation}\label{bisof32}
\begin{array}{ll}
&\ds\left|2^{-2j}\int_{\s}\left(\int_{\S}b^j_2a^2N F_{j,0}( u)d\o\right)\cdot\overline{\left(\int_{\S}\nabla b'F_{k,0}( u')d\o'\right)}d\s \right|\\
\ds\leq &\ds 2^{-2j}\normm{\int_{\S}b^j_2a^2N F_{j,0}( u)d\o}_{\ll{2}}\normm{\int_{\S}\nabla b'F_{k,0}( u')d\o'}_{\ll{2}}\\
\ds\leq &\ds 2^{-2j}\left(\int_{\S}\norm{b^j_2a^2N}_{\l{\infty}{2}}\norm{F_{j,0}}_{L^2_{ u}}d\o\right)\left(\int_{\S}\norm{\nabla b}_{\l{\infty}{2}}\norm{F_{k,0}}_{L^2_{ u}}d\o\right)\\
\lesssim &\ds D^22^{-(j-k)}\ga_j\ga_k.
\end{array}
\end{equation}
Finally, \eqref{bisof29}, \eqref{bisof30} and \eqref{bisof32} imply:
\begin{equation}\label{bisof33}
|A^1_{0,2}|\lesssim D^22^{-\frac{j-k}{2}}\ga_j\ga_k.
\end{equation}

\subsubsection{End of the proof of Proposition \ref{bisorthofreq}} Since $A^1_0=A^1_1+A^1_2$, the estimate \eqref{bisof28} of $A^1_{0,1}$ and the estimate \eqref{bisof33} of $A^1_{0,2}$ yield:
\begin{equation}\label{bisof34}
|A^1_{0}|\lesssim D^22^{-\frac{j-k}{2}}\ga_j\ga_k.
\end{equation}
Together with \eqref{bisof8} and \eqref{bisof18}, this implies:
\begin{equation}\label{bisof35}
\left|\int_{\s}U_jf(x)\overline{U_kf(x)}d\s \right|\lesssim D^22^{-\frac{|j-k|}{2}}\ga_j\ga_k\textrm{ for }|j-k|>2.
\end{equation}
Finally, \eqref{bisof35} together with Shur's Lemma yields:
\begin{equation}\label{bisof36}
\norm{Uf}_{\ll{2}}^2\lesssim\sum_{j\geq -1}\norm{U_jf}_{\ll{2}}^2+D^2\norm{f}^2_{\le{2}}.
\end{equation}
This concludes the proof of Proposition \ref{bisorthofreq}. \QED

\subsection{Proof of Proposition \ref{bisorthoangle} (almost orthogonality in angle)}\label{bissec:orthoangle}

We have to prove \eqref{bisorthoangle1}:
\begin{equation}\label{bisoa1}
\norm{U_jf}_{\ll{2}}^2\lesssim\sum_{\nu\in\Gamma}\norm{U^\nu_jf}_{\ll{2}}^2+D^2\ga_j^2.
\end{equation}
This will result from the following inequality:
\begin{equation}\label{bisoa2}
\left|\int_{\s}U^\nu_jf(x)\overline{U^{\nu'}f(x)}d\s \right|\lesssim \frac{D^2\ga_j^\nu\ga_j^{\nu'}}{2^{j\a/2}(2^{j/2}|\nu-\nu'|)^{2-\a}}+\frac{D^2\ga_j^\nu\ga_j^{\nu'}}{(2^{j/2}|\nu-\nu'|)^{3}},\,|\nu-\nu'|\neq 0,
\end{equation}
where $\a>0$. Indeed, since $\S$ is 2 dimensional and $1\leq 2^{j/2}|\nu-\nu'|\leq 2^{j/2}$ for $\nu, \nu'\in\Gamma$ and $\nu\neq\nu'$, we have:
\begin{equation}\label{bisoa3}
\sup_{\nu}\sum_{\nu'} \frac{1}{(2^{j/2}|\nu-\nu'|)^{3}}\leq C<+\infty,
\end{equation}
and 
\begin{equation}\label{bisoa3b}
\sup_{\nu}\sum_{\nu'} \frac{1}{2^{j\a/2}(2^{j/2}|\nu-\nu'|)^{2-\a}}\leq C_\a<+\infty\,\forall\a>0.
\end{equation}
Thus, \eqref{bisoa2}, \eqref{bisoa3} and \eqref{bisoa3b} together with Shur's Lemma imply \eqref{bisoa1}.

\begin{remark}
In \cite{SmTa}, the authors rely on a partial Fourier transform with respect to a coordinate system on $\p$ to prove almost orthogonality in angle for their parametrix. In our case, coordinate systems on $\p$ are not regular enough, which forces us to work invariantly. More precisely, we will use geometric integrations by parts in tangential directions to $\p$ in order to obtain \eqref{bisoa2}. 
\end{remark}

\subsubsection{A second decomposition in frequency} From now on, we focus on proving \eqref{bisoa2}. Integrating by parts twice in $\int_{\s}U^\nu_jf(x)\overline{U^{\nu'}_jf(x)}d\s $ would ultimately yield:
\begin{equation}\label{bisoa4}
\left|\int_{\s}U^\nu_jf(x)\overline{U^{\nu'}_jf(x)}d\s \right|\lesssim \frac{D^2\ga_j^\nu\ga_j^{\nu'}}{(2^{j/2}|\nu-\nu'|)^{2}},\,|\nu-\nu'|\neq 0.
\end{equation}
This corresponds to the case $\a=0$ in \eqref{bisoa3} and yields to a log-loss since we have:
\begin{equation}\label{bisoa5}
\sup_{\nu}\sum_{\nu'} \frac{1}{(2^{j/2}|\nu-\nu'|)^{2}}\sim j.
\end{equation}
To avoid this log-loss, we do a second decomposition in frequency. $\la$ belongs to the interval $[2^{j-1},2^{j+1}]$ which we decompose in intervals $I_k$:
\begin{equation}\label{bisoa6}
[2^{j-1},2^{j+1}]=\bigcup_{1\leq k\leq |\nu-\nu'|^{-\a}}I_k\textrm{ where }\textrm{diam}(I_k)\sim 2^j|\nu-\nu'|^\a.
\end{equation}
Let $\phi_k$ a partition of unity of the interval $[2^{j-1},2^{j+1}]$ associated to the $I_k$'s. We decompose $U^\nu_jf$ as follows:
\begin{equation}\label{bisoa7}
U^\nu_jf(x)=\sum_{1\leq k\leq |\nu-\nu'|^{-\a}}U^{\nu,k}_jf(x),
\end{equation}
where:
\begin{equation}\label{bisoa8}
U^{\nu,k}_jf(x)=\int_{\S}\int_{0}^{+\infty}e^{i\lambda u}b(x,\o)\psi(2^{-j}\lambda)\eta^\nu_j(\o)\phi_k(\la)f(\lambda\o)\lambda^2 d\lambda d\o.
\end{equation}
We also define:
\begin{equation}\label{bisoa9}
\begin{array}{l}
\ds\ga^{\nu,k}_j=\norm{\psi(2^{-j}\la)\eta^\nu_j(\o)\phi_k(\la)f}_{\le{2}},\,j\geq 0,\,\nu\in\Gamma,1\leq k\leq |\nu-\nu'|^{-\a}, 
\end{array}
\end{equation}
which satisfy:
\begin{equation}\label{bisoa10}
(\gamma^{\nu}_j)^2=\sum_{1\leq k\leq |\nu-\nu'|^{-\a}}(\ga^{\nu,k}_j)^2.
\end{equation}

\subsubsection{The two key estimates}\label{bissec:keyest} 

We will prove the following two estimates:
\begin{equation}\label{bisoa11}
\begin{array}{r}
\ds\left|\int_{\s}U^{\nu,k}_jf(x)\overline{U^{\nu',k}_jf(x)}d\s \right|\lesssim \frac{D^2\ga_j^{\nu,k}\ga_j^{\nu',k}}{2^{j\a/2}(2^{j/2}|\nu-\nu'|)^{2-\a}}+\frac{D^2\ga_j^{\nu,k}\ga_j^{\nu',k}}{(2^{j/2}|\nu-\nu'|)^{3}}\\
\ds\textrm{for }|\nu-\nu'|\neq 0,\,1\leq k\leq |\nu-\nu'|^{-\a},
\end{array}
\end{equation}
and
\begin{equation}\label{bisoa12}
\begin{array}{r}
\ds\left|\int_{\s}U^{\nu,k}_jf(x)\overline{U^{\nu',k'}_jf(x)}d\s \right|\lesssim \frac{D^2\ga_j^{\nu,k}\ga_j^{\nu',k'}}{|k-k'|2^{j/2(1-4\a)}(2^{j/2}|\nu-\nu'|)^{1+4\a}},\\
\ds\textrm{ for }|\nu-\nu'|\neq 0,\,1\leq k,k'\leq |\nu-\nu'|^{-\a},\,k\neq k'.
\end{array}
\end{equation}
\eqref{bisoa11} and \eqref{bisoa12} imply:
\begin{equation}\label{bisoa13}
\begin{array}{lll}
\ds\left|\int_{\s}U^{\nu}_jf(x)\overline{U^{\nu'}_jf(x)}d\s \right| & \ds\leq & \ds\sum_{1\leq k\leq |\nu-\nu'|^{-\a}}\left|\int_{\s}U^{\nu,k}_jf(x)\overline{U^{\nu',k}_jf(x)}d\s \right|\\
& & \ds +\sum_{1\leq k\neq k'\leq |\nu-\nu'|^{-\a}}\left|\int_{\s}U^{\nu,k}_jf(x)\overline{U^{\nu',k'}_jf(x)}d\s \right|\\
& \ds\lesssim & \ds\sum_{1\leq k\leq |\nu-\nu'|^{-\a}}\frac{D^2\ga_j^{\nu,k}\ga_j^{\nu',k}}{2^{j\a/2}(2^{j/2}|\nu-\nu'|)^{2-\a}}\\
& & \ds +\sum_{1\leq k\leq |\nu-\nu'|^{-\a}}\frac{D^2\ga_j^{\nu,k}\ga_j^{\nu',k}}{(2^{j/2}|\nu-\nu'|)^{3}}\\
& &\ds +\sum_{1\leq k\neq k'\leq |\nu-\nu'|^{-\a}}\frac{D^2\ga_j^{\nu,k}\ga_j^{\nu',k'}}{|k-k'|2^{\frac{j}{2}(1-4\a)}(2^{j/2}|\nu-\nu'|)^{1+4\a}}\\
& \lesssim & \ds\frac{D^2\ga_j^{\nu}\ga_j^{\nu'}}{2^{j\a/2}(2^{j/2}|\nu-\nu'|)^{2-\a}}+\frac{D^2\ga_j^{\nu}\ga_j^{\nu'}}{(2^{j/2}|\nu-\nu'|)^{3}},
\end{array}
\end{equation}
where we have used \eqref{bisoa10} and the fact that we may choose $0<\a<1/5$, together with the fact that:
\begin{equation}\label{bisoa14}
\sup_{1\leq k\leq |\nu-\nu'|^{-\a}}\sum_{1\leq k'\leq |\nu-\nu'|^{-\a},\,k'\neq k}\frac{1}{|k-k'|}\lesssim \a |\log(|\nu-\nu'|)|.
\end{equation}
Since \eqref{bisoa13} yields the wanted estimate \eqref{bisoa2}, we are left with proving \eqref{bisoa11} and \eqref{bisoa12}.

\subsubsection{Proof of \eqref{bisoa11}} The estimate \eqref{bisoa11} will result of two integrations by parts with respect to tangential derivatives. By definition of $\nabb$, we have $\nabb h=\nabla h-(\nabn h)N$ for any function $h$ on $\s$. In particular, we have $\nabb( u)=0$ and $\nabb( u')={a'}^{-1}N'-{a'}^{-1}g(N',N) N$. Now, since $|N'-(N'\cdot N)N|^2=1-(N'\cdot N)^2$, this yields:
\begin{equation}\label{bisoa15} 
e^{i\la u-i\la' u'}=\frac{i}{\la'(1-(N'\cdot N)^2)}\nabla_{N'-\gn N}(e^{i\la u-i\la' u'}),
\end{equation}
where we have used the fact that $N'-(N'\cdot N) N$ is a tangent vector with respect of the level surfaces of $ u$. Similarly, we have:
\begin{equation}\label{bisoa17} 
e^{i\la u-i\la' u'}=-\frac{i}{\la(1-(N'\cdot N)^2)}\nabla_{N-\gn N'}(e^{i\la u-i\la' u'}),
\end{equation}
where we have used the fact that $N-(N\cdot N') N'$ is a tangent vector with respect of the level surfaces of $ u'$. For $p\in\mathbb{Z}$, We introduce the notation $F_{j,k,p}( u)$:
\begin{equation}\label{bisoa19}
F_{j,k,p}( u)=\int_{0}^{+\infty}e^{i\lambda u}\psi(2^{-j}\lambda)\phi_k(\la)f(\lambda\o)(2^{-j}\la)^{p}\lambda^2 d\lambda.
\end{equation}
We integrate once by parts using \eqref{bisoa15} in $\int_{\s}U^{\nu,k}_jf(x)\overline{U^{\nu',k}_jf(x)}d\s $ and we obtain:
\begin{equation}\label{bisoa20}
\begin{array}{ll}
&\ds\int_{\s}U^{\nu,k}_jf(x)\overline{U^{\nu',k}_jf(x)}d\s \\
\ds = &\ds i2^{-j}\int_{\s\times\S\times\S}\textrm{div}\left(\frac{(N'-(N\cdot N')N)a'b\overline{b'}}{1-(N\cdot N')^2}\right)\\
&\ds \times F_{j,k,0}( u)\overline{F_{j,k,-1}( u')}\eta^\nu_j(\o)\eta^{\nu'}_j(\o')d\o d\o'd\s.
\end{array}
\end{equation}
We then integrate a second time by parts using \eqref{bisoa17} (so that there is at least one tangential derivative for each quantity where two derivatives are taken):
\begin{equation}\label{bisoa21}
\begin{array}{ll}
&\ds\int_{\s}U^{\nu,k}_jf(x)\overline{U^{\nu',k}_jf(x)}d\s \\
\ds = & \ds 2^{-2j}\int_{\s\times\S\times\S}
\textrm{div}\left(\frac{(N-(N\cdot N')N')a}{1-(N\cdot N')^2}\textrm{div}\left(\frac{(N'-(N\cdot N')N)a'b\overline{b'}}{1-(N\cdot N')^2}\right)\right)\\
&\ds\times F_{j,k,-1}( u)\overline{F_{j,k,-1}( u')}\eta^\nu_j(\o)\eta^{\nu'}_j(\o')d\o d\o'd\s.
\end{array}
\end{equation}

\vspace{0.2cm}

\noindent{\bf Computation of the right-hand side of \eqref{bisoa21}.} We would like to compute the double divergence term in the right-hand side of \eqref{bisoa21}. This is done in the following lemma.
\begin{lemma}\label{doublediv}
The double divergence term in the right-hand side of \eqref{bisoa21} is given by:
\begin{equation}\label{dd1} 
\begin{array}{ll}
&\ds\textrm{div}\left(\frac{(N-(N\cdot N')N')a}{1-(N\cdot N')^2}\textrm{div}\left(\frac{(N'-(N\cdot N')N)a'b\overline{b'}}{1-(N\cdot N')^2}\right)\right)\\
\ds = &\ds\frac{1}{|N_\nu-N_{\nu'}|^2}\left(\sum_{p,q\geq 0}c_{p,q}\left(\frac{N-N_\nu}{|N_\nu-N_{\nu'}|}\right)^p\left(\frac{N'-N_{\nu'}}{|N_\nu-N_{\nu'}|}\right)^q\right)F,
\end{array}
\end{equation}
where $F$ is a combination of terms in the following list:
\begin{equation}\label{dd2} 
\begin{array}{l}
\ds \frac{(\nabla\th-\nabla\th')aa'bb'}{|N_\nu-N_{\nu'}|},\,\frac{(\th-\th')\nabla(ab)a'b'}{|N_\nu-N_{\nu'}|},\,\th\nabla(ab)a'b',\,ab\th\nabla(a'b'),\,\nabb\nabla(ab)a'b',\\
\ds\nabla(a)\nabla(b)a'b',\,\nabla(ab)\nabla(a'b'),\,\frac{(\th-\th')^2aa'bb'}{|N_\nu-N_{\nu'}|^2},\,\th\th'aa'bb',\,\th^2aa'bb'.
\end{array}
\end{equation}
\end{lemma}
The proof of Lemma \ref{doublediv} is postponed to the Appendix A. The following lemma gives the structure of the terms in the list \eqref{dd2}.
\begin{lemma}\label{dd3}
The terms in the list \eqref{dd2} have the following form:
\begin{equation}\label{dd4}
H_1(x,\o,\nu,\nu')H_2(x,\o',\nu,\nu')+\frac{H_3(x,\o,\nu,\nu')H_4(x,\o',\nu,\nu')}{2^{j/2}|\nu-\nu'|},
\end{equation}
where $H_1, H_3, H_4$ satisfy:
\begin{equation}\label{dd5}
\norm{H_1}_{\ll{2}}+\norm{H_3}_{\l{\infty}{2}}+\norm{H_4}_{\lprime{\infty}{2}}\lesssim D,
\end{equation}
and where $H_2$ satisfies:
\begin{equation}\label{dd6}
\norm{H_2}_{\ll{\infty}}+\norm{\po H_2}_{\ll{2}}+\norm{\nabla\po H_2}_{\ll{2}}\lesssim D,
\end{equation}
for $\o$ in the support of $\eta_j^\nu$ and $\o'$ in the support of $\eta_j^{\nu'}$.
\end{lemma}
The proof of Lemma \ref{dd3} is postponed to the Appendix B. In the rest of this section, we show how Lemma \ref{doublediv} and Lemma \ref{dd3} yield the proof of \eqref{bisoa11}. 

\vspace{0.2cm}

\noindent{\bf End of the proof of \eqref{bisoa11}.} Using \eqref{bisoa21}, Lemma \ref{doublediv} and Lemma \ref{dd3}, we may rewrite $\int_{\s}U^{\nu,k}_jf(x)\overline{U^{\nu',k}_jf(x)}d\s$ as:
\bea\label{dd15}
&&\int_{\s}U^{\nu,k}_jf(x)\overline{U^{\nu',k}_jf(x)}d\s \\
\nn &=&\sum_{p,q\geq 0}c_{p,q}\int_{\s}\frac{2^{-j}}{(2^{j/2}|N_\nu-N_{\nu'}|)^2}\left(\int_{\S}\left(\frac{N-N_\nu}{|N_\nu-N_{\nu'}|}\right)^pH_1(x,\o,\nu,\nu')F_{j,k,-1}( u)\eta^\nu_j(\o)d\o\right)\\
\nn&&\times\left(\int_{\S}\left(\frac{N'-N_{\nu'}}{|N_\nu-N_{\nu'}|}\right)^qH_2(x,\o',\nu,\nu')\overline{F_{j,k,-1}( u')}\eta^{\nu'}_j(\o')d\o'\right)d\s\\
\nn&& +\sum_{p,q\geq 0}c_{p,q}\int_{\s}\frac{2^{-j}}{(2^{j/2}|N_\nu-N_{\nu'}|)^3}\left(\int_{\S}\left(\frac{N-N_\nu}{|N_\nu-N_{\nu'}|}\right)^pH_3(x,\o,\nu,\nu')F_{j,k,-1}( u)\eta^\nu_j(\o)d\o\right)\\
\nn&&\times\left(\int_{\S}\left(\frac{N'-N_{\nu'}}{|N_\nu-N_{\nu'}|}\right)^qH_4(x,\o',\nu,\nu')\overline{F_{j,k,-1}( u')}\eta^{\nu'}_j(\o')d\o'\right)d\s.
\eea
We estimate the two terms in the right-hand side of \eqref{dd15} starting with the second one. We have:
\bea\label{dd16}
\nn&&\ds \bigg|\int_{\s}\frac{2^{-j}}{(2^{j/2}|N_\nu-N_{\nu'}|)^3}\left(\int_{\S}\left(\frac{N-N_\nu}{|N_\nu-N_{\nu'}|}\right)^pH_3(x,\o,\nu,\nu')F_{j,k,-1}( u)\eta^\nu_j(\o)d\o\right)\\
\nn&&\ds\times\left(\int_{\S}\left(\frac{N'-N_{\nu'}}{|N_\nu-N_{\nu'}|}\right)^qH_4(x,\o',\nu,\nu')\overline{F_{j,k,-1}( u')}\eta^{\nu'}_j(\o')d\o'\right)d\s\bigg|\\
\nn&\leq &\ds\frac{2^{-j}}{(2^{j/2}|\nu-\nu'|)^{3+p+q}}\left(\int_{\S}\norm{H_3}_{\l{\infty}{2}}\norm{F_{j,k,-1}}_{L^2_u}\eta^\nu_j(\o)d\o\right)\\
\nn&&\ds\times\left(\int_{\S}\norm{H_4}_{\lprime{\infty}{2}}\norm{F_{j,k,-1}}_{L^2_{u'}}\eta^{\nu'}_j(\o')d\o'\right)\\
&\lesssim& \frac{D^2\gamma_j^{\nu,k}\gamma_j^{\nu',k}}{(2^{j/2}|\nu-\nu'|)^{3+p+q}},
\eea
where we have used {\bf Assumption 2} to estimate $|N_\nu-N_{\nu'}|$, Plancherel to estimate $\norm{F_{j,k,-1}}_{L^2_{ u}}$ and $\norm{F_{j,k,-1}}_{L^2_{ u'}}$, Cauchy-Schwartz in $\o$ and $\o'$, and the estimate \eqref{dd5} for $H_3$ and $H_4$.

We now estimate the second term in the right-hand side of \eqref{dd15}. We have:
\bea\label{dd17}
\nn&&\ds\bigg|\int_{\s}\frac{2^{-j}}{(2^{j/2}|N_\nu-N_{\nu'}|)^2}\left(\int_{\S}\left(\frac{N-N_\nu}{|N_\nu-N_{\nu'}|}\right)^pH_1(x,\o,\nu,\nu')F_{j,k,-1}( u)\eta^\nu_j(\o)d\o\right)\\
\nn&&\ds\times\left(\int_{\S}\left(\frac{N'-N_{\nu'}}{|N_\nu-N_{\nu'}|}\right)^qH_2(x,\o',\nu,\nu')\overline{F_{j,k,-1}( u')}\eta^{\nu'}_j(\o')d\o'\right)d\s\bigg|\\
\nn&\leq &\ds\frac{2^{-j}}{(2^{j/2}|\nu-\nu'|)^{2+p+q}}\left(\int_{\S}\norm{H_1}_{\ll{2}}\norm{F_{j,k,-1}}_{L^{\infty}_u}\eta^\nu_j(\o)d\o\right)\\
\nn&&\ds\times\normm{\int_{\S}(2^{j/2}(N'-N_{\nu'}))^qH_2(x,\o',\nu,\nu')F_{j,k,-1}( u')\eta^{\nu'}_j(\o')d\o'}_{\ll{2}}\\
\nn&\lesssim &\ds\frac{D|\nu-\nu'|^\a\gamma^{\nu,k}_j}{(2^{j/2}|\nu-\nu'|)^{2+p+q}}\\
&&\ds\times\normm{\int_{\S}(2^{j/2}(N'-N_{\nu'}))^qH_2(x,\o',\nu,\nu')F_{j,k,-1}( u')\eta^{\nu'}_j(\o')d\o'}_{\ll{2}},
\eea
where we have used {\bf Assumption 2} to estimate $|N_\nu-N_{\nu'}|$, Cauchy-Schwartz in $\o$, the estimate \eqref{dd5} for $H_1$ and the following estimate for $\norm{F_{j,k,-1}}_{L^\infty_{ u}}$:
\begin{equation}\label{bisoa29}
\norm{F_{j,k,-1}}_{L^\infty_{ u}}\lesssim 2^{3j/2}|\nu-\nu'|^{\frac{\a}{2}}\norm{\psi(2^{-j}\la)\phi_k(\la)f(\la\o)\la}_{L^2_\la},
\end{equation}
which follows form taking Cauchy Schwartz in $\la$ together with the size of the support of $\phi_k$.  Note that the symbol $F=(2^{j/2}(N'-N_{\nu'}))^qH_2(x,\o',\nu,\nu')$ satisfies the following assumptions:
\begin{equation}\label{dd18}
\norm{F}_{\ll{\infty}}\lesssim D,\,\norm{\po F}_{\ll{2}}\lesssim q D,\,\norm{\nabla\po F}_{\ll{2}}\lesssim q^2 D,
\end{equation}
where we have used {\bf Assumption 2} for $\po N$ and $\po^2N$, and the assumption \eqref{dd6} satisfied by $H_2$.
We will see in section \ref{bissec:diagonal} that assumptions \eqref{dd18} on a symbol is enough to control the diagonal term in $\ll{2}$ (i.e. to obtain the estimate \eqref{bisdiagonal1}). Thus, we obtain:
\begin{equation}\label{dd19}
\normm{\int_{\S}(2^{j/2}(N'-N_{\nu'}))^qH_2(x,\o',\nu,\nu')F_{j,k,-1}( u')\eta^{\nu'}_j(\o')d\o'}_{\ll{2}}\lesssim (1+q^2)D\gamma_j^{\nu',k}.
\end{equation}

\eqref{dd17} and \eqref{dd19} imply:
\bea\label{bisoa30}
\nn&&\ds\bigg|\int_{\s}\frac{2^{-j}}{(2^{j/2}|N_\nu-N_{\nu'}|)^2}\left(\int_{\S}\left(\frac{N-N_\nu}{|N_\nu-N_{\nu'}|}\right)^pH_1(x,\o,\nu,\nu')F_{j,k,-1}( u)\eta^\nu_j(\o)d\o\right)\\
\nn&&\ds\times\left(\int_{\S}\left(\frac{N'-N_{\nu'}}{|N_\nu-N_{\nu'}|}\right)^qH_2(x,\o',\nu,\nu')\overline{F_{j,k,-1}( u')}\eta^{\nu'}_j(\o')d\o'\right)d\s\bigg|\\
&\lesssim &\ds\frac{(1+q^2)D^2|\nu-\nu'|^{\frac{\a}{2}}}{(2^{j/2}|\nu-\nu'|)^{2+p+q}}\ga_j^{\nu,k}\ga_j^{\nu',k}.
\eea
 
Finally, \eqref{dd15}, \eqref{dd16} and \eqref{bisoa30} yield:
\begin{equation}\label{dd20}
\begin{array}{ll}
&\ds\left|\int_{\s}U^{\nu,k}_jf(x)\overline{U^{\nu',k}_jf(x)}d\s\right| \\
\ds \lesssim &\ds\sum_{p,q\geq 0}c_{p,q}\frac{(1+q)D^2|\nu-\nu'|^{\frac{\a}{2}}}{(2^{j/2}|\nu-\nu'|)^{2+p+q}}\ga_j^{\nu,k}\ga_j^{\nu',k}+\sum_{p,q\geq 0}
\frac{D^2\gamma_j^{\nu,k}\gamma_j^{\nu',k}}{(2^{j/2}|\nu-\nu'|)^{3+p+q}}\\
\ds\lesssim &\ds\frac{D^2|\nu-\nu'|^{\frac{\a}{2}}}{(2^{j/2}|\nu-\nu'|)^{2}}\ga_j^{\nu,k}\ga_j^{\nu',k}+\frac{D^2\gamma_j^{\nu,k}\gamma_j^{\nu',k}}{(2^{j/2}|\nu-\nu'|)^{3}},
\end{array}
\end{equation}
which concludes the proof of estimate \eqref{bisoa11}.

\subsubsection{Proof of \eqref{bisoa12}} The estimate \eqref{bisoa12} will result of two integrations by parts, one with respect to the normal derivative, and one with respect to tangential derivatives. We have:
\begin{equation}\label{bisoa43b}
e^{i\la u-i\la' u'}=-\frac{ia}{\la-\la'\frac{a}{a'}g(N,N')}\nabn(e^{i\la u-i\la' u'}).
\end{equation}
We integrate once by parts using \eqref{bisoa43b} in $\int_{\s}U^{\nu,k}_jf(x)\overline{U^{\nu',k'}_jf(x)}d\s $. We obtain:
\begin{equation}\label{bisoa43}
\begin{array}{l}
\ds \int_{\s}U^{\nu,k}_jf(x)\overline{U^{\nu',k'}_jf(x)}d\s  = i\int_{\s\times\S\times\S}\int_0^{+\infty}\int_0^{+\infty}\\
\ds \times\textrm{div}\left(\frac{aNbb'}{\la-\la'\frac{a}{a'}g(N,N')}\right)\eta_j^\nu(\o)\eta_j^{\nu'}(\o')\psi(2^{-j}\la)\psi(2^{-j}\la')\phi_k(\la)\phi_{k'}(\la') \\
\ds \times f(\la\o)f(\la'\o')\la^2 {\la'}^2d\la d\la' d\o d\o'd\s.
\end{array}
\end{equation}
We then expand the divergence term in the right-hand side of \eqref{bisoa43}:
\begin{equation}\label{bisoa43bb}
\ds\textrm{div}\left(\frac{aNbb'}{\la-\la'\frac{a}{a'}g(N,N')}\right)\\
\ds =D_1+D_2,
\end{equation}
where $D_1$ and $D_2$ are given by:
\begin{equation}\label{bisoa43bb1}
\begin{array}{ll}
\ds D_1= &\ds\frac{abb'\textrm{div}(N)+\nabn(ab)b'}{\la-\la'\frac{a}{a'}g(N,N')}\\
&\ds +\la'\frac{\nabn(a)a{a'}^{-1}bb'g(N,N')+\nabn(g(N,N'))a^2{a'}^{-1}bb'}{(\la-\la'\frac{a}{a'}g(N,N'))^2}
\end{array}
\end{equation}
and
\begin{equation}\label{bisoa43bb2}
D_2=\frac{ab\nabn(b')}{\la-\la'\frac{a}{a'}g(N,N')}-\la'\frac{\nabn(a')a^2{a'}^{-2}bb'g(N,N')}{(\la-\la'\frac{a}{a'}g(N,N'))^2}.
\end{equation}
We then integrate a second time by parts using \eqref{bisoa15} for $D_1$ and using \eqref{bisoa17} for $D_2$ (so that there is at least one tangential derivative on $a$, $a'$, $b$, $b'$ when two derivatives are taken). We obtain:
\begin{equation}\label{bisoa44}
\begin{array}{l}
\ds \int_{\s}U^{\nu,k}_jf(x)\overline{U^{\nu',k'}_jf(x)}d\s  = \int_{\s\times\S\times\S}\int_0^{+\infty}\int_0^{+\infty}\frac{1}{\la'}\\
\ds \times\textrm{div}\left(\frac{(N'-g(N,N')N)a'}{1-g(N,N')^2}D_1\right)\eta_j^\nu(\o)\eta_j^{\nu'}(\o')\psi(2^{-j}\la)\psi(2^{-j}\la')\phi_k(\la)\phi_{k'}(\la') \\
\ds \times f(\la\o)f(\la'\o')\la^2 {\la'}^2d\la d\la' d\o d\o'd\s+\int_{\s\times\S\times\S}\int_0^{+\infty}\int_0^{+\infty}\frac{1}{\la}\\
\ds \times\textrm{div}\left(\frac{(N-g(N,N')N')a}{1-g(N,N')^2}D_2\right)\eta_j^\nu(\o)\eta_j^{\nu'}(\o')\psi(2^{-j}\la)\psi(2^{-j}\la')\phi_k(\la)\phi_{k'}(\la') \\
\ds \times f(\la\o)f(\la'\o')\la^2 {\la'}^2d\la d\la' d\o d\o'd\s.
\end{array}
\end{equation}

\vspace{0.2cm}

\noindent{\bf Computation of the right-hand side of \eqref{bisoa44}.} We would like to compute the two divergence term in the right-hand side of \eqref{bisoa44}. This is done in the following lemma.

\begin{lemma}\label{twodiv}
The two divergence term in the right-hand side of \eqref{bisoa44} have the following form:
\begin{equation}\label{td1}
\begin{array}{l}
\ds \frac{1}{2^j(k-k')|N_\nu-N_{\nu'}||\nu-\nu'|^\a}\\
\times\ds\Bigg(\sum_{l,m,n,o,p,q\geq 0}\left(c^1_{p,q,l,m,n,o}F_1+
\frac{c^2_{p,q,l,m,n,o}F_2}{(k-k')|\nu-\nu'|^\a}+\frac{c^3_{p,q,l,m,n,o}F_3}{(k-k')^2|\nu-\nu'|^{2\a}}\right)\\
\times\ds\left(\frac{N-N_\nu}{|N_\nu-N_{\nu'}|}\right)^p\left(\frac{N'-N_{\nu'}}{|N_\nu-N_{\nu'}|}\right)^q\left(\frac{a-a_\nu}{(k-k')|\nu-\nu'|^\a}\right)^l\left(\frac{a'-a_{\nu'}}{(k-k')|\nu-\nu'|^\a}\right)^m\\
\ds\times\left(\frac{\la-k2^j|\nu-\nu'|^\a}{2^j(k-k')|\nu-\nu'|^\a}\right)^n\left(\frac{\la'-k'2^j|\nu-\nu'|^\a}{2^j(k-k')|\nu-\nu'|^\a}\right)^o\Bigg),

\end{array}
\end{equation}
where $F_j, j=1, 2, 3$ is a combination of terms in the following list:
\begin{equation}\label{td2} 
\begin{array}{l}
\ds \frac{(\th-\th')\th aa'bb'}{|N_\nu-N_{\nu'}|},\,\frac{(\th-\th')\nabla(ab)a'b'}{|N_\nu-N_{\nu'}|},\,\th\nabla(ab)a'b',\,ab\th\nabla(a'b'),\,\nabb\nabla(ab)a'b',\\
\ds\nabla(a)\nabla(b)a'b',\,\nabla(ab)\nabla(a'b'),\,\nabla(\th) aa'bb' ,\,\th\th'aa'bb',\,\th^2aa'bb'.
\end{array}
\end{equation}
\end{lemma}
The proof of Lemma \ref{twodiv} is postponed to the Appendix C. The following lemma gives the structure of the terms in the list \eqref{td2}.
\begin{lemma}\label{td3}
The terms in the list \eqref{td2} have the following form:
\begin{equation}\label{td4}
H_1(x,\o,\nu,\nu')H_2(x,\o',\nu,\nu')+H_3(x,\o,\nu,\nu')H_4(x,\o',\nu,\nu'),
\end{equation}
where $H_1, H_3, H_4$ satisfy:
\begin{equation}\label{td5}
\norm{H_1}_{\ll{2}}+\norm{H_3}_{\l{\infty}{2}}+\norm{H_4}_{\lprime{\infty}{2}}\lesssim D,
\end{equation}
and where $H_2$ satisfies:
\begin{equation}\label{td6}
\norm{H_2}_{\ll{\infty}}+\norm{\po H_2}_{\ll{2}}+\norm{\nabla\po H_2}_{\ll{2}}\lesssim D,
\end{equation}
for $\o$ in the support of $\eta_j^\nu$ and $\o'$ in the support of $\eta_j^{\nu'}$.
\end{lemma}
The proof of Lemma \ref{td3} follows the same line as the proof of Lemma \ref{dd3} and is left to the 
reader. In the rest of this section, we show how Lemma \ref{twodiv} and Lemma \ref{td3} yield the proof of \eqref{bisoa12}.

\vspace{0.2cm}

\noindent{\bf End of the proof of \eqref{bisoa12}.} Using \eqref{bisoa44}, Lemma \ref{twodiv} and Lemma \ref{td3}, we may rewrite $\int_{\s}U^{\nu,k}_jf(x)\overline{U^{\nu',k}_jf(x)}d\s$ as:
\bea\label{td15}
&&\ds\int_{\s}U^{\nu,k}_jf(x)\overline{U^{\nu',k}_jf(x)}d\s \\
\nn&= &\ds\sum_{l,m,n,o,p,q\geq 0}\left(\frac{c^1_{p,q,l,m,n,o}}{(k-k')|\nu-\nu'|^\a}+
\frac{c^2_{p,q,l,m,n,o}}{(k-k')^2|\nu-\nu'|^{2\a}}+\frac{c^3_{p,q,l,m,n,o}}{(k-k')^3|\nu-\nu'|^{3\a}}\right)\\
\nn&&\times\frac{1}{(k-k')^{l+m+n+o}}\int_{\s}\frac{2^{-3j/2}}{(2^{j/2}|N_\nu-N_{\nu'}|)^{1+p+q}}\\
\nn&&\times\left(\int_{\S}(2^{j/2}(N-N_\nu))^p\left(\frac{a-a_\nu}{|\nu-\nu'|^\a}\right)^l
H_1(x,\o,\nu,\nu')F_{j,k,\sigma_1,n}( u)\eta^\nu_j(\o)d\o\right)\\
\nn&&\times\left(\int_{\S}(2^{j/2}(N'-N_{\nu'}))^q\left(\frac{a'-a_{\nu'}}{|\nu-\nu'|^\a}\right)^m H_2(x,\o',\nu,\nu')\overline{F_{j,k,\sigma_2,o}( u')}\eta^{\nu'}_j(\o')d\o'\right)d\s\\
\nn&& +\sum_{l,m,n,o,p,q\geq 0}\left(\frac{c^1_{p,q,l,m,n,o}}{(k-k')|\nu-\nu'|^\a}+
\frac{c^2_{p,q,l,m,n,o}}{(k-k')^2|\nu-\nu'|^{2\a}}+\frac{c^3_{p,q,l,m,n,o}}{(k-k')^3|\nu-\nu'|^{3\a}}\right)\\
\nn&&\times\frac{1}{(k-k')^{l+m+n+o}}\int_{\s}\frac{2^{-3j/2}}{(2^{j/2}|N_\nu-N_{\nu'}|)^{1+p+q}}\\
\nn&&\times\left(\int_{\S}(2^{j/2}(N-N_\nu))^p\left(\frac{a-a_\nu}{|\nu-\nu'|^\a}\right)^l
H_3(x,\o,\nu,\nu')F_{j,k,\sigma_1,n}( u)\eta^\nu_j(\o)d\o\right)\\
\nn&&\times\left(\int_{\S}(2^{j/2}(N'-N_{\nu'}))^q\left(\frac{a'-a_{\nu'}}{|\nu-\nu'|^\a}\right)^m H_4(x,\o',\nu,\nu')\overline{F_{j,k,\sigma_2,o}( u')}\eta^{\nu'}_j(\o')d\o'\right)d\s,
\eea
where $(\sigma_1,\sigma_2)=(0,-1)$ in the case of the term involving $D_1$, and $(\sigma_1,\sigma_2)=(-1,0)$ in the case of the term involving $D_2$. We estimate the two terms in the right-hand side of \eqref{td15} starting with the second one. We have:
\bea\label{td16}
\nn&& \bigg|\int_{\s}\frac{2^{-3j/2}}{(2^{j/2}|N_\nu-N_{\nu'}|)^{1+p+q}}\\
\nn&&\times\left(\int_{\S}(2^{j/2}(N-N_\nu))^p\left(\frac{a-a_\nu}{|\nu-\nu'|^\a}\right)^l H_3(x,\o,\nu,\nu')F_{j,k,\sigma_1,n}( u)\eta^\nu_j(\o)d\o\right)\\
\nn&&\times\left(\int_{\S}(2^{j/2}(N'-N_{\nu'}))^q\left(\frac{a'-a_{\nu'}}{|\nu-\nu'|^\a}\right)^mH_4(x,\o',\nu,\nu')\overline{F_{j,k,\sigma_2,o}( u')}\eta^{\nu'}_j(\o')d\o'\right)d\s\bigg|\\
\nn&\leq& \frac{2^{-3j/2}}{(2^{j/2}|\nu-\nu'|)^{1+p+q}}\left(\int_{\S}\norm{H_3}_{\l{\infty}{2}}\norm{F_{j,k,\sigma_1,n}}_{L^2_u}\eta^\nu_j(\o)d\o\right)\\
\nn&&\times\left(\int_{\S}\norm{H_4}_{\lprime{\infty}{2}}\norm{F_{j,k,\sigma_2,o}}_{L^2_{u'}}\eta^{\nu'}_j(\o')d\o'\right)\\
&\lesssim& \frac{2^{-j/2}D^2\gamma_j^{\nu,k}\gamma_j^{\nu',k}}{(2^{j/2}|\nu-\nu'|)^{1+p+q}},
\eea
where we have used {\bf Assumption 2} to estimate $|N_\nu-N_{\nu'}|$, Plancherel to estimate $\norm{F_{j,k,\sigma_1,n}}_{L^2_{ u}}$ and $\norm{F_{j,k,\sigma_2,o}}_{L^2_{ u'}}$, Cauchy-Schwartz in $\o$ and $\o'$, the estimate \eqref{td5} for $H_3$ and $H_4$, and:
\begin{equation}\label{td16b}
\frac{|a-a_{\nu}|}{|\nu-\nu'|^\a}\leq \frac{|\o-\nu|^\a}{|\nu-\nu'|^\a}\norm{\po^\a a}_{\ll{\infty}}\lesssim\frac{\ep}{(2^{j/2}|\nu-\nu'|)^\a}\lesssim 1
\end{equation}
on the support of $\eta_j^\nu$ thanks to {\bf Assumption 2}.

We now estimate the second term in the right-hand side of \eqref{td15}. We have:
\begin{equation}\label{td17}
\begin{array}{l}
\ds \bigg|\int_{\s}\frac{2^{-3j/2}}{(2^{j/2}|N_\nu-N_{\nu'}|)^{1+p+q}}\\
\ds\times\left(\int_{\S}(2^{j/2}(N-N_\nu))^p\left(\frac{a-a_\nu}{|\nu-\nu'|^\a}\right)^l H_1(x,\o,\nu,\nu')F_{j,k,\sigma_1,n}( u)\eta^\nu_j(\o)d\o\right)\\
\ds\times\left(\int_{\S}(2^{j/2}(N'-N_{\nu'}))^q\left(\frac{a'-a_{\nu'}}{|\nu-\nu'|^\a}\right)^mH_2(x,\o',\nu,\nu')\overline{F_{j,k,\sigma_2,o}( u')}\eta^{\nu'}_j(\o')d\o'\right)d\s\bigg|\\
\ds\leq \frac{2^{-3j/2}}{(2^{j/2}|\nu-\nu'|)^{1+p+q}}\left(\int_{\S}\norm{H_1}_{\ll{2}}\norm{F_{j,k,\sigma_1,n}}_{L^{\infty}_u}\eta^\nu_j(\o)d\o\right)\\
\ds\times\normm{\int_{\S}(2^{j/2}(N'-N_{\nu'}))^q\left(\frac{a'-a_{\nu'}}{|\nu-\nu'|^\a}\right)^mH_2(x,\o',\nu,\nu')F_{j,k,\sigma_2,o}( u')\eta^{\nu'}_j(\o')d\o'}_{\ll{2}}\\
\ds\lesssim \frac{2^{-j/2}D\gamma^{\nu,k}_j}{(2^{j/2}|\nu-\nu'|)^{1+p+q}}\\
\ds\times\normm{\int_{\S}(2^{j/2}(N'-N_{\nu'}))^q\left(\frac{a'-a_{\nu'}}{|\nu-\nu'|^\a}\right)^mH_2(x,\o',\nu,\nu')F_{j,k,\sigma_2,o}( u')\eta^{\nu'}_j(\o')d\o'}_{\ll{2}}
\end{array}
\end{equation}
where we have used {\bf Assumption 2} to estimate $|N_\nu-N_{\nu'}|$, Cauchy-Schwartz in $\o$,  the  estimate \eqref{bisoa29} for $\norm{F_{j,k,\sigma_1,n}}_{L^\infty_{ u}}$, and the estimate \eqref{td5} for $H_1$.  Note that the symbol 
$$F=(2^{j/2}(N'-N_{\nu'}))^q\left(\frac{a'-a_{\nu'}}{|\nu-\nu'|^\a}\right)^mH_2(x,\o',\nu,\nu')$$ 
satisfies the following assumptions:
\begin{equation}\label{td18}
\begin{array}{l}
\ds\norm{F}_{\ll{\infty}}\lesssim D,\,\norm{\po F}_{\ll{2}}\lesssim \left(q+\frac{m}{|\nu-\nu'|^\a}\right) D,\\
\ds\textrm{and }\norm{\nabla\po F}_{\ll{2}}\lesssim \left(q^2+\frac{m^2+mq}{|\nu-\nu'|^\a}\right) D,
\end{array}
\end{equation}
where we have used {\bf Assumption 2} for $\po N$, $\po^2N$, $\po a$ and $\po^\a a$, and the assumption \eqref{td6} satisfied by $H_2$. 
We will see in section \ref{bissec:diagonal} that assumptions \eqref{td18} on a symbol is enough to control the diagonal term in $\ll{2}$ (i.e. to obtain the estimate \eqref{bisdiagonal1}). Thus, we obtain:
\bea\label{td19}
\nn&&\normm{\int_{\S}(2^{j/2}(N'-N_{\nu'}))^q\left(\frac{a'-a_{\nu'}}{|\nu-\nu'|^\a}\right)^mH_2(x,\o',\nu,\nu')F_{j,k,\sigma_2,o}( u')\eta^{\nu'}_j(\o')d\o'}_{\ll{2}}\\
&\lesssim &\!\!\!\ds\left(1+q^2+\frac{m^2+mq}{|\nu-\nu'|^\a}\right)D\gamma_j^{\nu',k}.
\eea
\eqref{td17} and \eqref{td19} imply:
\bea\label{bisoa30b}
\nn&&\bigg|\int_{\s}\frac{2^{-3j/2}}{(2^{j/2}|N_\nu-N_{\nu'}|)^{1+p+q}}\\
\nn&&\times\left(\int_{\S}(2^{j/2}(N-N_\nu))^p\left(\frac{a-a_\nu}{|\nu-\nu'|^\a}\right)^l H_1(x,\o,\nu,\nu')F_{j,k,\sigma_1,n}( u)\eta^\nu_j(\o)d\o\right)\\
\nn&&\times\left(\int_{\S}(2^{j/2}(N'-N_{\nu'}))^q\left(\frac{a'-a_{\nu'}}{|\nu-\nu'|^\a}\right)^mH_2(x,\o',\nu,\nu')\overline{F_{j,k,\sigma_2,o}( u')}\eta^{\nu'}_j(\o')d\o'\right)d\s\bigg|\\
&\lesssim& \frac{(1+q^2+m^2)2^{-j/2}D^2|\nu-\nu'|^{-\a}}{(2^{j/2}|\nu-\nu'|)^{1+p+q}}\ga_j^{\nu,k}\ga_j^{\nu',k}.
\eea
 
Finally, \eqref{td15}, \eqref{td16} and \eqref{bisoa30b} yield:
\begin{equation}\label{td20}
\begin{array}{ll}
&\ds\left|\int_{\s}U^{\nu,k}_jf(x)\overline{U^{\nu',k}_jf(x)}d\s\right| \\
\ds\lesssim &\ds\sum_{l,m,n,o,p,q\geq 0}\left(\frac{c^1_{p,q,l,m,n,o}}{|k-k'||\nu-\nu'|^\a}+
\frac{c^2_{p,q,l,m,n,o}}{|k-k'|^2|\nu-\nu'|^{2\a}}+\frac{c^3_{p,q,l,m,n,o}}{|k-k'|^3|\nu-\nu'|^{3\a}}\right)\\
&\ds\times \frac{(1+q^2+m^2)2^{-j/2}D^2|\nu-\nu'|^{-\a}}{|k-k'|^{l+m+n+o}(2^{j/2}|\nu-\nu'|)^{1+p+q}}\ga_j^{\nu,k}\ga_j^{\nu',k}\\
& +\ds\sum_{l,m,n,o,p,q\geq 0}\left(\frac{c^1_{p,q,l,m,n,o}}{|k-k'||\nu-\nu'|^\a}+
\frac{c^2_{p,q,l,m,n,o}}{|k-k'|^2|\nu-\nu'|^{2\a}}+\frac{c^3_{p,q,l,m,n,o}}{|k-k'|^3|\nu-\nu'|^{3\a}}\right)\\
&\ds\times \frac{2^{-j/2}D^2}{|k-k'|^{l+m+n+o}(2^{j/2}|\nu-\nu'|)^{1+p+q}}\ga_j^{\nu,k}\ga_j^{\nu',k}\\
\ds\lesssim &\ds\frac{D^2\ga_j^{\nu,k}\ga_j^{\nu',k}}{|k-k'|2^{j/2(1-4\a)}(2^{j/2}|\nu-\nu'|)^{1+4\a}},
\end{array}
\end{equation}
which concludes the proof of estimate \eqref{bisoa12}.

\subsubsection{End of the proof of Proposition \ref{bisorthoangle}}

We have proved the estimates \eqref{bisoa11} and \eqref{bisoa12} in the two previous sections. Since \eqref{bisoa11} and \eqref{bisoa12} yield \eqref{bisoa2} (see section \ref{bissec:keyest}), this concludes the proof of Proposition \ref{bisorthoangle}. \QED

\subsection{Proof of Proposition \ref{bisdiagonal} (control of the diagonal term)}\label{bissec:diagonal}

We have to prove \eqref{bisdiagonal1}:
\begin{equation}\label{bisdi1}
\norm{U^\nu_jf}_{\ll{2}}\lesssim D\ga^\nu_j.
\end{equation}
Recall that $U^\nu_j$ is given by:
\begin{equation}\label{bisdi2}
U^\nu_jf(x)=\int_{\S}bF_j( u)\eta_j^\nu(\o)d\o,
\end{equation}
where $F_j( u)$ is defined by:
\begin{equation}\label{bisdi3}
F_j( u)=\int_0^{+\infty}e^{i\la u}\psi(2^{-j}\la)f(\la\o)\la^2d\la.
\end{equation}
We decompose $U^\nu_j$ in the sum of two terms:
\begin{equation}\label{bisdi4}
U^\nu_jf(x)=b(x,\nu)\int_{\S}F_j( u)\eta_j^\nu(\o)d\o+\int_{\S}(b(x,\o)-b(x,\nu))F_j( u)\eta_j^\nu(\o)d\o.
\end{equation}

We start with the first term. The assumption \eqref{thregx1s} on $b$ implies:
\begin{equation}\label{bisdi6}
\normm{b(x,\nu)\int_{\S}F_j( u)\eta_j^\nu(\o)d\o}_{\ll{2}}\lesssim D\normm{\int_{\S}F_j( u)\eta_j^\nu(\o)d\o}_{\ll{2}}.
\end{equation}
The following proposition allows us to estimate the right-hand side of \eqref{bisdi6}.
\begin{proposition}\label{bisdiprop1}
The right-hand side of \eqref{bisdi6} satisfies the following bound:
\begin{equation}\label{bisdiprop2}
\normm{\int_{\S}F_j( u)\eta_j^\nu(\o)d\o}_{\ll{2}}\lesssim\gamma_j^\nu.
\end{equation}
\end{proposition}
The proof of Proposition \ref{bisdiprop1} is postponed to section \ref{sec:bisdiprop1}. In the rest of this section, we show how Proposition \ref{bisdiprop1} yields the proof of \eqref{bisdi1}. In particular, \eqref{bisdiprop2} together with \eqref{bisdi6} implies the following bound for the first term in the right-hand side of \eqref{bisdi4}:
\begin{equation}\label{bisdi6b}
\normm{b(x,\nu)\int_{\S}F_j( u)\eta_j^\nu(\o)d\o}_{\ll{2}}\lesssim D\gamma_j^\nu.
\end{equation}

We turn to the second term in the right-hand side of \eqref{bisdi4}. We have:
\begin{equation}\label{bisdi7}
\begin{array}{ll}
&\ds\normm{\int_{\S}(b(x,\o)-b(x,\nu))F_j( u)\eta_j^\nu(\o)d\o}_{\ll{2}}\\
\ds\leq &\ds \int_{\S}\norm{b(x,\o)-b(x,\nu)}_{\l{\infty}{2}}\norm{F_j}_{L^2_{ u}}\eta_j^\nu(\o)d\o\\
\ds\leq &\ds\int_{\S}|\o-\nu|(\norm{\po b}_{\ll{2}}+\norm{\nabla\po b}_{\ll{2}})\norm{F_j}_{L^2_{ u}}\eta_j^\nu(\o)d\o\\
\lesssim &\ds D\gamma_j^\nu,
\end{array}
\end{equation}
where we have used Plancherel to estimate $\norm{F_j}_{L^2_{ u}}$, Cauchy-Schwartz in $\o$, the fact that $H^1(\s)$ embeds in $\l{\infty}{2}$ (see \cite{param1} Corollary 3.6 for a proof only using the regularity given by {\bf Assumption 1}), the assumption \eqref{threomega1s} on $b$, and the fact that $|\o-\nu|\lesssim 2^{-j/2}$ on the support of $\eta_j^\nu$.

Finally, \eqref{bisdi4}, \eqref{bisdi6b} and \eqref{bisdi7} yield the wanted estimate \eqref{bisdi1} which concludes the proof of Proposition \ref{bisdiagonal}. \QED.

\subsubsection{Proof of Proposition \ref{bisdiprop1}}\label{sec:bisdiprop1}

Recall that $\int_{\S}F_j( u)\eta_j^\nu(\o)d\o$ is given by:
\begin{equation}\label{di10}
\int_{\S}F_j( u)\eta_j^\nu(\o)d\o=\int_{\S}\int_0^{+\infty}e^{i\la u}\psi(2^{-j}\la)\eta_j^\nu(\o)f(\la\o)\la^2d\la d\o.
\end{equation}
Relying on the classical $TT^*$ argument, \eqref{bisdiprop2} is equivalent to proving the boundedness in $\ll{2}$ of the operator whose kernel $K$ is given by:
\begin{equation}\label{di11}
K(x,y)=\int_{\S}\int_0^{+\infty}e^{i\la u(x,\o)-i\la u(y,\o)}\psi(2^{-j}\la)\eta_j^\nu(\o)\la^2d\la d\o,\,x, y\in\s.
\end{equation}
The decay satisfied by this kernel is stated in the following proposition.
\begin{proposition}\label{di12}
The kernel $K$ defined in \eqref{di11} satisfies the following decay estimate for all $x, y$ in $\s$:
\begin{equation}\label{di13}
\begin{array}{ll}
\ds |K(x,y)|\lesssim & \ds\frac{2^j}{(1+|2^j|u(x,\nu)-u(y,\nu)|-2^{j/2}|\po u(x,\nu)-\po u(y,\nu)||)^2}\\
& \ds\times\frac{2^j}{(1+2^{j/2}|\po u(x,\nu)-\po u(y,\nu)|)^3}.
\end{array}
\end{equation}
\end{proposition}
The proof of Proposition \ref{di12} is postponed to section \ref{sec:di20}. In the rest of this section, we show how \eqref{di13} implies Proposition \ref{bisdiprop1}. According to Schur's Lemma, the operator whose kernel is $K$ is bounded on $\ll{2}$ provided we can prove the following bound:
\begin{equation}\label{di14}
\sup_{x\in\s}\int_\s |K(x,y)|dy<+\infty,\,\sup_{y\in\s}\int_\s |K(x,y)|dx<+\infty.
\end{equation}
Due to the symmetry of $K$ in $x, y$, the two bounds in \eqref{di14} are obtained in the same way. We 
focus on establishing the first bound. In view of \eqref{di13}, we have:
\begin{equation}\label{di15}
\begin{array}{ll}
\ds \int_\s |K(x,y)|dy\lesssim & \ds\int_\s\frac{2^j}{(1+|2^j|u(x,\nu)-u(y,\nu)|-2^{j/2}|\po u(x,\nu)-\po u(y,\nu)||)^2}\\
& \ds\times\frac{2^j}{(1+2^{j/2}|\po u(x,\nu)-\po u(y,\nu)|)^3}dy.
\end{array}
\end{equation}
Now, according to {\bf Assumption 4}, there is a global change of variable on $\s$ $\phi_\nu:\s\rightarrow\R^3$ defined by:
\begin{equation}\label{di16}
\phi_\nu(x):=u(x,\nu)\nu+\po u(x,\nu),
\end{equation}
such that $\phi_\nu$ is a bijection, and the determinant of its Jacobian satisfies the following estimate:
\begin{equation}\label{di17}
\norm{|\det(\textrm{Jac}\phi_\nu)|-1}_{\ll{\infty}}\lesssim\ep.
\end{equation}
Using the change of variable $y\rightarrow\underline{y}=\phi_\nu(y)-\phi_\nu(x)\in\R^3$ in the right-hand side of \eqref{di15} together with \eqref{di17}, we obtain:
\begin{equation}\label{di18}
\ds \int_\s |K(x,y)|dy\lesssim\ds\int_{\R^3}\frac{2^j}{(1+|2^j|\underline{y}\c\nu|-2^{j/2}|\underline{y}'||)^2}\frac{2^j}{(1+2^{j/2}|\underline{y}'|)^3}d\underline{y},
\end{equation}
where $y=y\c\nu+y'$ and $y'\c\nu=0$. Making the change of variable $y\rightarrow z$ where $z$ is defined by $z\c\nu=2^j\underline{y}\c\nu$ and $z'=2^{j/2}\underline{y}'$ in the right-hand side of \eqref{di18}, and remarking that $z\cdot\nu$ is one dimensional, and $z'$ is two dimensional, we obtain:
\begin{equation}\label{di19}
\ds \int_\s |K(x,y)|dy\lesssim\ds\int_{\R^3}\frac{dz}{(1+||z\c\nu|-|z'||)^2(1+|z'|)^3}\lesssim 1.
\end{equation}
\eqref{di19} implies the first bound in \eqref{di14}. $K$ being symmetric with respect to $x, y$, the second bound in \eqref{di14} is also true. Thus, the operator whose kernel is $K$ is bounded on $\ll{2}$ which concludes the proof of Proposition \ref{bisdiprop1}. \QED

\subsubsection{Proof of Proposition \ref{di12}}\label{sec:di20}

Recall the definition of $K$:
\begin{equation}\label{di21}
K(x,y)=\int_{\S}\int_0^{+\infty}e^{i\la u(x,\o)-i\la u(y,\o)}\psi(2^{-j}\la)\eta_j^\nu(\o)\la^2d\la d\o,\,x, y\in\s.
\end{equation}
We need to prove that $K$ satisfies the following decay estimate for all $x, y$ in $\s$:
\begin{equation}\label{di22}
\begin{array}{ll}
\ds |K(x,y)|\lesssim & \ds\frac{2^j}{(1+|2^j|u(x,\nu)-u(y,\nu)|-2^{j/2}|\po u(x,\nu)-\po u(y,\nu)||)^2}\\
& \ds\times\frac{2^j}{(1+2^{j/2}|\po u(x,\nu)-\po u(y,\nu)|)^3}.
\end{array}
\end{equation}

\vspace{0.2cm}

\noindent{\bf Proof of \eqref{di22}.} Recall from Remark \ref{reduccompact} that $u(x,\o)$ is exactly equal to $\xo$ in $|x|\geq 2$. Thus, we may restrict ourselves to $|x|\leq 2$ where we have in view of {\bf Assumption 2}:
\begin{equation}\label{reguomega1}
|u(x,\o)|+|\po u(x,\o)|+|\po^2u(x,\o)|+|\po^3u(x,\o)|\lesssim 1,\,\forall x\textrm{ with }|x|\leq 2,\,\forall \o\in\S.
\end{equation}
We will obtain \eqref{di22} as a consequence of the following estimate:
\begin{equation}\label{di23}
\ds |K(x,y)|\lesssim \ds\int_{\S}\int_0^{+\infty}\frac{1}{(1+2^j|u(x,\o)-u(y,\o)|)^3}\widetilde{\psi}(2^{-j}\la)\widetilde{\eta}_j^\nu(\o)\la^2d\la d\o,
\end{equation}
where $\widetilde{\psi}$ is smooth and compactly supported in $(0,+\infty)$ and $\widetilde{\eta}^\nu_j$ is bounded on $\S$ and has the same support as $\eta^\nu_j$.
Indeed, we have:
\begin{equation}\label{di24}
\begin{array}{rl}
\ds u(x,\o)-u(y,\o)= &\ds u(x,\nu)-u(y,\nu)+(\po u(x,\nu)-\po u(y,\nu))(\o-\nu)\\
&\ds +O(|\nu-\o|^2),\\
\ds\po u(x,\o)-\po u(y,\o)= &\ds\po u(x,\nu)-\po u(y,\nu)+O(|\nu-\o|),
\end{array}
\end{equation}
where we have used a Taylor expansion in $\o$ together with \eqref{reguomega1}. Using the fact that $2^{j/2}|\o-\nu|\lesssim 1$ on the support of $\widetilde{\eta}_j^\nu$ together with \eqref{di24}, we obtain for $\o$ in the support of $\widetilde{\eta}_j^\nu$:
\bea\label{di25}
&&\nn 1+|2^j|u(x,\nu)-u(y,\nu)|-2^{j/2}|\po u(x,\nu)-\po u(y,\nu)||\lesssim 1+2^j|u(x,\o)-u(y,\o)|,\\
&&\textrm{and } 1+2^{j/2}|\po u(x,\nu)-\po u(y,\nu)|\lesssim 1+2^{j/2}|\po u(x,\o)-\po u(y,\o)|.
\eea
\eqref{di23} and \eqref{di25} imply:
\begin{equation}\label{di26}
\begin{array}{r}
\ds |K(x,y)|\lesssim\ds\frac{1}{(1+|2^j|u(x,\nu)-u(y,\nu)|-2^{j/2}|\po u(x,\nu)-\po u(y,\nu)||)^2}\\
\ds\times\frac{1}{(1+2^{j/2}|\po u(x,\nu)-\po u(y,\nu)|)^3}\int_{\S}\int_0^{+\infty}\widetilde{\psi}(2^{-j}\la)\widetilde{\eta}_j^\nu(\o)\la^2d\la d\o.
\end{array}
\end{equation}
Now, we have:
\begin{equation}\label{di27}
\int_{\S}\int_0^{+\infty}\widetilde{\psi}(2^{-j}\la)\widetilde{\eta}_j^\nu(\o)\la^2d\la d\o=\left(\int_0^{+\infty}\widetilde{\psi}(2^{-j}\la)\la^2d\la\right)\left(\int_{\S}\widetilde{\eta}_j^\nu(\o) d\o\right)\lesssim 2^{2j},
\end{equation}
where we have used the fact that $\widetilde{\eta}^\nu_j$ is bounded on $\S$ and the fact that the support of $\widetilde{\eta}_j^\nu$ is two dimensional with a diameter of size $\sim 2^{-j/2}$. Finally, \eqref{di26} and \eqref{di27} imply \eqref{di22} which is the wanted estimate.

\vspace{0.2cm}

\noindent{\bf Proof of \eqref{di23}.} To conclude the proof of Proposition \ref{di12}, it remains to prove \eqref{di23}. This will follow by performing three integrations by parts with respect to $\o$ and two integrations by parts with respect to $\la$. We start with the integrations by parts with respect to $\o$. 
Our goal is to show that $K(x,y)$ is a sum of terms of the form:
\begin{equation}\label{di47}
\begin{array}{r}
\ds\int_{\S}\int_0^{+\infty}\frac{e^{i\la u(x,\o)-i\la u(y,\o)}F(x,y,\o,\nu)(\la(\po u(x,\o)-\po u(y,\o))^2)^l}{(1+\la|\po u(x,\o)-\po u(y,\o)|^2)^m}\\
\times\ds\psi(2^{-j}\la)\widetilde{\eta}_j^\nu(\o)\la^2d\la d\o,
\end{array}
\end{equation}
where $l, m$ are integers, where $F$ does not depend on $\la$, where $\widetilde{\eta}^\nu_j$ is bounded on $\S$ and has the same support as $\eta^\nu_j$ and where the integrand in \eqref{di47} satisfies:
\begin{equation}\label{di48}
\left|\frac{F(x,y,\o,\nu)(\la(\po u(x,\o)-\po u(y,\o))^2)^l}{(1+\la|\po u(x,\o)-\po u(y,\o)|^2)^m}\right|\lesssim\frac{1}{(1+2^{j/2}|\po u(x,\o)-\po u(y,\o)|)^3}.
\end{equation}
To this end, we use:
\begin{equation}\label{di33}
e^{i\la(u(x,\o)-u(y,\o))}=\frac{(1-i(\po u(x,\o)-\po u(y,\o))\po)e^{i\la(u(x,\o)-u(y,\o))}}{(1+\la|\po u(x,\o)-\po u(y,\o)|^2)}.
\end{equation}
We integrate by parts once in the integral \eqref{di21} defining $K$ using \eqref{di33}. We obtain:
\begin{equation}\label{di34}
\begin{array}{l}
\ds K(x,y)= \ds\int_{\S}\int_0^{+\infty}\frac{e^{i\la u(x,\o)-i\la u(y,\o)}}{1+\la|\po u(x,\o)-\po u(y,\o)|^2}\psi(2^{-j}\la)\eta_j^\nu(\o)\la^2d\la d\o\\
\ds +\int_{\S}\int_0^{+\infty}\frac{e^{i\la u(x,\o)-i\la u(y,\o)}A_1}{1+\la|\po u(x,\o)-\po u(y,\o)|^2}\psi(2^{-j}\la)\eta_j^\nu(\o)\la^2d\la d\o\\
\ds +\int_{\S}\int_0^{+\infty}\frac{e^{i\la u(x,\o)-i\la u(y,\o)}A_2}{(1+\la|\po u(x,\o)-\po u(y,\o)|^2)^2}\psi(2^{-j}\la)\eta_j^\nu(\o)\la^2d\la d\o\\
\ds +\int_{\S}\int_0^{+\infty}\frac{e^{i\la u(x,\o)-i\la u(y,\o)}A_3}{1+\la|\po u(x,\o)-\po u(y,\o)|^2}\psi(2^{-j}\la)2^{-j/2}\po\eta_j^\nu(\o)\la^2d\la d\o,
\end{array}
\end{equation}
where $A_1, A_2, A_3$ are given by:
\begin{equation}\label{di35}
\begin{array}{l}
\ds A_1=i(\po^2 u(x,\o)-\po^2 u(y,\o)),\\ 
\ds A_2=-2i\la(\po^2 u(x,\o)-\po^2 u(y,\o))(\po u(x,\o)-u(y,\o))^2,\\
\ds A_3=i2^{j/2}(\po u(x,\o)-\po u(y,\o)).
\end{array}
\end{equation}

\vspace{0.2cm}

\noindent{\bf The first term in the right-hand side of \eqref{di34}.} Integrating by parts in the first term of the right-hand side of \eqref{di34} using \eqref{di33} yields:
\begin{equation}\label{di36}
\begin{array}{ll}
&\ds\int_{\S}\int_0^{+\infty}\frac{e^{i\la u(x,\o)-i\la u(y,\o)}}{1+\la|\po u(x,\o)-\po u(y,\o)|^2}\psi(2^{-j}\la)\eta_j^\nu(\o)\la^2d\la d\o\\
\ds = &\ds\int_{\S}\int_0^{+\infty}\frac{e^{i\la u(x,\o)-i\la u(y,\o)}A_0^1}{(1+\la|\po u(x,\o)-\po u(y,\o)|^2)^2}\psi(2^{-j}\la)\eta_j^\nu(\o)\la^2d\la d\o\\
&\ds +\int_{\S}\int_0^{+\infty}\frac{e^{i\la u(x,\o)-i\la u(y,\o)}A_0^2}{(1+\la|\po u(x,\o)-\po u(y,\o)|^2)^3}\psi(2^{-j}\la)\eta_j^\nu(\o)\la^2d\la d\o\\
&\ds +\int_{\S}\int_0^{+\infty}\frac{e^{i\la u(x,\o)-i\la u(y,\o)}A_0^3}{(1+\la|\po u(x,\o)-\po u(y,\o)|^2)^2}\psi(2^{-j}\la)2^{-j/2}\po\eta_j^\nu(\o)\la^2d\la d\o,
\end{array}
\end{equation}
where $A_0^1, A_0^2, A_0^3$ are given by:
\begin{equation}\label{di37}
\begin{array}{l}
\ds A_0^1=1+i(\po^2 u(x,\o)-\po^2 u(y,\o)),\\ 
\ds A_0^2=-4i\la(\po^2 u(x,\o)-\po^2 u(y,\o))(\po u(x,\o)-u(y,\o))^2,\\
\ds A_0^3=i2^{j/2}(\po u(x,\o)-\po u(y,\o)).
\end{array}
\end{equation}
\eqref{reguomega1} implies the following bound for the integrand in the right-hand side of \eqref{di36}:
\begin{equation}\label{di39}
\begin{array}{l}
\ds\frac{|A_0^1|}{(1+\la|\po u(x,\o)-\po u(y,\o)|^2)^2}+\frac{|A_0^2|}{(1+\la|\po u(x,\o)-\po u(y,\o)|^2)^3}\\
\ds +\frac{|A_0^3|}{(1+\la|\po u(x,\o)-\po u(y,\o)|^2)^2}\lesssim \frac{1}{(1+2^{j/2}|\po u(x,\o)-\po u(y,\o)|)^{3}}.
\end{array}
\end{equation}
Also,  we have:
\begin{equation}\label{di38}
|\po\eta_j^\nu(\o)|\lesssim 2^{j/2}\textrm{ for all }\o\in\S
\end{equation}
so that $2^{-j/2}\po\eta^\nu_j$ satisfies the assumptions of $\widetilde{\eta}^\nu_j$. 
Thus, the first term of the right-hand side of \eqref{di34} satisfies \eqref{di47} \eqref{di48}. 

\vspace{0.2cm}

\noindent{\bf The terms involving $A_1$, $A_2$ and $A_3$ in the right-hand side of \eqref{di34}.} 

Integrating by parts in the term involving $A_1$ of the right-hand side of \eqref{di34} using \eqref{di33} yields:
\bea\label{di36:back}
&&\ds\int_{\S}\int_0^{+\infty}\frac{e^{i\la u(x,\o)-i\la u(y,\o)}A_1}{1+\la|\po u(x,\o)-\po u(y,\o)|^2}\psi(2^{-j}\la)\eta_j^\nu(\o)\la^2d\la d\o\\
\nn &=& \ds\int_{\S}\int_0^{+\infty}\frac{e^{i\la u(x,\o)-i\la u(y,\o)}A_0^1A_1}{(1+\la|\po u(x,\o)-\po u(y,\o)|^2)^2}\psi(2^{-j}\la)\eta_j^\nu(\o)\la^2d\la d\o\\
\nn&&\ds +\int_{\S}\int_0^{+\infty}\frac{e^{i\la u(x,\o)-i\la u(y,\o)}A_0^2A_1}{(1+\la|\po u(x,\o)-\po u(y,\o)|^2)^3}\psi(2^{-j}\la)\eta_j^\nu(\o)\la^2d\la d\o\\
\nn&&\ds +\int_{\S}\int_0^{+\infty}\frac{e^{i\la u(x,\o)-i\la u(y,\o)}A_0^3A_1}{(1+\la|\po u(x,\o)-\po u(y,\o)|^2)^2}\psi(2^{-j}\la)2^{-j/2}\po\eta_j^\nu(\o)\la^2d\la d\o\\
\nn&&\ds+\int_{\S}\int_0^{+\infty}\frac{e^{i\la u(x,\o)-i\la u(y,\o)}i(\po u(x,\o)-\po u(y,\o))\po A_1}{(1+\la|\po u(x,\o)-\po u(y,\o)|^2)^2}\psi(2^{-j}\la)\eta_j^\nu(\o)\la^2d\la d\o,
\eea
where $A_0^1, A_0^2, A_0^3$ are given by \eqref{di37}. In view of \eqref{di39}, we have the following bound for the integrand in the right-hand side of \eqref{di36:back}:
\begin{equation}\label{di39:back}
\begin{array}{l}
\ds\frac{|A_0^1A_1|}{(1+\la|\po u(x,\o)-\po u(y,\o)|^2)^2}+\frac{|A_0^2A_1|}{(1+\la|\po u(x,\o)-\po u(y,\o)|^2)^3}\\
\ds +\frac{|A_0^3A_1|}{(1+\la|\po u(x,\o)-\po u(y,\o)|^2)^2}+\frac{|\po u(x,\o)-\po u(y,\o)||\po A_1|}{(1+\la|\po u(x,\o)-\po u(y,\o)|^2)^2}\\
\ds\lesssim \frac{|A_1|}{(1+2^{j/2}|\po u(x,\o)-\po u(y,\o)|)^{3}}+\frac{|\po u(x,\o)-\po u(y,\o)||\po A_1|}{(1+2^{j/2}|\po u(x,\o)-\po u(y,\o)|)^{3}}.
\end{array}
\end{equation}
Now, in view of the definition \eqref{di35} of $A_1$ and the estimate \eqref{reguomega1}, we have:
\begin{equation}\label{back2back}
|A_1|\lesssim 1,\textrm{ and }|\po u(x,\o)-\po u(y,\o)||\po A_1|\lesssim 1.
\end{equation}
In view of \eqref{di39:back} and \eqref{back2back} the term involving $A_1$ of the right-hand side of \eqref{di34} satisfies \eqref{di47} \eqref{di48}. 

We proceed similarly for $A_2$ and $A_3$. In particular, in view of the definition \eqref{di35} of $A_2, A_3$ and the estimate \eqref{reguomega1}, we may replace \eqref{back2back} with the following estimates:
$$|A_2|\lesssim (2^{j/2}|\po u(x,\o)-\po u(y,\o)|)^2,$$
$$|\po u(x,\o)-\po u(y,\o)||\po A_2|\lesssim (2^{j/2}|\po u(x,\o)-\po u(y,\o)|)^2,$$
$$|A_3|\lesssim 2^{j/2}|\po u(x,\o)-\po u(y,\o)|,$$
and
$$|\po u(x,\o)-\po u(y,\o)||\po A_3|\lesssim 2^{j/2}|\po u(x,\o)-\po u(y,\o)|.$$
Finally, the four terms in the right-hand side of \eqref{di34} satisfy the estimates \eqref{di47} \eqref{di48}. Thus $K(x,y)$ satisfies \eqref{di47} \eqref{di48}. 

\vspace{0.2cm}

\noindent{\bf Integration by parts with respect to $\la$ and end of the proof of \eqref{di23}.} In order to obtain \eqref{di23}, we still need to perform two integration by parts with respect to $\la$ in \eqref{di47}. We have:
\begin{equation}\label{di28}
e^{i\la(u(x,\o)-u(y,\o))}=\frac{(1-i2^{j}(u(x,\o)-u(y,\o))2^{j}\partial_\la)e^{i\la(u(x,\o)-u(y,\o))}}{(1+2^{2j}|u(x,\o)-u(y,\o)|^2)}.
\end{equation}
Notice that the only term depending on $\la$ under the integral \eqref{di47} is:
\begin{equation}\label{di28b}
\frac{\psi(2^{-j}\la)\la^{2+l}}{(1+\la|\po u(x,\o)-\po u(y,\o)|^2)^m}.
\end{equation}
Now, we have:
\begin{equation}\label{di29}
\begin{array}{ll}
&\ds 2^j\partial_\la\left(\frac{\psi(2^{-j}\la)\la^{2+l}}{(1+\la|\po u(x,\o)-\po u(y,\o)|^2)^m}\right)\\
= &\ds\frac{\left((2+l-m)\overline{\psi}(2^{-j}\la)+\psi'(2^{-j}\la)\right)\la^{2+l}}{(1+\la|\po u(x,\o)-\po u(y,\o)|^2)^m}\\
&\ds +\frac{m\overline{\psi}(2^{-j}\la)\la^{2+l}}{(1+\la|\po u(x,\o)-\po u(y,\o)|^2)^{m+1}}\textrm{ where }\overline{\psi}(\la)=\frac{\psi(\la)}{\la}.
\end{array}
\end{equation}
Thus integrating by parts in $\la$  in the integral \eqref{di47} using \eqref{di28} essentially 
divides the integrand by $1+2^{j}|u(x,\o)-u(y,\o)|$. In particular, after two integrations by parts 
using \eqref{di28} in the integral \eqref{di47}, and together with the estimate \eqref{di48}, we obtain 
that $K(x,y)$ is a sum of terms of the form:
\begin{equation}\label{di47b}
\begin{array}{r}
\ds\int_{\S}\int_0^{+\infty}\frac{e^{i\la u(x,\o)-i\la u(y,\o)}F(x,y,\o,\nu)(\la(\po u(x,\o)-\po u(y,\o))^2)^l}{(1+\la|\po u(x,\o)-\po u(y,\o)|^2)^m(1+2^{j}|u(x,\o)-u(y,\o)|)^2}\\
\ds\times\widetilde{\psi}(2^{-j}\la)\widetilde{\eta}_j^\nu(\o)\la^2d\la d\o,
\end{array}
\end{equation}
where $l, m$ are integers, where $\widetilde{\psi}$ is smooth and compactly supported in $(0,+\infty)$,  where $\widetilde{\eta}^\nu_j$ is bounded on $\S$ and has the same support as $\eta^\nu_j$ and where the integrand in \eqref{di47b} satisfies:
\begin{equation}\label{di48b}
\begin{array}{ll}
&\ds\left|\frac{F(x,y,\o,\nu)(\la(\po u(x,\o)-\po u(y,\o))^2)^l}{(1+\la|\po u(x,\o)-\po u(y,\o)|^2)^m(1+2^{j}|u(x,\o)-u(y,\o)|)^2}\right|\\
\ds\lesssim &\ds \frac{1}{(1+2^{j}|u(x,\o)-u(y,\o)|)^2(1+2^{j/2}|\po u(x,\o)-\po u(y,\o)|)^3}.
\end{array}
\end{equation}
 Finally, \eqref{di47b} and \eqref{di48b} yield \eqref{di23} which is the wanted estimate. This concludes the proof of Proposition \ref{di12}. \QED

\section{Proof of Theorem \ref{th2}}\label{sec:th2}

In order to prove Theorem \ref{th2}, we first show that the Fourier integral operator $U$ of Theorem 
\ref{th1} almost preserve the $L^2$ norm provided we make additional assumptions on its symbol. 
We then use this observation to prove the estimate \eqref{choicef4bis}. Finally, we conclude the proof of Theorem \ref{th2} by establishing the existence of $(f_+,f_-)$ solution of the system \eqref{choicef3bis}. 

\subsection{A refinement of Theorem \ref{th1}}

In Theorem \ref{th1}, we have proved that the Fourier integral operator $U$ with phase $u$ and symbol 
$b$ is bounded on $\ll{2}$ provided $u$ satisfies {\bf Assumption 1}, {\bf Assumption 2} and {\bf Assumption 4}, and the symbol $b$ satisfies \eqref{thregx1s} \eqref{threomega1s}. We now would like to prove that $U$ satisfies the following bound from below:
\begin{equation}\label{re1}
\norm{f}_{\le{2}}\lesssim\norm{Uf}_{\ll{2}},
\end{equation}
provided $u$ also satisfies {\bf Assumption 5} and under additional assumptions on the symbol $b$. This is the aim of the following proposition.

\begin{proposition}\label{re2}
Let $u$ be a function on $\s\times\S$ satisfying {\bf Assumption 1}, {\bf Assumption 2}, {\bf Assumption 4} and {\bf Assumption 5}. Let $U$ the Fourier integral operator with phase $u(x,\o)$ and symbol $b(x,\o)$:
\begin{equation}\label{fio1}
Uf(x)=\int_{\S}\int_{0}^{+\infty}e^{i\lambda u(x,\o)}b(x,\o)f(\lambda\o)\lambda^2 d\lambda d\o.
\end{equation}
We assume furthermore that $b(x,\o)$ satisfies:
\begin{equation}\label{re3}
\norm{\po b}_{\ll{2}}+\norm{\nabla\po b}_{\ll{2}}\lesssim 1,
\end{equation}
\begin{equation}\label{re4}
\norm{b-1}_{\ll{\infty}}+\norm{\nabla b}_{\l{\infty}{2}}+\norm{\nabb\nabla b}_{\ll{2}}\lesssim \ep,
\end{equation}
and
\begin{equation}\label{re5}
\begin{array}{l}
\ds\nabn b=b^j_1+b^j_2\textrm{ where }\norm{b^j_1}_{\ll{2}}\lesssim 2^{-\frac{j}{2}}\ep,\,\norm{b^j_2}_{\l{\infty}{2}}\lesssim\ep\\
\ds\textrm{and }\norm{\nabn b^j_2}_{\ll{2}}+\norm{b^j_2}_{\l{2}{\infty}}\lesssim 2^{\frac{j}{2}}\ep.
\end{array}
\end{equation}
Then, $U$ is bounded on $L^2$ and satisfies the estimate:
\begin{equation}\label{re6}
\norm{f}_{\le{2}}\lesssim\norm{Uf}_{\ll{2}}.
\end{equation}
\end{proposition}

\begin{remark}
Notice that the only difference in the assumptions with respect to Theorem \ref{th1} lies in the fact that $u$ also satisfies {\bf Assumption 5} and in the constant $D$ which has been replaced by 1 in \eqref{re3} and by $\ep$ in \eqref{re4} \eqref{re5}.
\end{remark}

We now turn to the proof of Proposition \ref{re2}. We review the three steps of Theorem \ref{th1}   - decomposition in frequency, decomposition in angle, and control of the diagonal term -
 indicating each time how to refine the estimates.

\subsubsection{Step 1: decomposition in frequency}

As in step 1 of the proof of Theorem \ref{th1}, we decompose $Uf$ in frequency:
\begin{equation}\label{re7}
Uf(x)=\sum_{j\geq -1}U_jf(x),
\end{equation}
where the operators $U_j$ are defined by \eqref{bisb9} \eqref{bisb10}. We have:
\begin{equation}\label{re8}
\norm{Uf}^2_{\ll{2}}=\sum_{|j-l|\leq 2}\int_{\s}U_jf(x)\overline{U_lf(x)}d\s+\sum_{|j-l|> 2}\int_{\s}U_jf(x)\overline{U_lf(x)}d\s.
\end{equation}
Now, the proof of Proposition \ref{bisorthofreq} together with the fact that $b$ satisfies \eqref{re4} \eqref{re5} immediately yields:
\begin{equation}\label{re9}
\left|\sum_{|j-l|> 2}\int_{\s}U_jf(x)\overline{U_lf(x)}d\s\right|\lesssim \ep\norm{f}^2_{\le{2}}.
\end{equation}
Thus, together with \eqref{re8}, we obtain:
\begin{equation}\label{re10}
\norm{Uf}^2_{\ll{2}}=\sum_{|j-l|\leq 2}\int_{\s}U_jf(x)\overline{U_lf(x)}d\s+O(\ep)\norm{f}_{\le{2}}^2.
\end{equation}

\begin{remark}
The sum over $|j-l|\leq 2$ in the right-hand side of \eqref{re10} corresponds to the terms such that the support of $\psi(2^{-j}\la)$ and the support of $\psi(2^{-j}\la')$ have a non empty intersection, where $\psi$ has been introduced in \eqref{bisb7}. 
\end{remark}

\subsubsection{Step 2: decomposition in angle}

As in step 2 of the proof of Theorem \ref{th1}, we decompose $U_jf$ in angle:
\begin{equation}\label{re11}
U_jf(x)=\sum_{\nu\in\Gamma}U^\nu_jf(x),
\end{equation}
where the operators $U_j^\nu$ are defined by \eqref{bisb16}. 
In order to control the diagonal term in a third step (see next section), we have to modify slightly the size of the support of our partition of unity $\eta_j^\nu$ on $\S$ introduced in \eqref{bisb14}. Let $\de>0$ such that:
\begin{equation}\label{re12}
0<\sqrt{\ep} <\!\!< \de <\!\!< 1.
\end{equation}
We now require that the support of $\eta^\nu_j$ is a patch on $\S$ of diameter $\sim \de 2^{-j/2}$. We have:
\begin{equation}\label{re13}
\begin{array}{ll}
\ds\sum_{|j-l|\leq 2}\int_{\s}U_jf(x)\overline{U_lf(x)}d\s= & \ds\sum_{|j-l|\leq 2}\sum_{|\nu-\nu'|\leq 2\de 2^{-j/2}}\int_{\s}U_j^\nu f(x)\overline{U_l^{\nu'}f(x)}d\s\\& \ds +\sum_{|j-l|\leq 2}\sum_{|\nu-\nu'|> 2\de 2^{-j/2}}\int_{\s}U_j^\nu f(x)\overline{U_l^{\nu'}f(x)}d\s.
\end{array}
\end{equation}
The proof of Proposition \ref{bisorthoangle}  yields:
\begin{equation}\label{re14}
\left|\sum_{|\nu-\nu'|> 2\de 2^{-j/2}}\int_{\s}U_j^\nu f(x)\overline{U_l^{\nu'}f(x)}d\s\right|\lesssim \frac{\ep}{\de^2}\gamma_j\gamma_l,
\end{equation}
where $\gamma_j, \gamma_l$ have been defined in \eqref{bisdecf}. Indeed, \eqref{re14} follows from the equivalent of the two key estimates \eqref{bisoa11} \eqref{bisoa12}. For example, let us consider the equivalent of \eqref{bisoa11}. We obtain:
\begin{equation}\label{re14:1}
\begin{array}{r}
\ds\left|\int_{\s}U^{\nu,k}_jf(x)\overline{U^{\nu',k}_jf(x)}d\s \right|\lesssim \frac{\de\ep\ga_j^{\nu,k}\ga_j^{\nu',k}}{2^{j\a/2}(2^{j/2}|\nu-\nu'|)^{2-\a}}+\frac{\de\ep\ga_j^{\nu,k}\ga_j^{\nu',k}}{(2^{j/2}|\nu-\nu'|)^{3}}\\
\ds\textrm{for }|\nu-\nu'|\neq 0,\,1\leq k\leq |\nu-\nu'|^{-\a},
\end{array}
\end{equation}
where $\ep$ comes from the fact that $b$ satisfies \eqref{re4}, and $\de$ from the fact that the square root of the volume of the support of $\eta^\nu_j$ now yields $\de 2^{-j/2}$ instead of $2^{-j/2}$. The worst term in the right-hand side of \eqref{re14:1} is the second one. It may be rewritten:
\begin{equation}\label{re14:2}
\ds\frac{\de\ep\ga_j^{\nu,k}\ga_j^{\nu',k}}{(2^{j/2}|\nu-\nu'|)^{3}}=\frac{\ep}{\de^2}\frac{\ga_j^{\nu,k}\ga_j^{\nu',k}}{(2^{j/2}\de^{-1}|\nu-\nu'|)^{3}},
\end{equation}
and yields the factor $\ep\de^{-2}$ in the right-hand side of \eqref{re14}.

Finally, \eqref{re13} and \eqref{re14} yield:
\begin{equation}\label{re15}
\begin{array}{ll}
\ds\sum_{|j-l|\leq 2}\int_{\s}U_jf(x)\overline{U_lf(x)}d\s= & \ds\sum_{|j-l|\leq 2}\sum_{|\nu-\nu'|\leq 2\de 2^{-j/2}}\int_{\s}U_j^\nu f(x)\overline{U_l^{\nu'}f(x)}d\s\\
&\ds +O\left(\frac{\ep}{\de^2}\right)\norm{f}^2_{\le{2}}.
\end{array}
\end{equation}

\begin{remark}
The sum over $|\nu-\nu'|\leq 2\de 2^{-j/2}$ in the right-hand side of \eqref{re15} corresponds to the terms such that the support of $\eta_j^\nu$ and the support of $\eta_j^{\nu'}$ have a non empty intersection. The number of terms in this sum only depends on the dimension of $\S$ and is therefore a universal constant.
\end{remark}

\subsubsection{Step 3: control of the diagonal term}

The goal of this section is to estimate the term $\sum_{|j-l|\leq 2}\sum_{|\nu-\nu'|\leq 2\de 2^{-j/2}}\int_{\s}U_j^\nu f(x)\overline{U_l^{\nu'}f(x)}d\s$. 

\vspace{0.2cm}

\noindent{\bf A first reduction.} Remark first that the proof of Proposition \ref{bisdiagonal} together with the fact that $b$ satisfies \eqref{re3} immediately yields:
\begin{equation}\label{re16}
\norm{U_j^\nu f}_{\ll{2}}\lesssim\gamma_j^\nu.
\end{equation}
We introduce the operator $S_j^\nu$ defined on $\s$ by:
\begin{equation}\label{re17}
S_j^\nu f(x)=\int_{\S}\int_{0}^{+\infty}e^{i\la u(x,\o)}\psi(2^{-j}\la)\eta_j^\nu(\o) f(\la\o)\la^2 d\la d\o.
\end{equation}
By Proposition \ref{bisdiprop1}, we have:
\begin{equation}\label{re18}
\norm{S_j^\nu f}_{\ll{2}}\lesssim\gamma_j^\nu.
\end{equation}
The estimate \eqref{bisdi7} together with the assumption \eqref{re3} on $b$, \eqref{re16}, \eqref{re18} and the fact that $|\o-\nu|\lesssim\de 2^{-j/2}$ on the support of $\eta_j^\nu$ yields:
\begin{equation}\label{re19}
\begin{array}{ll}
&\ds\sum_{|j-l|\leq 2}\sum_{|\nu-\nu'|\leq 2\de 2^{-j/2}}\int_{\s}U_j^\nu f(x)\overline{U_l^{\nu'}f(x)}d\s\\
\ds = &\ds\sum_{|j-l|\leq 2}\sum_{|\nu-\nu'|\leq 2\de 2^{-j/2}}\int_{\s}b(x,\nu)S_j^\nu f(x)\overline{b(x,\nu')S_l^{\nu'}f(x)}d\s+O(\de)\norm{f}^2_{\le{2}},
\end{array}
\end{equation}
which together with the assumption \eqref{re4} on $b$ and \eqref{re18} implies:
\begin{equation}\label{re20}
\begin{array}{ll}
&\ds\sum_{|j-l|\leq 2}\sum_{|\nu-\nu'|\leq 2\de 2^{-j/2}}\int_{\s}U_j^\nu f(x)\overline{U_l^{\nu'}f(x)}d\s\\
= &\ds\sum_{|j-l|\leq 2}\sum_{|\nu-\nu'|\leq 2\de 2^{-j/2}}\int_{\s}S_j^\nu f(x)\overline{S_l^{\nu'}f(x)}d\s+O(\de+\ep)\norm{f}^2_{\le{2}}.
\end{array}
\end{equation}
We want to estimate the term $\sum_{|j-l|\leq 2}\sum_{|\nu-\nu'|\leq 2\de 2^{-j/2}}\int_{\s}U_j^\nu f(x)\overline{U_l^{\nu'}f(x)}d\s$. In view of \eqref{re20}, we may estimate instead the term $\sum_{|j-l|\leq 2}\sum_{|\nu-\nu'|\leq 2\de 2^{-j/2}}\int_{\s}S_j^\nu f(x)\overline{S_l^{\nu'}f(x)}d\s$. 

\vspace{0.2cm}

\noindent{\bf End of the proof of Proposition \ref{re2}.} Recall {\bf Assumption 4} which states that the map $\phi_\nu:\s\rightarrow\R^3$ defined by:
\begin{equation}\label{re21}
\phi_\nu(x):=u(x,\nu)\nu+\po u(x,\nu),
\end{equation}
is a bijection, such that the determinant of its Jacobian satisfies the following estimate:
\begin{equation}\label{re22}
\norm{|\det(\textrm{Jac}\phi_\nu)|-1}_{\ll{\infty}}\lesssim\ep.
\end{equation}
Let us note $\f$ the inverse Fourier transform on $\R^3$. We introduce the operator $\widetilde{S}_j^\nu$ on $\s$ defined by:
\begin{equation}\label{re24}
\widetilde{S}_j^\nu f(x)=\f(\psi(2^{-j}\c)\eta_j^\nu f)(\phi_\nu(x))=\int_{\R^3}e^{i\la\phi_\nu(x)\c\o}\psi(2^{-j}\la)\eta_j^\nu(\o)f(\la\o)\la^2 d\la d\o.
\end{equation}
The following proposition shows that the term $\int_{\s}S_j^\nu f(x)\overline{S_l^{\nu'}f(x)}d\s$ is close to 
the term $\int_{\s}\widetilde{S}_j^\nu f(x)\overline{\widetilde{S}_l^{\nu'}f(x)}d\s$.
\begin{proposition}\label{re25}
We have the following bound:
\begin{equation}\label{re26}
\norm{S_j^\nu f-\widetilde{S}_j^\nu f}_{\ll{2}}\lesssim\de^{\frac{1}{2}}\gamma_j^\nu.
\end{equation}
\end{proposition}
We postponed the proof of Proposition \ref{re25} to the next section. Let us show how Proposition \ref{re25} allows us to conclude the proof of Proposition \ref{re2}. \eqref{re10}, \eqref{re15}, \eqref{re20} and \eqref{re26} yield:
\begin{equation}\label{re27}
\norm{Uf}_{\ll{2}}^2=\sum_{|j-l|\leq 2}\sum_{|\nu-\nu'|\leq 2\de 2^{-j/2}}\int_{\s}\widetilde{S}_j^\nu f(x)\overline{\widetilde{S}_l^{\nu'}f(x)}d\s+O\left(\frac{\ep}{\de^2}+\de^{\frac{1}{2}}\right)\norm{f}^2_{\le{2}}.
\end{equation}
Making the change of variable $y=\phi_\nu(x)$ in $\int_{\s}\widetilde{S}_j^\nu f(x)\overline{\widetilde{S}_l^{\nu'}f(x)}d\s$ and using \eqref{re22} and \eqref{re24} implies:
\begin{equation}\label{re28}
\begin{array}{ll}
&\ds\sum_{|j-l|\leq 2}\sum_{|\nu-\nu'|\leq 2\de 2^{-j/2}}\int_{\s}\widetilde{S}_j^\nu f(x)\overline{\widetilde{S}_l^{\nu'}f(x)}d\s\\
\ds = &\ds\sum_{|j-l|\leq 2}\sum_{|\nu-\nu'|\leq 2\de 2^{-j/2}}\int_{\R^3}\f(\psi(2^{-j}\c)\eta_j^\nu f)(y)\overline{\f(\psi(2^{-l}\c)\eta_j^{\nu'} f)(y)}dy\\
&\ds +O(\ep)\norm{f}^2_{\le{2}}\\
\ds = &\ds\sum_{|j-l|\leq 2}\sum_{|\nu-\nu'|\leq 2\de 2^{-j/2}}\int_{\R^3}\psi(2^{-j}\la)\eta_j^\nu(\o) f(\la\o)\overline{\psi(2^{-l}\la)\eta_j^{\nu'}(\o) f(\la\o)}dy\\
&\ds +O(\ep)\norm{f}^2_{\le{2}},
\end{array}
\end{equation}
where we have used the fact that $\f$ is an isomorphism on $\le{2}$ in the last equality of \eqref{re28}. 
Now, we have:
\begin{equation}\label{re29}
\ds\sum_{|j-l|\leq 2}\sum_{|\nu-\nu'|\leq 2\de 2^{-j/2}}\int_{\R^3}\psi(2^{-j}\la)\eta_j^\nu(\o) f(\la\o)\overline{\psi(2^{-l}\la)\eta_j^{\nu'}(\o) f(\la\o)}dy=\norm{f}^2_{\le{2}},
\end{equation}
which together with \eqref{re27} and \eqref{re28} yields:
\begin{equation}\label{re30}
\norm{Uf}_{\ll{2}}^2=\norm{f}_{\le{2}}^2+O\left(\frac{\ep}{\de^2}+\de^{\frac{1}{2}}\right)\norm{f}^2_{\le{2}}.
\end{equation}
Choosing $\de^{\frac{1}{2}}$ and $\ep\de^{-2}$ small enough, we deduce from \eqref{re30}:
\begin{equation}\label{re31}
\norm{f}_{\le{2}}\lesssim\norm{Uf}_{\ll{2}},
\end{equation}
which is the wanted estimate. This conclude the proof of Proposition \ref{re2}. \QED

\subsubsection{Proof of Proposition \ref{re25}}

\vspace{0.2cm}

\noindent{\bf Reduction to a decay estimate.} Relying on the classical $TT^*$ argument, \eqref{re26} is equivalent to proving the boundedness 
in $\ll{2}$ with a norm $O(\de)$ of the operator whose kernel $K$ is given by:
\begin{equation}\label{re32}
\begin{array}{ll}
\ds K(x,y)= & \ds\int_{\S}\int_{0}^{+\infty}e^{i\la u(x,\o)-i\la u(y,\o)}\psi(2^{-j}\la)\eta_j^\nu(\o)\la^2 d\la d\o\\
&\ds +\int_{\S}\int_{0}^{+\infty}e^{i\la\phi_\nu(x)\c\o-i\la \phi_\nu(y)\c\o}\psi(2^{-j}\la)\eta_j^\nu(\o)\la^2 d\la d\o\\
&\ds -\int_{\S}\int_{0}^{+\infty}e^{i\la u(x,\o)-i\la \phi_\nu(y)\c\o}\psi(2^{-j}\la)\eta_j^\nu(\o)\la^2 d\la d\o\\
&\ds -\int_{\S}\int_{0}^{+\infty}e^{i\la \phi_\nu(x)\c\o-i\la u(y,\o)}\psi(2^{-j}\la)\eta_j^\nu(\o)\la^2 d\la d\o.
\end{array}
\end{equation}
\eqref{re26} then reduces to proving the following decay for the kernel $K$ in \eqref{re32}:
\begin{equation}\label{re33}
\begin{array}{ll}
\ds |K(x,y)|\lesssim & \ds\de^{\frac{1}{2}}\frac{2^j}{(1+|2^j|u(x,\nu)-u(y,\nu)|-2^{j/2}|\po u(x,\nu)-\po u(y,\nu)||)^2}\\
& \ds\times\frac{2^j}{(1+2^{j/2}|\po u(x,\nu)-\po u(y,\nu)|)^3}.
\end{array}
\end{equation}
The proof of the fact that \eqref{re33} implies \eqref{re26} is identical to the proof in section \ref{sec:bisdiprop1} of the fact that the decay estimate \eqref{di13} implies \eqref{bisdiprop2}. In fact, 
performing the exact same changes of variables leads to:
\begin{equation}\label{re33b}
\sup_{x\in\s}\int_\s |K(x,y)|dy\lesssim\de,\,\sup_{y\in\s}\int_\s |K(x,y)|dx\lesssim\de.
\end{equation}
Finally, \eqref{re33b} yields \eqref{re26} by Schur's Lemma.

\vspace{0.2cm}

\noindent{\bf Proof of the decay estimate \eqref{re33}.} The proof of \eqref{re33} follows from the proof 
of Proposition \ref{di12} in section \ref{sec:di20}. In fact, let us consider the following quantity $A$ defined by:
\begin{equation}\label{re34}
\ds A=\ds\int_{\S}\int_{0}^{+\infty}e^{i\la\rho(\o)}\psi(2^{-j}\la)\eta_j^\nu(\o)\la^2 d\la d\o
\end{equation}
where $\rho$ is a function defined on $\S$. Then applying 3 integrations by parts with respect to 
$\o$ and 2 integrations by parts with respect to $\la$ as in the proof of Proposition \ref{di12} yields to the following equality:
\begin{equation}\label{re35}
\begin{array}{rr}
\ds A= &\ds\int_{\S}\int_{0}^{+\infty}F_0(2^j\rho(\o),2^{j/2}\po\rho(\o),\po\rho(\o),\po^2\rho(\o),\po^3\rho(\o),2^{-j}\la)\\
&\ds\times\psi_0(2^{-j}\la)\eta_j^\nu(\o)\la^2 d\la d\o\\
&\ds +\de^{-1}\int_{\S}\int_{0}^{+\infty}F_1(2^j\rho(\o),2^{j/2}\po\rho(\o),\po\rho(\o),\po^2\rho(\o),\po^3\rho(\o),2^{-j}\la)\\
&\ds\times\psi_1(2^{-j}\la)(\de 2^{-j/2}\po)\eta_j^\nu(\o)\la^2 d\la d\o\\
& +\ds\de^{-2}\int_{\S}\int_{0}^{+\infty}F_2(2^j\rho(\o),2^{j/2}\po\rho(\o),\po\rho(\o),\po^2\rho(\o),2^{-j}\la)\\
&\ds\times\psi_2(2^{-j}\la)(\de 2^{-j/2}\po)^{2}\eta_j^\nu(\o)\la^2 d\la d\o\\
& +\ds\de^{-3}\int_{\S}\int_{0}^{+\infty}F_3(2^j\rho(\o),2^{j/2}\po\rho(\o),\po\rho(\o),2^{-j}\la)\\
&\ds\times\psi_3(2^{-j}\la)(\de 2^{-j/2}\po)^{3}\eta_j^\nu(\o)\la^2 d\la d\o,
\end{array}
\end{equation}
where $\psi_l, l=0, 1, 2,3$ is smooth and compactly supported in $(0,+\infty)$, $(\de 2^{-j/2}\po)^{l}\eta^\nu_j, l=0, 1, 2,3$ is bounded on $\S$ and has the same support as $\eta^\nu_j$, and $F_l, l=0,1 ,2 ,3$ are smooth function satisfying the following estimates:
\begin{equation}\label{re36}
\begin{array}{ll}
&\ds |F_0(2^j\rho(\o),2^{j/2}\po\rho(\o),\po\rho(\o),\po^2\rho(\o),\po^3\rho(\o),2^{-j}\la)|\\
&\ds +|F_1(2^j\rho(\o),2^{j/2}\po\rho(\o),\po\rho(\o),\po^2\rho(\o),\po^3\rho(\o),2^{-j}\la)|\\
&\ds +|F_2(2^j\rho(\o),2^{j/2}\po\rho(\o),\po\rho(\o),\po^2\rho(\o),2^{-j}\la)|\\
&\ds+ |F_3(2^j\rho(\o),2^{j/2}\po\rho(\o),\po\rho(\o),2^{-j}\la)|\\
\ds\lesssim &\ds\frac{1}{(1+2^j|\rho(\o)|)^2(1+2^{j/2}|\po\rho(\o)|)^3}.
\end{array}
\end{equation}
Indeed, this has been done in the proof of Proposition \ref{di12} for the particular case $\rho(\o)=u(x,\o)-u(y,\o)$ but is easily seen to hold in the general case with the exact same proof. Applying \eqref{re34} to the 4 terms in the right-hand side of \eqref{re32} respectively with 
\begin{equation}\label{re38}
\begin{array}{l}
\ds\rho_1(\o)=u(x,\o)-u(y,\o),\, \rho_2(\o)=\phi_\nu(x)\c\o-\phi_\nu(y)\c\o,\\ 
\ds\rho_3(\o)=u(x,\o)-\phi_\nu(y)\c\o\textrm{ and }\rho_4(\o)=\phi_\nu(y)\c\o-u(y,\o)
\end{array}
\end{equation}
 yields:
\begin{equation}\label{re39}
\begin{array}{rr}
\ds K(x,y)= & \ds\sum_{q=1}^2\sum_{l=0}^3\de^{-l}\int_{\S}\int_{0}^{+\infty}F_l[\rho_q]\psi_l(2^{-j}\la)(\de 2^{-j/2}\po)^{l}\eta_j^\nu(\o)\la^2 d\la d\o\\
& -\ds\sum_{q=3}^4\sum_{l=0}^4\de^{-l}\int_{\S}\int_{0}^{+\infty}F_l[\rho_q]\psi_l(2^{-j}\la)(\de 2^{-j/2}\po)^{l}\eta_j^\nu(\o)\la^2 d\la d\o,
\end{array}
\end{equation}
where $F_l[\rho_q]$ is defined for $q=1, 2, 3, 4$ by:
\begin{equation}\label{re40}
\begin{array}{l}
\ds F_l[\rho_q]=F_l(2^j\rho_q(\o),2^{j/2}\po\rho_q(\o),\po\rho_q(\o),\po^2\rho_q(\o),\po^3\rho_q(\o),2^{-j}\la),\,l=0, 1,\\
F_2[\rho_q]=F_2(2^j\rho_q(\o),2^{j/2}\po\rho_q(\o),\po\rho_q(\o),\po^2\rho_q(\o),2^{-j}\la),\\
F_3[\rho_q]=F_3(2^j\rho_q(\o),2^{j/2}\po\rho_q(\o),\po\rho_q(\o),2^{-j}\la).
\end{array}
\end{equation}
We rewrite \eqref{re39} as:
\begin{equation}\label{re41}
\begin{array}{rr}
\ds K(x,y)= &\ds\sum_{l=0}^3\de^{-l}\int_{\S}\int_{0}^{+\infty}(F_l[\rho_1]-F_l[\rho_3])\psi_l(2^{-j}\la)(\de 2^{-j/2}\po)^{l}\eta_j^\nu(\o)\la^2 d\la d\o\\
& +\ds\sum_{l=0}^3\de^{-l}\int_{\S}\int_{0}^{+\infty}(F_l[\rho_2]-F_l[\rho_4])\psi_l(2^{-j}\la)(\de 2^{-j/2}\po)^{l}\eta_j^\nu(\o)\la^2 d\la d\o.
\end{array}
\end{equation}
We now estimate $F_l[\rho_1]-F_l[\rho_3]$ and $F_l[\rho_2]-F_l[\rho_4]$. One 
easily checks that the first order derivatives of $F_l$ satisfy the same estimate as the estimates \eqref{re36} satisfied by $F_l$. Together with {\bf Assumption 5} on $u(x,\o)-\phi_\nu(x)\c\o$, \eqref{di24} and \eqref{di25}, we deduce the following estimates on the support of $\eta_j^\nu$:
\begin{equation}\label{re42}
\begin{array}{l}
\ds |F_0[\rho_1]-F_0[\rho_3]|+|F_0[\rho_2]-F_0[\rho_4]|+|F_1[\rho_1]-F_1[\rho_3]|+|F_1[\rho_2]-F_1[\rho_4]|\\
\ds\leq \ds |F_0[\rho_1]|+|F_0[\rho_2]|+|F_0[\rho_3]|+|F_0[\rho_4]|+|F_1[\rho_1]|+|F_1[\rho_2]|+|F_1[\rho_3]|+|F_1[\rho_4]|\\
\ds\lesssim \frac{1}{(1+|2^j|u(x,\nu)-u(y,\nu)|-2^{j/2}|\po u(x,\nu)-\po u(y,\nu)||)^2}\\
\ds\times\frac{1}{(1+2^{j/2}|\po u(x,\nu)-\po u(y,\nu)|)^3},
\end{array}
\end{equation}
\begin{equation}\label{re43}
\begin{array}{l}
|F_2[\rho_1]-F_2[\rho_3]|+|F_2[\rho_2]-F_2[\rho_4]|\\
\ds\lesssim \big(2^j|u(x,\o)-\phi_\nu(x)\c\o|+2^j|u(y,\o)-\phi_\nu(y)\c\o|\\
\ds +2^{j/2}|\po u(x,\o)-\po(\phi_\nu(x)\c\o)|+2^{j/2}|\po u(y,\o)-\po(\phi_\nu(y)\c\o)|\\
\ds +|\po^2u(x,\o)-\po^2(\phi_\nu(x)\c\o)|+|\po^2u(y,\o)-\po^2(\phi_\nu(y)\c\o)|\big)\\
\ds\times\frac{1}{(1+|2^j|u(x,\nu)-u(y,\nu)|-2^{j/2}|\po u(x,\nu)-\po u(y,\nu)||)^2}\\
\ds\times\frac{1}{(1+2^{j/2}|\po u(x,\nu)-\po u(y,\nu)|)^3}\\
\ds\lesssim \ds\de\frac{1}{(1+|2^j|u(x,\nu)-u(y,\nu)|-2^{j/2}|\po u(x,\nu)-\po u(y,\nu)||)^2}\\
\ds\times\frac{1}{(1+2^{j/2}|\po u(x,\nu)-\po u(y,\nu)|)^3},
\end{array}
\end{equation}
and
\begin{equation}\label{re44}
\begin{array}{l}
|F_3[\rho_1]-F_3[\rho_3]|+|F_3[\rho_2]-F_3[\rho_4]|\\
\ds\lesssim \big(2^j|u(x,\o)-\phi_\nu(x)\c\o|+2^j|u(y,\o)-\phi_\nu(y)\c\o|\\
\ds +2^{j/2}|\po u(x,\o)-\po(\phi_\nu(x)\c\o)|+2^{j/2}|\po u(y,\o)-\po(\phi_\nu(y)\c\o)|\big)\\
\ds\times\frac{1}{(1+|2^j|u(x,\nu)-u(y,\nu)|-2^{j/2}|\po u(x,\nu)-\po u(y,\nu)||)^2}\\
\ds\times\frac{1}{(1+2^{j/2}|\po u(x,\nu)-\po u(y,\nu)|)^3}\\
\ds\lesssim \ds\de^2\frac{1}{(1+|2^j|u(x,\nu)-u(y,\nu)|-2^{j/2}|\po u(x,\nu)-\po u(y,\nu)||)^2}\\
\ds\times\frac{1}{(1+2^{j/2}|\po u(x,\nu)-\po u(y,\nu)|)^3},
\end{array}
\end{equation}
where we have used the fact that $|\o-\nu|\lesssim \de 2^{-j/2}$ on the support of $\eta_j^\nu$, and $\ep\lesssim \delta$ in view of \eqref{re12}. Now, we have:
\begin{equation}\label{re45}
\begin{array}{ll}
&\ds\int_{\S}\int_{0}^{+\infty}\psi_l(2^{-j}\la)(\de 2^{-j/2}\po)^{l}\eta_j^\nu(\o)\la^2 d\la d\o\\[3mm]
\ds = &\ds\left(\int_{0}^{+\infty}\psi_l(2^{-j}\la)\la^2 d\la\right)
\left(\int_{\S}(\de 2^{-j/2}\po)^{l}\eta_j^\nu(\o)d\o\right)\lesssim \de^2 2^{2j},
\end{array}
\end{equation}
where we have used the fact that $(\de 2^{-j/2}\po)^{l}\eta_j^\nu(\o)$ is bounded on $\S$ and the fact that its support is two dimensional with a diameter $\sim \de 2^{-j/2}$. \eqref{re41}-\eqref{re45} immediately yield the decay estimate \eqref{re33}. Finally, as explained after \eqref{re33}, \eqref{re33} yields \eqref{re33b} which implies \eqref{re26}. This concludes the proof of Proposition \ref{re25}. \QED

\begin{remark}
Note that {\bf Assumption 5} does not contain any estimate for the term $\po^3u-\po^3(\phi_\nu(x)\cdot\o)$. Instead, this term is estimated using  {\bf Assumption 2}:
$$|\po^3u-\po^3(\phi_\nu(x)\cdot\o)|\lesssim 1,$$
and thus is not bounded from above by $O(\ep)$ unlike the corresponding estimate for $\po^2u-\po^2(\phi_\nu(x)\cdot\o)$ in {\bf Assumption 5}. As a consequence, \eqref{re42} cannot be improved, and is responsible for the  introduction of the extra smallness parameter $\delta$ in the decomposition in angle.
\end{remark}

\subsection{Proof of the estimate \eqref{choicef4bis}}

Recall the definition of the Fourier integral operators $M_\pm$ and $Q_\pm$ introduced in section \ref{sec:choicef}: 
\begin{equation}\label{et1}
M_\pm f(x)=\int_{\S}\int_{0}^{+\infty}e^{\pm i\lambda u(x,\pm\o)}f(\lambda\o)\lambda^2 d\lambda d\o,
\end{equation}
and 
\begin{equation}\label{et2}
Q_\pm f(x)=\int_{\S}\int_{0}^{+\infty}e^{\pm i\lambda u(x,\pm\o)}a(x,\pm\o)^{-1} f(\lambda\o)\lambda^2 d\lambda d\o.
\end{equation}
Let $(f_+,f_-)$ satisfying:
\begin{equation}\label{et3}
\left\{\begin{array}{l}
M_+f_++M_-f_-=\phi_0,\\
Q_+(\la f_+)-Q_-(\la f_-)=i\phi_1.
\end{array}\right.
\end{equation}
The goal of this section is to prove that $(f_+,f_-)$ satisfies the following estimate:
\begin{equation}\label{et4}
\norm{\la f_+}_{L^2(\R^3)}+\norm{\la f_-}_{L^2(\R^3)}\lesssim \norm{\nabla\phi_0}_{\ll{2}}+\norm{\phi_1}_{\ll{2}}.
\end{equation}
Using Proposition \ref{re2} in the case of a symbol $b\equiv 1$, we obtain:
\begin{equation}\label{et5}
\norm{\la f_+}_{L^2(\R^3)}\lesssim \norm{M_+(\la f_+)}_{\ll{2}}\textrm{ and }\norm{\la f_-}_{L^2(\R^3)}\lesssim\norm{M_-(\la f_-)}_{\ll{2}},
\end{equation}
which yields:
\begin{equation}\label{et6}
\begin{array}{ll}
&\ds\norm{\la f_+}_{L^2(\R^3)}+\norm{\la f_-}_{L^2(\R^3)}\\
\ds\lesssim &\ds\norm{M_+(\la f_+)}_{\ll{2}}+\norm{M_-(\la f_-)}_{\ll{2}}\\
\ds\lesssim &\ds\norm{M_+(\la f_+)+M_-(\la f_-)}_{\ll{2}}+\norm{M_+(\la f_+)-M_-(\la f_-)}_{\ll{2}}.
\end{array}
\end{equation}
We have:
\begin{equation}\label{et7}
(Q_\pm-M_\pm) f(x)=\int_{\S}\int_{0}^{+\infty}e^{\pm i\lambda u(x,\pm\o)}(a(x,\pm\o)^{-1}-1) f(\lambda\o)\lambda^2 d\lambda d\o.
\end{equation}
Due to {\bf Assumption 1-3} on $a$, the symbol $a(x,\pm\o)^{-1}-1$ of $Q_\pm-M_\pm$ satisfies the assumptions \eqref{thregx1s}-\eqref{cordecfr1s} of Theorem \ref{th1} with $D=\ep$. Thus, we obtain from \eqref{l2} that:
\begin{equation}\label{et8}
\norm{(Q_\pm-M_\pm)f}_{\ll{2}}\lesssim\ep\norm{f}_{\le{2}}.
\end{equation}
The second equation of \eqref{et3}, \eqref{et6} and \eqref{et8} yield:
\begin{equation}\label{et9}
\ds\norm{\la f_+}_{L^2(\R^3)}+\norm{\la f_-}_{L^2(\R^3)}\lesssim \norm{M_+(\la f_+)+M_-(\la f_-)}_{\ll{2}}+\norm{\phi_1}_{\ll{2}}.
\end{equation}
The following lemma will allow us to bound the first term in the right-hand side of \eqref{et9}.
\begin{lemma}\label{et10}
For any $(f_+,f_-)$, we have the following bound: 
\begin{equation}\label{et11}
\begin{array}{ll}
\ds\norm{M_+(\la f_+)+M_-(\la f_-)}_{\ll{2}}\lesssim & \ds\norm{\nabla M_+(f_+)+\nabla M_-( f_-)}_{\ll{2}}\\
& \ds +\left(\de+\frac{\ep}{\de^2}\right)(\norm{\la f_+}_{\le{2}}+\norm{\la f_-}_{\le{2}}),
\end{array}
\end{equation}
where $\de$ may be chosen as in \eqref{re12}.
\end{lemma}
Before proving Lemma \ref{et10}, we first conclude the proof of the estimate \eqref{choicef4bis}. \eqref{re12}, \eqref{et9} and \eqref{et11} yield:
\begin{equation}\label{et12}
\ds\norm{\la f_+}_{L^2(\R^3)}+\norm{\la f_-}_{L^2(\R^3)}\lesssim \norm{\nabla M_+(f_+)+\nabla M_-( f_-)}_{\ll{2}}+\norm{\phi_1}_{\ll{2}}.
\end{equation}
Applying $\nabla$ to the first equation of \eqref{et3} and using \eqref{et12} implies:
\begin{equation}\label{et13}
\norm{\la f_+}_{L^2(\R^3)}+\norm{\la f_-}_{L^2(\R^3)}\lesssim \norm{\nabla\phi_0}_{\ll{2}}+\norm{\phi_1}_{\ll{2}},
\end{equation}
which is the wanted estimate \eqref{choicef4bis}. 

\vspace{0.2cm}

\noindent{\bf Proof of Lemma \ref{et10}.} Since $\nabla u=a^{-1}N$, we have:
\begin{equation}\label{et14}
\nabla M_\pm f(x)=\pm i\int_{\S}\int_{0}^{+\infty}e^{\pm i\lambda u(x,\pm\o)}a(x,\pm\o)^{-1}N(x,\pm\o)\la f(\lambda\o)\lambda^2 d\lambda d\o.
\end{equation}
We introduce the operator $P_\pm$ defined by:
\begin{equation}\label{et15}
P_\pm f(x)=\int_{\S}\int_{0}^{+\infty}e^{\pm i\lambda u(x,\pm\o)}N(x,\pm\o) f(\lambda\o)\lambda^2 d\lambda d\o.
\end{equation}
Due to {\bf Assumption 1-3} on $a$ and $N$, the symbol $i\pm(a(x,\pm\o)^{-1}-1)N_\pm$ of $\nabla M_\pm\mp i P_\pm(\la .)$ satisfies the assumptions \eqref{thregx1s}-\eqref{cordecfr1s} of Theorem \ref{th1} with $D=\ep$. Thus, we obtain from \eqref{l2} that:
\begin{equation}\label{et16}
\norm{\nabla M_\pm(f)\mp i P_\pm(\la f)}_{\ll{2}}\lesssim\ep\norm{\la f}_{\le{2}}.
\end{equation}
Thus, the proof of Lemma \ref{et10} reduces to the proof of the following estimate:
\begin{equation}\label{et17}
\begin{array}{ll}
\ds\norm{M_+(f_+)+M_-(f_-)}_{\ll{2}}\lesssim & \ds\norm{P_+(f_+)-P_-( f_-)}_{\ll{2}}\\
& \ds +\left(\de+\frac{\ep}{\de^2}\right)(\norm{f_+}_{\le{2}}+\norm{f_-}_{\le{2}}),
\end{array}
\end{equation}
for any $(f_+,f_-)$ in $\le{2}\times\le{2}$. To prove \eqref{et17}, we decompose in frequency and angle 
as in the proof of Proposition \ref{re2}, in order to reduce ourselves to diagonal terms.

\vspace{0.2cm}

\noindent{\bf Decomposition in frequency.} As in step 1 of the proof of Theorem \ref{th1}, we decompose $M_\pm(f_\pm)$ and $P_\pm(f_\pm)$ in frequency:
\begin{equation}\label{et18}
M_\pm(f_\pm)(x)=\sum_{j\geq -1}(M_{\pm})_jf_\pm(x)\textrm{ and }P_\pm(f_\pm)(x)=\sum_{j\geq -1}(P_{\pm})_jf_\pm(x)
\end{equation}
where the operators $(M_\pm)_j$, $(P_\pm)_j$ are defined as in \eqref{bisb9} \eqref{bisb10}. Following step 1 of the proof of Proposition \ref{re2}, we obtain the equivalent of \eqref{re10}:
\begin{equation}\label{et19}
\begin{array}{ll}
&\ds\norm{M_+(f_+)+M_-(f_-)}^2_{\ll{2}}\\
\ds = &\ds\sum_{|j-l|\leq 2}\int_{\s}((M_+)_jf_+(x)+(M_-)_jf_-(x))\overline{((M_+)_lf_+(x)+(M_-)_lf_-(x))}d\s\\
&\ds +O(\ep)(\norm{f_+}_{\le{2}}^2+\norm{f_-}_{\le{2}}^2),
\end{array}
\end{equation}
and 
\begin{equation}\label{et20}
\begin{array}{ll}
&\ds\norm{P_+(f_+)-P_-(f_-)}^2_{\ll{2}}\\
\ds = &\ds\sum_{|j-l|\leq 2}\int_{\s}((P_+)_jf_+(x)-(P_-)_jf_-(x))\c\overline{((P_+)_lf_+(x)-(P_-)_lf_-(x))}d\s\\
&\ds +O(\ep)(\norm{f_+}_{\le{2}}^2+\norm{f_-}_{\le{2}}^2).
\end{array}
\end{equation}

\vspace{0.2cm}

\noindent{\bf Decomposition in angle.} As in step 2 of the proof of Proposition \ref{re2}, we decompose $(M_\pm)_jf_\pm$ and $(P_\pm)_jf_\pm$ in angle:
\begin{equation}\label{et21}
(M_\pm)_jf_\pm (x)=\sum_{\nu\in\Gamma}(M_\pm)^\nu_jf(x)\textrm{ and }(P_\pm)_jf_\pm (x)=\sum_{\nu\in\Gamma}(P_\pm)^\nu_jf(x),
\end{equation}
where the operators $(M_\pm)_j^\nu$ and $(P_\pm)_j^\nu$ are defined as in \eqref{bisb16} and  
 where  the support of our partition of unity $\eta_j^\nu$ on $\S$ is  a patch of diameter $\sim \de 2^{-j/2}$ 
  with $\de$ chosen as in \eqref{re12}. Following step 2 of the proof of Proposition \ref{re2}, we obtain the equivalent of \eqref{re15}:
\begin{equation}\label{et22}
\begin{array}{ll}
&\ds\sum_{|j-l|\leq 2}\int_{\s}((M_+)_jf_+(x)+(M_-)_jf_-(x))\overline{((M_+)_lf_+(x)+(M_-)_lf_-(x))}d\s\\
\ds = &\ds\sum_{|j-l|\leq 2}\sum_{|\nu-\nu'|\leq 2\de 2^{-j/2}}\int_{\s}((M_+)^\nu_jf_+(x)+(M_-)^{\nu}_jf_-(x))\\
&\ds\times \overline{((M_+)^{\nu'}_lf_+(x)+(M_-)^{\nu'}_lf_-(x))}d\s+O\left(\frac{\ep}{\de^2}\right)(\norm{f_+}_{\le{2}}^2+\norm{f_-}_{\le{2}}^2),
\end{array}
\end{equation}
and
\bea\label{et23}
&&\sum_{|j-l|\leq 2}\int_{\s}((P_+)_jf_+(x)-(P_-)_jf_-(x))\c\overline{((P_+)_lf_+(x)-(P_-)_lf_-(x))}d\s\\
\nn& = &\ds\sum_{|j-l|\leq 2}\sum_{|\nu-\nu'|\leq 2\de 2^{-j/2}}\int_{\s}((P_+)^\nu_jf_+(x)-(P_-)^{\nu}_jf_-(x))\c\overline{((P_+)^{\nu'}_lf_+(x)-(P_-)^{\nu'}_lf_-(x))}d\s\\
\nn&&+O\left(\frac{\ep}{\de^2}\right)(\norm{f_+}_{\le{2}}^2+\norm{f_-}_{\le{2}}^2).
\eea

\vspace{0.2cm}

\noindent{\bf End of the proof of Lemma \ref{et10}.} \eqref{et19} and \eqref{et22} yield:
\begin{equation}\label{et24}
\begin{array}{r}
\ds\norm{M_+(f_+)+M_-(f_-)}_{\ll{2}}^2= \ds\sum_{|j-l|\leq 2}\sum_{|\nu-\nu'|\leq 2\de 2^{-j/2}}\int_{\s}((M_+)^\nu_jf_+(x)+(M_-)^{\nu}_jf_-(x))\\\ds\times \overline{((M_+)^{\nu'}_lf_+(x)+(M_-)^{\nu'}_lf_-(x))}d\s+O\left(\frac{\ep}{\de^2}\right)(\norm{f_+}_{\le{2}}^2+\norm{f_-}_{\le{2}}^2),
\end{array}
\end{equation}
and \eqref{et20} and \eqref{et23} yield:
\bea\label{et25}
&&\norm{P_+(f_+)-P_-(f_-)}_{\ll{2}}^2\\
\nn&=& \sum_{|j-l|\leq 2}\sum_{|\nu-\nu'|\leq 2\de 2^{-j/2}}\int_{\s}((P_+)^\nu_jf_+(x)-(P_-)^{\nu}_jf_-(x))\c \overline{((P_+)^{\nu'}_lf_+(x)-(P_-)^{\nu'}_lf_-(x))}d\s\\
\nn&&+O\left(\frac{\ep}{\de^2}\right)(\norm{f_+}_{\le{2}}^2+\norm{f_-}_{\le{2}}^2).
\eea
The operator $(P_\pm)^\nu_j-N(x,\pm\nu)(M_+)^\nu_j$ has a symbol given by $N(x,\pm\o)-N(x,\pm\nu)$. 
Thus, the estimate \eqref{bisdi7} together with {\bf Assumption 2} on $\po N$, and the fact that $|\o-\nu|\lesssim\de 2^{-j/2}$ on the support of $\eta_j^\nu$ yields:
\bea\label{et26}
\nn&&\sum_{|j-l|\leq 2}\sum_{|\nu-\nu'|\leq 2\de 2^{-j/2}}\int_{\s}((P_+)^\nu_jf_+(x)-(P_-)^{\nu}_jf_-(x))\c \overline{((P_+)^{\nu'}_lf_+(x)-(P_-)^{\nu'}_lf_-(x))}d\s\\
& = &\ds\sum_{|j-l|\leq 2}\sum_{|\nu-\nu'|\leq 2\de 2^{-j/2}}\int_{\s}(N(x,\nu)(M_+)^\nu_jf_+(x)-N(x,-\nu)(M_-)^{\nu}_jf_-(x))\\
\nn&&\ds\c \overline{(N(x,\nu)(M_+)^{\nu'}_lf_+(x)-N(x,-\nu)(M_-)^{\nu'}_lf_-(x))}d\s+O(\de)(\norm{f_+}_{\le{2}}^2+\norm{f_-}_{\le{2}}^2).
\eea
Now, {\bf Assumption 6} yields $|N(x,\nu)+N(x,-\nu)|\lesssim\ep$ which together with \eqref{et26} and the fact that $N$ is a unit vector implies:
\bea\label{et27}
\nn&&\ds\sum_{|j-l|\leq 2}\sum_{|\nu-\nu'|\leq 2\de 2^{-j/2}}\int_{\s}((P_+)^\nu_jf_+(x)-(P_-)^{\nu}_jf_-(x))\c \overline{((P_+)^{\nu'}_lf_+(x)-(P_-)^{\nu'}_lf_-(x))}d\s\\
\nn& = &\ds\sum_{|j-l|\leq 2}\sum_{|\nu-\nu'|\leq 2\de 2^{-j/2}}\int_{\s}((M_+)^\nu_jf_+(x)+
(M_-)^{\nu}_jf_-(x))\overline{((M_+)^{\nu'}_lf_+(x)+(M_-)^{\nu'}_lf_-(x))}d\s\\
&&\ds +O(\de+\ep)(\norm{f_+}_{\le{2}}^2+\norm{f_-}_{\le{2}}^2).
\eea
Finally, \eqref{et24}, \eqref{et25} and \eqref{et27} yield:
\begin{equation}\label{et28}
\begin{array}{ll}
\ds\norm{P_+(f_+)-P_-(f_-)}_{\ll{2}}^2= & \ds \norm{M_+(f_+)+M_-(f_-)}_{\ll{2}}^2 \\
&\ds +O\left(\de+\frac{\ep}{\de^2}\right)(\norm{f_+}_{\le{2}}^2+\norm{f_-}_{\le{2}}^2),
\end{array}
\end{equation}
which implies \eqref{et17}. As noticed at the beginning of the proof, \eqref{et17} yields the wanted 
estimate \eqref{et11}. This concludes the proof of Lemma \ref{et10}. \QED

\subsection{Existence of $(f_+,f_-)$}

In the previous section, we have proved the estimate \eqref{choicef4bis}:
\begin{equation}\label{su1}
\norm{\la f_+}_{L^2(\R^3)}+\norm{\la f_-}_{L^2(\R^3)}\lesssim \norm{\nabla\phi_0}_{\ll{2}}+\norm{\phi_1}_{\ll{2}},
\end{equation}
for any $(f_+,f_-)$ satisfying the following system:
\begin{equation}\label{su2}
\left\{\begin{array}{l}
M_+f_++M_-f_-=\phi_0,\\
Q_+(\la f_+)-Q_-(\la f_-)=i\phi_1.
\end{array}\right.
\end{equation}
Notice that \eqref{su1} implies the uniqueness of $(f_+,f_-)$ solution of \eqref{su2}. In this section, we 
complete the proof of Theorem \ref{th2} by proving the existence of $(f_+,f_-)$ solution of \eqref{su2}. 

Recall that the phase $u(x,\o)$ of our Fourier integral operators has been constructed in \cite{param1} on $\s\times\S$ under the assumption that $(\s,g,k)$ satisfies the following bounds consistent with the assumptions on $\Sigma$ for $R$ and $k$ in Theorem \ref{th:mainbl2}:
\be{l2bounds}
\norm{R}_{\ll{2}}\leq\ep,\,\norm{\nabla k}_{\ll{2}}\leq\ep.
\eeq
$(\s,g,k)$ also satisfies the constraint equations:
\be{const1}
\left\{\begin{array}{l}
\nabla^j k_{ij}=0,\\
 R=|k|^2,\\
 \textrm{Tr}k=0, 
\end{array}\right.
\eeq
where the last equation in \eqref{const1} comes from the fact that we work with a maximal foliation. We introduce two sets $V$ and $W$:
\be{defv}
V=\{(\s,g,k)\textrm{ such that \eqref{l2bounds} and \eqref{const1} are satisfied}\},
\eeq
and 
\be{defw}
\begin{array}{ll}
\ds W=\{ & \ds(\s,g,k)\in V\textrm{ such that }(f_+,f_-)\textrm{ solution of }\eqref{su2}\textrm{ exist for all }(\phi_0,\phi_1)\\
& \ds\textrm{ such that }\nabla\phi_0\in\ll{2}\textrm{ and }\phi_1\in\ll{2}\}.
\end{array}
\eeq
In order to prove the existence of $(f_+,f_-)$ solution of \eqref{su2}, we will show that $V=W$ by a connectedness argument. This will result from the following two lemmas.
\begin{lemma}\label{su3}
Let $N\geq 0$ an integer. Then, the set $V$ is connected for the topology of $(g,k)\in C^q(\s)\times C^{q-1}(\s)$.
\end{lemma}

\begin{lemma}\label{su4}
Let $N\geq 0$ an integer. Then, the set $W$ is open and closed in $V$ for the topology of $(g,k)\in C^q(\s)\times C^{q-1}(\s)$ provided $q$ is chosen sufficiently large.
\end{lemma}

\begin{remark}
The assumptions on the regularity on $(\s,g,k)$ in Lemma \ref{su3} and \ref{su4} are much stronger than the ones appearing in the bounded $L^2$ curvature conjecture. We would like to insist on the fact that this smoothness is only assumed to obtain the existence of $(f_+,f_-)$ solution of \eqref{su2}. On the other hand, we only rely on the control of $\norm{R}_{\ll{2}}$ and $\norm{\nabla k}_{\ll{2}}$ given by \eqref{l2bounds} to prove the estimate \eqref{su1}. 
\end{remark}

We postpone the proof of Lemma \ref{su3} and Lemma \ref{su4} respectively to section \ref{proofsu3} and section \ref{proofsu4}. Let us now conclude the proof of Theorem \ref{th2}. Note first that $W$ is not empty. In fact, the flat initial data set $(\s,g,k)=(\R^3,\de,0)$ belongs to $V$, where $\de$ denotes the euclidean metric. In that case, our construction in \cite{param1} yields the usual Fourier phase $u(x,\o)=\xo$. 
Then, the system \eqref{su2} reduces to:
\begin{equation}\label{su5}
\left\{\begin{array}{l}
\f(f_+)+\f(f_-)=\phi_0,\\
\f(\la f_+)-\f(\la f_-)=i\phi_1,
\end{array}\right.
\end{equation}
which admits the solution:
\begin{equation}\label{su6}
f_\pm=\frac{1}{2}\left(\mathcal{F}(\phi_0)\pm i\frac{\mathcal{F}(\phi_1)}{\la}\right),
\end{equation}
where $\mathcal{F}$ denotes the Fourier transform on $\R^3$. Thus, $(\s,g,k)=(\R^3,\de,0)$ belongs to $W$, which implies that $W$ is not empty. It is also open and closed in $V$ for the topology of $(g,k)\in C^q(\s)\times C^{q-1}(\s)$ by Lemma \ref{su4} for $q$ sufficiently large. Since $V$ is connected for the topology of $(g,k)\in C^q(\s)\times C^{q-1}(\s)$ by Lemma \ref{su3}, this implies that $W=V$. This proves the existence of $(f_+,f_-)$ solution of \eqref{su2} and concludes the proof of Theorem \ref{th2}. \QED 

\subsubsection{Proof of Lemma \ref{su3}}\label{proofsu3}

\vspace{0.2cm}

\noindent{\bf The conformal method of Lichnerowicz.} We start by reviewing the conformal method of Lichnerowicz for constructing solutions to the constraint equations \eqref{const1} on $\s$. Let $\underline{g}$ a Riemannian metric on $\s$. We define the Riemannian metric $g$ and the symmetric 2-tensor $k$ as:
\begin{equation}\label{su7}
\left\{\begin{array}{l}
\ds g=\phi^4\underline{g},\\
\ds k=\phi^{-2}\sigma,
\end{array}\right.
\end{equation}
where $\sigma$ is a traceless symmetric 2-tensor and $\phi$ a conformal factor tending to 1 at infinity. Then, $(g,k)$ defined in \eqref{su7} satisfies the constraint equations \eqref{const1} provided that $(\phi,\sigma)$ satisfy the following system:
\begin{equation}\label{su8}
\left\{\begin{array}{l}
\ds -8\Delta\phi+R\phi-|\sigma|^2\phi^{-7}=0,\\
\ds\textrm{div}\sigma=0,
\end{array}\right.
\end{equation}
where $R$ is the scalar curvature of $\underline{g}$ and where the divergence and the Laplacian are taken with respect to $\underline{g}$. 

\vspace{0.2cm}

\noindent{\bf The existence of $\sigma$.} We now turn to the question of the existence of $\phi$ and $\sigma$ solution to \eqref{su8}. In order to exploit the smallness condition \eqref{l2bounds}, we need an existence theory for rough solutions to the constraint equations \eqref{const1}. We will follow the exposition in \cite{Max} (we refer to \cite{ChBrYo} for the smooth case). Let $l\in\N$ and $\rho\in\R$. We introduce the spaces $H^l_\rho(\s)$ defined by:
\begin{equation}\label{defhkrho}
H^l_\rho(\s)=\left\{h\,/\,\sum_{|\a|\leq l}\norm{(1+|x|)^{-\rho-3/2+|\a|}h}_{\ll{2}}<+\infty\right\}.
\end{equation}
We recall first the construction of a symmetric traceless divergence free 2-tensor $\sigma$ on $\s$. To this end, we introduce the conformal Killing operator $\mathbb{L}$ and the vector Laplacian $\Delta_{\mathbb{L}}$:
\begin{equation}\label{su9}
\left\{\begin{array}{l}
\ds\mathbb{L}X=\mathcal{L}_X\underline{g}-\frac{2}{3}\textrm{div}(X)\underline{g},\\
\ds\Delta_{\mathbb{L}}=\textrm{div}(\mathbb{L}X),
\end{array}\right.
\end{equation}
where $X$ is a vectorfield on $\s$, $\mathcal{L}_X$ is the Lie derivative with respect to $X$, and the divergence is taken with respect to the metric $\underline{g}$. If $S$ is a symmetric traceless 2-tensor, and if we can solve 
\begin{equation}\label{su10}
\Delta_{\mathbb{L}}X=-\textrm{div}(S),
\end{equation} 
then setting $\sigma=S+\mathbb{L}X$ yields div$\sigma=0$ which solves the second equation of \eqref{su8}. The fact that this is always possible is known as the York decomposition. In the context of a rough metric $\underline{g}$, the following result holds (see \cite{Max}):
\begin{equation}\label{su10b}
\begin{array}{l}
\ds\textrm{Let }1-<\rho<0,\,\underline{g}\in H^2_\rho(\s)\textrm{ and }S\in H^1_{\rho-1}(\s). \textrm{ There is a unique }X\\
\ds\textrm{solution to \eqref{su10}, and }X\textrm{ satisfies }\norm{X}_{H^2_\rho(\s)}\lesssim \norm{S}_{H^1_{\rho-1}(\s)}.\\
\textrm{This yields a solution }\sigma\textrm{ to div}\sigma=0\textrm{ such that }\sigma=S+\mathbb{L}X\\
\ds\textrm{and }\norm{\sigma}_{H^1_{\rho-1}(\s)}\lesssim \norm{S}_{H^1_{\rho-1}(\s)}.
\end{array}
\end{equation} 

\vspace{0.2cm}

\noindent{\bf The existence of $\phi$.} We then have to solve the first equation of \eqref{su8} which is the Lichnerowicz equation. This is 
not an easy task in general since one has to show that $\underline{g}$ is conformally related to a metric with vanishing scalar curvature. However, we are in the particular case of small data in view of \eqref{l2bounds}, and we will obtain the existence of $\phi$ by a fixed point method. Let for $-1<\rho<0$ and let $g\in H^2_\rho(s)$. Then, recall from \cite{Max} that $-\Delta$ is invertible as an operator from $H^2_\rho(\s)$ to $H^0_{\rho-2}(\s)$ so that the following estimate holds:
\begin{equation}\label{su11}
\norm{(-\Delta)^{-1}h}_{H^2_\rho(\s)}\lesssim\norm{h}_{H^0_{\rho-2}(\s)},\,-1<\rho<0. 
\end{equation} 
This allows us to rewrite the first equation of \eqref{su8} in the form of a fixed point for $\psi=\phi-1$:
\begin{equation}\label{su11:1}
\psi=\frac{1}{8}(-\Delta)^{-1}\left(-R+|\sigma|^2-R\psi+|\sigma|^2((1+\psi)^{-7}-1)\right).
\end{equation} 
Now, we deduce from the embedding of $H^2_\rho(\s)$ in $L^\infty(\s)$ for $\rho<0$ and from the properties of the spaces $H^l_\rho(\s)$ with respect to the pointwise multiplication proved in \cite{Max} the following inequality:
\begin{equation}\label{su11:2}
\begin{array}{r}
\ds\norm{-R+|\sigma|^2-R\psi+|\sigma|^2((1+\psi)^{-7}-1)}_{H^0_{\rho-2}(\s)}\lesssim\ds\norm{R}_{H^0_{\rho-2}(\s)}(1+\norm{\psi}_{H^2_{\rho}(\s)})\\
\ds +\norm{\sigma}^2_{H^1_{\rho-1}(\s)}(1+\vartheta(\norm{\psi}_{H^2_{\rho}(\s)})),
\end{array}
\end{equation} 
where $\vartheta$ is an increasing function, and where we assume that the control $\norm{\psi}_{H^2_\rho(\s)}\leq 1/2$ holds. Thus, in view of \eqref{su10b}, we have for $-1<\rho<0$ and $\norm{\psi}_{H^2_\rho(\s)}\leq 1/2$:
\begin{equation}\label{su11:3}
\begin{array}{r}
\ds\norm{-R+|\sigma|^2-R\psi+|\sigma|^2((1+\psi)^{-7}-1)}_{H^0_{\rho-2}(\s)}\lesssim\ds\norm{R}_{H^0_{\rho-2}(\s)}(1+\norm{\psi}_{H^2_{\rho}(\s)})\\
\ds +\norm{S}^2_{H^1_{\rho-1}(\s)}(1+\vartheta(\norm{\psi}_{H^2_{\rho}(\s)})),
\end{array}
\end{equation} 
where $\vartheta$ is an increasing function. In view of \eqref{su11:3}, we immediately obtain the existence of $\psi$ solution to \eqref{su11:1} provided $\norm{R}_{H^0_\rho(\s)}+\norm{S}_{H^1_{\rho-1}(\s)}\lesssim\ep$ for a sufficiently small $\ep$.

\vspace{0.2cm}

\noindent{\bf Proof of Lemma \ref{su3}.} Let us come back to the proof of Lemma \ref{su3}. We will prove that all solutions $(\s,g,k)$ of the constraint equations \eqref{const1} satisfying the bound \eqref{l2bounds} are connected to $(\R^3,\de,0)$ by a continuous path. For $0\leq\tau\leq 1$, we introduce:
\begin{equation}\label{su12}
\underline{g}_\tau=\tau g+(1-\tau)\de\textrm{ and }S_\tau=\tau k -\frac{\tau\textrm{Tr}_\tau k}{3}\underline{g}_\tau, 
\end{equation} 
where Tr$_\tau$ denotes the trace with respect to the metric $\underline{g}_\tau$. Let $-1<\rho<0$. From the smallness assumptions \eqref{l2bounds} and the definition \eqref{su12} of $\underline{g}_\tau$ and $S_\tau$, we immediately obtain:
\begin{equation}\label{su12:1}
\norm{\underline{R}_\tau}_{H^0_{\rho-2}(\s)}+\norm{S_\tau}_{H^1_{\rho-1}(\s)}\lesssim\ep,
\end{equation} 
where $\underline{R}_\tau$ is the scalar curvature of $\underline{g}_\tau$. In view of \eqref{su10b}, \eqref{su11:1} \eqref{su11:3} and \eqref{su12:1}, we obtain the existence of $(\sigma_\tau, \phi_\tau)$ 
in $H^1_{\rho-1}(\s)\times H^2_\rho(s)$ solution to:
\begin{equation}\label{su12:2}
\left\{\begin{array}{l}
\ds -8\Delta\phi_\tau+\underline{R}_\tau\phi_\tau-|\sigma_\tau|^2\phi_\tau^{-7}=0,\\
\ds\textrm{div}\sigma_\tau=0,
\end{array}\right.
\end{equation}
where $\underline{R}_\tau$ is the scalar curvature of $\underline{g}_\tau$ and where the divergence and the Laplacian are taken with respect to $\underline{g}_\tau$. Finally, setting 
\begin{equation}\label{su13}
\left\{\begin{array}{l}
\ds g_\tau=\phi_\tau^4\underline{g}_\tau,\\
\ds k_\tau=\phi_\tau^{-2}\sigma_\tau,
\end{array}\right.
\end{equation}
we obtain a solution $(\s,g_\tau,k_\tau)$ to the constraint equations \eqref{const1} which satisfies the following bound:
\begin{equation}\label{su13:1}
\norm{g_\tau-\de}_{H^2_{\rho}(\s)}+\norm{k_\tau}_{H^1_{\rho-1}(\s)}\lesssim\ep.
\end{equation} 
Thus $(g_\tau,k_\tau)$ satisfies the bound \eqref{l2bounds} so that $(\s,g_\tau, k_\tau)$ belongs to the set $V$ defined by \eqref{defv}. Furthermore, recall from \eqref{su11:1} that $\phi_\tau$ is obtained by a fixed point argument. This implies in particular the uniqueness of $(\sigma_\tau,\phi_\tau)$ so that $(g_\tau,k_\tau)=(\de,0)$ at $\tau=0$ and $(g_\tau,k_\tau)=(g,k)$ at $\tau=1$. Using standard results in elliptic regularity, we also obtain that the path $\tau\rightarrow (g_\tau,k_\tau)$ is continuous for the topology of $C^q(\s)\times C^{q-1}(\s)$ provided $(g,k)\in C^q(\s)\times C^{q-1}(\s)$. Thus, all solutions $(\s,g,k)$ of the constraint equations \eqref{const1} satisfying the bound \eqref{l2bounds} are connected to $(\R^3,\de,0)$ by a continuous path, which concludes the proof of Lemma \ref{su3}. \QED

\begin{remark}
In general, the connectedness of the set of all solutions $(\s,g,k)$ of the constraint equations \eqref{const1} is an open problem (see \cite{SW} for a partial answer). Here, the smallness condition \eqref{l2bounds} makes the problem much easier, as the solutions are obtained by a fixed point argument in this case.
\end{remark}

\subsubsection{Proof of Lemma \ref{su4}}\label{proofsu4}

\vspace{0.2cm}

\noindent{\bf The operator $\Lambda$.} We start by rewriting the system \eqref{su2} as:
\begin{equation}\label{su14}
\left\{\begin{array}{l}
\nabla M_+f_++\nabla M_-f_-=\nabla\phi_0,\\
Q_+(\la f_+)-Q_-(\la f_-)=i\phi_1.
\end{array}\right.
\end{equation}
We define the operators $\widetilde{M}_\pm$ as:
\begin{equation}\label{su15}
\widetilde{M}_\pm f(x)=\pm \int_{\S}\int_{0}^{+\infty}e^{\pm i\lambda u(x,\pm\o)}a(x,\pm\o)^{-1}N(x,\pm\o) f(\lambda\o)\lambda^2 d\lambda d\o,
\end{equation}
so that \eqref{su14} becomes:
\begin{equation}\label{su16}
\left\{\begin{array}{l}
\widetilde{M}_+(\la f_+)-\widetilde{M}_-(\la f_-)=-i\nabla\phi_0,\\
Q_+(\la f_+)-Q_-(\la f_-)=i\phi_1.
\end{array}\right.
\end{equation}
We define the linear operator $\Lambda$ as:
\begin{equation}\label{su17}
\Lambda(f_+,f_-)=(\widetilde{M}_+(f_+)-\widetilde{M}_-(f_-),Q_+(f_+)-Q_-(f_-)).
\end{equation}
By Theorem \ref{th1} and {\bf Assumption 1-4} on $u$, $a$ and $N$, $\Lambda$ is a bounded operator from $\le{2}\times\le{2}$ to $\ll{2}^3\times\ll{2}$. In view of \eqref{choicef4bis}, it satisfies the following estimate:
\begin{equation}\label{su18}
\norm{f_+}_{\le{2}}+\norm{f_-}_{\le{2}}\lesssim\norm{\Lambda(f_+,f_-)}_{\ll{2}^3\times\ll{2}}.
\end{equation}
Finally, in view of \eqref{su16} and the definition \eqref{su17} of $\Lambda$, we may rewrite the set $W$ as:
\be{su19}
W=\{(\s,g,k)\in V\textrm{ such that }\Lambda\textrm{ is surjective}\}.
\eeq

\vspace{0.2cm}

\noindent{\bf $W$ is closed.} We have to show that the set $W$ given by \eqref{su19} is both open and closed for the $C^q(\s)\times C^{q-1}(\s)$ topology when $q$ is sufficiently large. Let us first show that $W$ is closed. Let $(g_n,k_n), n\in\N$ a sequence in $W$ such that it has a limit $(g,k)$ for the $C^q(\s)\times C^{q-1}(\s)$ topology. Let $\Lambda_n$ be the operator associated to $(g_n,k_n)$, and $\Lambda$ the operator associated to $\Lambda$. $\Lambda_n$ is surjective for all $n\geq 0$, and we would like to prove that $\Lambda$ is surjective. Notice first that the fact that $\Lambda$ is a bounded operator from $\le{2}\times\le{2}$ to $\ll{2}^3\times\ll{2}$ together with the estimate \eqref{su18} implies that the image of $\Lambda$ is closed in $\ll{2}^3\times\ll{2}$. Thus, we may reduce the problem to showing that a dense subset of $\ll{2}^3\times\ll{2}$ belongs to the image of $\Lambda$. Let us consider $(\phi_0,\phi_1)$ in $C^7_c(\s)\times C^6_c(\s)$ which is dense in $\ll{2}^3\times\ll{2}$. Since $\Lambda_n$ is surjective, there are $(f_+^n,f_-^n)$ such that:
\begin{equation}\label{su20}
\Lambda_n(f^n_+,f^n_-)=(\nabla\phi_0,\phi_1).
\end{equation}
Differentiating \eqref{su20} six times and using \eqref{su18}, we obtain:
\begin{equation}\label{su21}
\begin{array}{ll}
&\ds\norm{(1+\la^6)f^n_+}_{\le{2}}+\norm{(1+\la^6)f^n_-}_{\le{2}}\\
\ds\lesssim &\ds(\norm{g_n}_{C^q(\s)}+\norm{k_n}_{C^{q-1}(\s)})(\norm{\phi_0}_{C^7_c(\s)}+\norm{\phi_1}_{C^6_c}),
\end{array}
\end{equation}
for a sufficiently large $q$. We deduce from \eqref{su21} the existence of a constant $C>0$ independent of $n$ such that:
\begin{equation}\label{su22}
\ds\norm{(1+\la^6)f^n_+}_{\le{2}}+\norm{(1+\la^6)f^n_-}_{\le{2}}\leq C<+\infty.
\end{equation}
In particular, we may assume in up to a subsequence that $(f^n_+,f^n_-)$ converges up to a subsequence to $(f_+,f_-)$ weakly in $\le{2}\times\le{2}$. We have:
\begin{equation}\label{su23}
\ds\Lambda_n(f^n_+,f^n_-)-\Lambda(f_+,f_-)=(\Lambda_n-\Lambda)(f^n_+,f^n_-)+\Lambda(f^n_+-f_+,f^n_--f_-).
\end{equation}
We will show that both terms in the right-hand side of \eqref{su23} converge to 0 weakly in $\ll{2}^3\times\ll{2}$. We start with the first term. For $(H,h)\in C^0_c(\s)^3\times C^0_c(\s)$, we have in view of the definition \eqref{su17} of $\Lambda$:
\bea\label{su23:2}
&&\ds\left|\int_\s ((\Lambda_n-\Lambda)(f^n_+,f^n_-),(H,h))d\s\right|\\[3mm]
\nn&= &\ds\left|\int_{\S}\int_{0}^{+\infty}((f^n_+,f^n_-),(\Lambda_n-\Lambda)^*(H,h))\la^2d\la d\o\right|\\
\nn&\lesssim &\!\!\!\!\ds\norm{(H,h)}_{C^0_c(\s)^3\times C^0_c(\s)}\int_{\S}\int_{0}^{+\infty}(|f^n_+(\la\o)|+|f^n_-(\la\o)|)(\la\norm{u_n(.,\o)-u(.,\o)}_{\ll{\infty}}\\
\nn&&\ds +\norm{a^{-1}_n(.,\o)-a^{-1}(.,\o)}_{\ll{\infty}}+\norm{N_n(.,\o)-N(.,\o)}_{\ll{\infty}})\la^2d\la d\o\\
\nn&\lesssim &\ds\norm{(H,h)}_{C^0_c(\s)^3\times C^0_c(\s)}(\norm{(1+\la^6)f^n_+}_{\le{2}}+\norm{(1+\la^6)f^n_-}_{\le{2}})\\
\nn&&\ds \times\sup_{\o\in\S}(\norm{u_n(.,\o)-u(.,\o)}_{\ll{\infty}}+\norm{a^{-1}_n(.,\o)-a^{-1}(.,\o)}_{\ll{\infty}}\\
\nn&&\ds +\norm{N_n(.,\o)-N(.,\o)}_{\ll{\infty}}).
\eea
Since $(g_n,k_n)$ converges to $(g,k)$ in $C^{q}(\s)\times C^{q-1}(\s)$, we have for $q$ large enough:
\bea\label{su23:3}
\lim_{n\rightarrow +\infty}\sup_{\o\in\S}(\norm{u_n(.,\o)-u(.,\o)}_{\ll{\infty}}+\norm{a^{-1}_n(.,\o)-a^{-1}(.,\o)}_{\ll{\infty}}&&\\
\nn\ds +\norm{N_n(.,\o)-N(.,\o)}_{\ll{\infty}})& = &0.
\eea
Using \eqref{su22}, \eqref{su23:2} and \eqref{su23:3} implies that $(\Lambda_n-\Lambda)(f^n_+,f^n_-)$  converges weakly in $\ll{2}^3\times\ll{2}$ to 0. Also, using the fact that $\Lambda$ is a bounded operator from $\le{2}\times\le{2}$ to $\ll{2}^3\times\ll{2}$ and that $(f^n_+,f^n_-)$ converges to $(f_+,f_-)$ weakly in $\le{2}\times\le{2}$, we obtain that $\Lambda(f^n_+-f_+,f^n_--f_-)$ converges weakly in $\ll{2}^3\times\ll{2}$ to 0. In view of \eqref{su23}, this implies that $\Lambda_n(f^n_+,f^n_-)$ converges weakly to $\Lambda(f_+,f_-)$ in $\ll{2}^3\times\ll{2}$. Together with \eqref{su20}, this implies 
\begin{equation}\label{su24}
\Lambda(f_+,f_-)=(\nabla\phi_0,\phi_1).
\end{equation}
Thus, $\Lambda$ is surjective which concludes the proof of the fact that $W$ is closed.

\vspace{0.2cm}

\noindent{\bf $W$ is open.} To conclude the proof of Lemma \ref{su4}, we need to prove that $W$ is open for the $C^q(\s)\times C^{q-1}(\s)$ topology when $q$ is sufficiently large. Let $(g,k)\in W$ and let $\Lambda$ the operator associated to $(g,k)$. Then $\Lambda$ is surjective which together with the estimate \eqref{su18} implies that $\Lambda$ is an isomorphism from $\le{2}\times\le{2}$ to $\ll{2}^3\times\ll{2}$. In turn, this implies that $\Lambda\Lambda^*$ is an isomorphism of $\ll{2}^3\times\ll{2}$. Let $(\widetilde{g},\widetilde{k})\in W$ such that:
\begin{equation}\label{su25}
\norm{\widetilde{g}-g}_{C^q(\s)}+\norm{\widetilde{k}-k}_{C^{q-1}(\s)}\leq\de
\end{equation}
for a small constant $\de>0$ to be chosen later, and let $\widetilde{\Lambda}$ the operator associated to $(\widetilde{g},\widetilde{k})$. Then, $\widetilde{\Lambda}$ and $\Lambda$ consist of Fourier integral operators whose phase and symbol are $O(\de)$ close to each other in the $C^q(\s)$ topology. Integrating by parts several times in the kernel of $\Lambda\Lambda^*$ and $\widetilde{\Lambda}\widetilde{\Lambda}^*$, we deduce the following bound provided $q$ is sufficiently large:
\begin{equation}\label{su26}
\norm{\widetilde{\Lambda}\widetilde{\Lambda}^*-\Lambda\Lambda^*}_{\mathcal{L}(\ll{2}^3\times\ll{2})}\lesssim\de.
\end{equation}
Since the isomorphism of $\ll{2}^3\times\ll{2}$ form an open set, we deduce from \eqref{su26} that 
$\widetilde{\Lambda}\widetilde{\Lambda}^*$ is an isomorphism of $\ll{2}^3\times\ll{2}$ for $\de>0$ small enough. In particular, $\widetilde{\Lambda}$ is surjective for $\de>0$ small enough. Therefore, 
$\widetilde{\Lambda}\in W$ provided $\de>0$ defined in \eqref{su25} is chosen small enough. Thus, 
$W$ is open. This concludes the proof of Lemma \ref{su4}. \QED

\appendix

\section{Proof of Lemma \ref{doublediv}}\label{sec:doublediv}

We would like to compute the double divergence term in the right-hand side of \eqref{bisoa21}:
\begin{equation}\label{ap1}
\textrm{div}\left(\frac{(N-(N\cdot N')N')a}{1-(N\cdot N')^2}\textrm{div}\left(\frac{(N'-(N\cdot N')N)a'bb'}{1-(N\cdot N')^2}\right)\right).
\end{equation}
We recall the structure equations for $N$:
\begin{equation}\label{frame}
\left\{\begin{array}{l}
\nabla_AN=\th_{AB}e_B,\\[1mm]
\nabn N=-\nabb\lg.
\end{array}\right.
\end{equation}
In particular, \eqref{frame} implies:
\begin{equation}\label{ap2}
\textrm{div}(N)=\trt.
\end{equation}
Using \eqref{ap2}, we have:
\begin{equation}\label{ap3}
\begin{array}{ll}
&\ds\textrm{div}\left(\frac{(N'-(N\cdot N')N)a'bb'}{1-(N\cdot N')^2}\right)\\
\ds = &\ds\frac{(\trt'-\gn\trt-\nabn(\gn))a'bb'+\nabla_{N'-\gn N}(a'bb')}{1-g(N,N')^2}\\
&\ds +\frac{2a'bb'\nabla_{N'-\gn N}(\gn)\gn}{(1-g(N,N')^2)^2}.
\end{array}
\end{equation}
Differentiating again, we obtain:
\begin{equation}\label{ap4}
\begin{array}{ll}
&\ds\textrm{div}\left(\frac{(N-(N\cdot N')N')a}{1-(N\cdot N')^2}\textrm{div}\left(\frac{(N'-(N\cdot N')N)a'bb'}{1-(N\cdot N')^2}\right)\right)\\
\ds = &\ds\frac{A_1}{(1-g(N,N')^2)^2}+\frac{A_2}{(1-g(N,N')^2)^3}+\frac{A_3}{(1-g(N,N')^2)^4}
\end{array}
\end{equation}
where $A_1, A_2, A_3$ are given by:
\begin{equation}\label{ap5}
\begin{array}{ll}
\ds A_1= &\ds (\nabla_{N-\gn N'}\trt'-\gn \nabla_{N-\gn N'}\trt)aa'bb'\\
&\ds -(\nabla_{N-\gn N'}(\gn)\trt+\nabla_{N-\gn N'}\nabn(\gn))aa'bb'\\
&\ds +(\trt'-\gn\trt-\nabn(\gn))a\nabla_{N-\gn N'}(a'bb')\\
&\ds +a\nabla_{N-\gn N'}\nabla_{N'-\gn N}(a'bb')\\
&\ds +((\trt-\gn\trt'-\nabla_{N'}(\gn))a+\nabla_{N-\gn N'}(a))\\
&\ds \times((\trt'-\gn\trt-\nabn(\gn))a'bb'+\nabla_{N'-\gn N}(a'bb')),
\end{array}
\end{equation}
\begin{equation}\label{ap6}
\begin{array}{ll}
\ds A_2= &\ds 2a\nabla_{N-\gn N'}(\gn)\gn\\
&\ds\times((\trt'-\gn\trt-\nabn(\gn))a'bb'+\nabla_{N'-\gn N}(a'bb'))\\
&\ds +2aa'bb'\nabla_{N-\gn N'}\nabla_{N'-\gn N}(\gn)\gn\\
&\ds +2aa'bb'\nabla_{N'-\gn N}(\gn)\nabla_{N-\gn N'}(\gn)\\
&\ds +2a\nabla_{N'-\gn N}(\gn)\gn\nabla_{N-\gn N'}(a'bb')\\
&\ds +((\trt-\gn\trt'-\nabla_{N'}(\gn))a+\nabla_{N-\gn N'}(a))\\
&\ds\times 2a'bb'\nabla_{N'-\gn N}(\gn)\gn\\
& +2a\nabla_{N-\gn N'}(\gn)\gn\\
&\ds \times((\trt'-\gn\trt-\nabn(\gn))a'bb'+\nabla_{N'-\gn N}(a'bb')),
\end{array}
\end{equation}
and 
\begin{equation}\label{ap7}
\ds A_3=8aa'bb'\nabla_{N'-\gn N}(\gn)\nabla_{N-\gn N'}(\gn)\gn^2.
\end{equation}

Notice that $N-(N\c N')N'$ is tangent to $P_{u'}$ and that $N'-(N\c N')N$ is tangent to $\p$. Notice also that $1-\gn$ is of order two in $N-N'$:
\begin{equation}\label{ap8}
1-\gn = \frac{g(N-N',N-N')}{2}.
\end{equation} 
In view of \eqref{ap4}-\eqref{ap8}, one easily checks that the double divergence \eqref{ap1} takes the wanted form \eqref{dd1} \eqref{dd2} provided that we are able to control all the terms in the following list:
\begin{equation}\label{ap9}
\begin{array}{l}
\ds \frac{\nabn(\gn)}{(1-\gn)^{1/2}},\, \frac{\nabla_{N-\gn N'}(N'-\gn N)}{1-\gn},\,  \frac{\nabla_{N-\gn N'}(\gn)}{(1-\gn)^{3/2}},\\
\ds\frac{\nabla_{N-\gn N'}\nabn(\gn)}{1-\gn},\, \frac{\nabla_{N-\gn N'}\nabla_{N'-\gn N}(\gn)}{(1-\gn)^{2}}.
\end{array}
\end{equation}

\vspace{0.2cm}

\noindent{\bf Control of the first term of \eqref{ap9}.} Using the structure equation for $N$ \eqref{frame}, we have:
\bea\label{ap10}
&&\nabn(\gn)\\
\nn&=&g(\nabn N,N')+\gn g(N,\nabla_{N'}N')+g(N,\nabla_{N-\gn N'}N')\\
\nn& =&-g(\nabb\log(a),N')+\gn g(N,\nabb'\log(a'))+\th'(N-\gn N',N-\gn N')\\
\nn& =&-g(\nabla\log(a),N'-\gn N)+\gn g(N-\gn N',\nabla\log(a'))\\
\nn&& +\th'(N-\gn N',N-\gn N').
\eea
Using \eqref{ap8}, we have:
\begin{equation}\label{ap11}
\frac{|N-\gn N'|}{(1-\gn)^{1/2}}+\frac{|N'-\gn N|}{(1-\gn)^{1/2}} \lesssim 1.
\end{equation} 
In view of \eqref{ap8}, \eqref{ap10} and \eqref{ap11}, the term $\frac{\nabn(\gn)}{(1-\gn)^{1/2}}$ is under control and involves terms in the list \eqref{dd2}.

\vspace{0.2cm}

\noindent{\bf Control of the second term of \eqref{ap9}.} Using the structure equation for $N$ \eqref{frame}, we have:
\begin{equation}\label{ap12}
\begin{array}{l}
\ds \nabla_{N-\gn N'}(N'-\gn N)=\nabla_{N-\gn N'}N'-\gn \nabla_{N-\gn N'}N\\
\ds -\nabla_{N-\gn N'}(\gn)N\\
\ds = \th'(N-\gn N',e_{A'})e_{A'}\\
\ds -\gn((1-\gn^2)\nabn N-\gn\nabla_{N'-\gn N}N)\\
\ds -((1-\gn^2)g(-\nabb\log(a),N')\\
\ds -\gn\th(N'-\gn N,N'-\gn N)\\
\ds +\th'(N-\gn N',N-\gn N'))N\\
\ds = \th'(N-\gn N',e_{A'})e_{A'}+\th(N'-\gn N, e_A)e_A \\
\ds -(1-\gn^2)\th(N'-\gn N,e_A)e_A
+\gn (1-\gn^2)\nabb\log(a)\\
\ds +(1-\gn^2)\nabla_{N'-\gn N}\log(a) N\\
\ds +\gn\th(N'-\gn N,N'-\gn N)N\\
\ds -\th'(N-\gn N',N-\gn N')N,
\end{array}
\end{equation}
where we have used the fact that:
\begin{equation}\label{ap13}
\ds N-\gn N'=(1-\gn^2)N-\gn(N'-\gn N).
\end{equation}
Note that the tangential components $N-(N\c N')N'$ and $N'-(N\c N')N$ satisfy:
\begin{equation}\label{ap14}
(N-\gn N')+(N'-\gn N)=(N+N')(1-\gn),
\end{equation}
so that we may divide $\th'(N-\gn N',e_{A'})e_{A'}+\th(N'-\gn N, e_A)e_A$ by $1-\gn$. Thus, 
in view of \eqref{ap8}, \eqref{ap11}, \eqref{ap12} and \eqref{ap14}, the term $\frac{\nabla_{N-\gn N'}(N'-\gn N)}{1-\gn}$ is under control and involves terms in the list \eqref{dd2}.

\vspace{0.2cm}

\noindent{\bf Control of the third term of \eqref{ap9}.} Using the structure equation for $N$ \eqref{frame} together with \eqref{ap13}, we have:
\begin{equation}\label{ap17}
\begin{array}{l}
\ds \nabla_{N-\gn N'}(\gn)=g(\nabla_{N-\gn N'}N,N')+g(N,\nabla_{N-\gn N'}N')\\
\ds =-(1-\gn^2)g(\nabla\log(a),N'-\gn N)\\
\ds -\gn\th(N'-\gn N,N'-\gn N)\\
\ds +\th'(N-\gn N',N-\gn N').
\end{array}
\end{equation}
In view of \eqref{ap8}, \eqref{ap11}, \eqref{ap14} and \eqref{ap17}, the term $\frac{\nabla_{N-\gn N'}(\gn)}{(1-\gn)^{3/2}}$ is under control and involves terms in the list \eqref{dd2}.

\vspace{0.2cm}

\noindent{\bf Control of the fourth term of \eqref{ap9}.} Differentiating \eqref{ap10} with respect to $\nabla_{N-\gn N'}$, we have:
\begin{equation}\label{ap15}
\begin{array}{l}
\ds \nabla_{N-\gn N'}\nabn(\gn)=-\nabla^2\log(a)(N-\gn N',N'-\gn N)\\
\ds -g(\nabla\log(a),\nabla_{N-\gn N'}(N'-\gn N))\\
\ds +\nabla_{N-\gn N'}(\gn)g(\nabla\log(a'),N-\gn N')\\
\ds +\gn \nabla^2\log(a')(N-\gn N',N-\gn N')\\
\ds +\gn g(\nabla\log(a'),\nabla_{N-\gn N'}(N-\gn N'))\\
\ds +\nabla_{N-\gn N'}\th'(N-\gn N',N-\gn N')\\
\ds +2\th'(\nabla_{N-\gn N'}(N-\gn N'),N-\gn N').
\end{array}
\end{equation}
In view of \eqref{ap15}, we need to control the term $\frac{\nabla_{N-\gn N'}(N-\gn N')}{1-\gn}$. This is very similar to \eqref{ap12}. Using the structure equation for $N$ \eqref{frame}, we obtain:
\begin{equation}\label{ap16}
\begin{array}{l}
\ds \nabla_{N-\gn N'}(N-\gn N')=\nabla_{N-\gn N'}N-\gn \nabla_{N-\gn N'}N'\\
\ds -\nabla_{N-\gn N'}(\gn)N'\\
\ds = (1-\gn^2)\nabn N-\gn\nabla_{N'-\gn N}N\\
\ds -\gn\th'(N-\gn N',e_{A'})e_{A'}\\
\ds -((1-\gn^2)g(-\nabb\log(a),N')\\
\ds -\gn\th(N'-\gn N,N'-\gn N)\\
\ds +\th'(N-\gn N',N-\gn N'))N'\\
\ds = -(1-\gn^2)\nabb\log(a)-\gn\th(N'-\gn N, e_A)e_A\\
\ds -\gn\th'(N-\gn N',e_{A'})e_{A'} \\
\ds +(1-\gn^2)\nabla_{N'-\gn N}\log(a) N'\\
\ds +\gn\th(N'-\gn N,N'-\gn N)N'\\
\ds -\th'(N-\gn N',N-\gn N')N',
\end{array}
\end{equation}
In view of \eqref{ap8}, \eqref{ap11}, \eqref{ap14} and \eqref{ap16}, the term $\frac{\nabla_{N-\gn N'}(N-\gn N')}{1-\gn}$ is under control and involves terms in the list \eqref{dd2}. Note also that the terms $\nabla^2\log(a)(N-\gn N',N'-\gn N)$ and $\nabla^2\log(a')(N-\gn N',N-\gn N')$ appearing in \eqref{ap15} 
both contain at least one tangential derivative. Together with \eqref{ap8}, \eqref{ap11}, \eqref{ap14}, \eqref{ap17}, \eqref{ap15} and \eqref{ap16}, this yields that the term $\frac{\nabla_{N-\gn N'}\nabn(\gn)}{1-\gn}$ is under control and involves terms in the list \eqref{dd2}.

\vspace{0.2cm}

\noindent{\bf Control of the fifth term of \eqref{ap9}.} Exchanging the role of $N$ and $N'$ in \eqref{ap17}, we obtain:
\begin{equation}\label{ap18}
\begin{array}{l}
\ds \nabla_{N'-\gn N}(\gn) =-(1-\gn^2)g(\nabla\log(a'),N-\gn N')\\
\ds -\gn\th'(N-\gn N',N-\gn N')\\
\ds +\th(N'-\gn N,N'-\gn N).
\end{array}
\end{equation}
Differentiating \eqref{ap18} with respect to $\nabla_{N-\gn N'}$, we obtain:
\begin{equation}\label{ap19}
\begin{array}{l}
\ds \nabla_{N-\gn N'}\nabla_{N'-\gn N}(\gn) \\
\ds =-(1-\gn^2)\nabla^2\log(a')(N-\gn N',N-\gn N')\\
\ds  -(1-\gn^2)\nabla_{\nabla_{N-\gn N'}(N-\gn N')}\log(a')\\
\ds +2\gn\nabla_{N-\gn N'}(\gn)\nabla_{N-\gn N'}\log(a')\\
\ds -\nabla_{N-\gn N'}(\gn)\th'(N-\gn N',N-\gn N')\\
\ds -\gn\nabla_{N-\gn N'}\th'(N-\gn N',N-\gn N')\\
\ds -2\gn\th'(\nabla_{N-\gn N'}(N-\gn N'),N-\gn N')\\
\ds +\nabla_{N-\gn N'}\th(N'-\gn N,N'-\gn N)\\
\ds +2\th(\nabla_{N-\gn N'}(N'-\gn N), N'-\gn N).
\end{array}
\end{equation}
Together with \eqref{ap12}, \eqref{ap17} and \eqref{ap16}, we get:
\begin{equation}\label{ap20}
\begin{array}{l}
\ds \nabla_{N-\gn N'}\nabla_{N'-\gn N}(\gn) \\
\ds =-(1-\gn^2)\nabla^2\log(a')(N-\gn N',N-\gn N')\\
\ds  +(1-\gn^2)^2 g(\nabb\log(a),\nabla\log(a'))\\
\ds +(1-\gn^2)\gn\th(N'-\gn N,\nabla\log(a'))\\
\ds +(1-\gn^2)\gn\th'(N-\gn N',\nabla\log(a'))\\
\ds -(1-\gn^2)^2\nabla_{N'-\gn N}\log(a)\nabla_{N'}\log(a')\\
\ds -(1-\gn^2)\gn\th(N'-\gn N,N'-\gn N)\nabla_{N'}\log(a')\\
\ds +(1-\gn^2)\th'(N-\gn N',N-\gn N')\nabla_{N'}\log(a')\\
\ds -2\gn\nabla_{N-\gn N'}\log(a')(1-\gn^2)\nabla_{N'-\gn N}\log(a)\\
\ds -2\gn^2\nabla_{N-\gn N'}\log(a')\th(N'-\gn N,N'-\gn N)\\
\ds +2\gn\nabla_{N-\gn N'}\log(a')\th'(N-\gn N',N-\gn N')\\
\ds +(1-\gn^2)\nabla_{N'-\gn N}\log(a)\th'(N-\gn N',N-\gn N')\\
\ds +\gn\th(N'-\gn N,N'-\gn N)\\
\times \th'(N-\gn N',N-\gn N')\\
\ds -\th'(N-\gn N',N-\gn N')^2\\
\ds -\gn\nabla_{N-\gn N'}\th'(N-\gn N',N-\gn N')\\
\ds +\nabla_{N-\gn N'}\th(N'-\gn N,N'-\gn N)\\
\ds +2\gn(1-\gn^2)\th'(\nabb\log(a),N-\gn N')\\
\ds +2\gn(1-\gn^2)\th(\nabb\log(a),N'-\gn N)\\
\ds -2(1-\gn^2)\th(e_A,N'-\gn N)^2\\
\ds +2\gn^2\th'(e_{A'},N-\gn N')\th(N'-\gn N,e_{A'})\\
\ds +2\gn^2\th'(e_{A'},N-\gn N')^2\\
\ds +2\gn^2\th'(e_{A},N-\gn N')\th(N'-\gn N,e_{A})\\
\ds +2\gn^2\th(e_{A},N'-\gn N)^2
\end{array}
\end{equation}
Note that the term $\nabla^2\log(a')(N-\gn N',N-\gn N')$ appearing in \eqref{ap20} contains at least one tangential derivative (it actually contains two tangential derivatives). Note also that the terms:
\begin{equation}\label{ap21}
\begin{array}{l}
\ds 2\gn^2\th'(e_{A'},N-\gn N')\th(N'-\gn N,e_{A'})\\
\ds +2\gn^2\th'(e_{A'},N-\gn N')^2\\
\ds +2\gn^2\th'(e_{A},N-\gn N')\th(N'-\gn N,e_{A})\\
\ds +2\gn^2\th(e_{A},N'-\gn N)^2,
\end{array}
\end{equation}
appearing in \eqref{ap20} may be rewritten:
\begin{equation}\label{ap22}
\begin{array}{l}
\ds 2\gn^2(\th'(N-\gn N',.)+\th(N'-\gn N,.))^2.
\end{array}
\end{equation}
Together with \eqref{ap8}, \eqref{ap11}, \eqref{ap14} and \eqref{ap20}, this yields that the term 
$$\frac{\nabla_{N-\gn N'}\nabla_{N'-\gn N}(\gn)}{(1-\gn)^2}$$ 
is under control and involves terms in the list \eqref{dd2}. This concludes the proof of Lemma \ref{doublediv}. \QED


\section{Proof of Lemma \ref{dd3}}

We start with the terms $\nabb\nabla(ab)a'b',\,\th\nabla(ab)a'b',\,\nabla(a)\nabla(b)a'b',\,\th^2aa'bb'$ in the list \eqref{dd2}. They all take the form \eqref{dd4} with $H_3=H_4=0$, $H_2=a'b'$ and taking respectively $H_1=\nabb\nabla(ab)$, $H_1=\th\nabla(ab)$, $H_1=\nabla(a)\nabla(b)$ and $H_1=\th^2ab$. Thanks to {\bf Assumption 1} and {\bf Assumption 2} on $a, \th$ and the assumptions \eqref{thregx1s} \eqref{threomega1s} on $b$, the estimates \eqref{dd5} and \eqref{dd6} are satisfied. 

We now consider the other terms: 
$$\frac{(\nabla\th-\nabla\th')aa'bb'}{|N_\nu-N_{\nu'}|}, \frac{(\th-\th')\nabla(ab)a'b'}{|N_\nu-N_{\nu'}|}, ab\th\nabla(a'b'), \nabla(ab)\nabla(a'b'), \frac{(\th-\th')^2aa'bb'}{|N_\nu-N_{\nu'}|^2}, \th\th'aa'bb'.$$ 
We focus on $\frac{(\nabla\th-\nabla\th')aa'bb'}{|N_\nu-N_{\nu'}|}$ and $\frac{(\th-\th')^2aa'bb'}{|N_\nu-N_{\nu'}|^2}$ the others being similar. For $\frac{(\nabla\th-\nabla\th')aa'bb'}{|N_\nu-N_{\nu'}|}$, we have:
\begin{equation}\label{dd7}
\frac{(\nabla\th-\nabla\th')aa'bb'}{|N_\nu-N_{\nu'}|}=\frac{(\nabla\th-\nabla\th_\nu)aa'bb'}{|N_\nu-N_{\nu'}|}+\frac{(\nabla\th_\nu-\nabla\th')aa'bb'}{|N_\nu-N_{\nu'}|}
\end{equation}
and the two terms in \eqref{dd7} are of the form \eqref{dd4} with $H_3=H_4=0$, and respectively $H_1=\frac{(\nabla\th-\nabla\th_\nu)ab}{|N_\nu-N_{\nu'}|}$, $H_2=a'b'$ and $H_1=\frac{(\nabla\th_{\nu}-\nabla\th')a'b'}{|N_\nu-N_{\nu'}|}$, $H_2=ab$. Thanks to {\bf Assumption 1} and {\bf Assumption 2} on $a, \th$ and the assumptions \eqref{thregx1s} \eqref{threomega1s} on $b$, the estimates \eqref{dd5} and \eqref{dd6} are satisfied. In particular, we have:
\begin{equation}\label{dd8}
\normm{\frac{(\nabla\th-\nabla\th_\nu)ab}{|N_\nu-N_{\nu'}|}}_{\ll{2}}\lesssim D\frac{\norm{\nabla\po\th}_{\ll{2}}}{2^{j/2}|\nu-\nu'|}\lesssim D,
\end{equation}
and 
\begin{equation}\label{dd9}
\normm{\frac{(\nabla\th_\nu-\nabla\th')ab}{|N_\nu-N_{\nu'}|}}_{\ll{2}}\lesssim D\norm{\nabla\po\th}_{\ll{2}}\lesssim D,
\end{equation}
where we have used {\bf Assumption 2} to estimate $|N_\nu-N_{\nu'}|$
, the fact that $|\nu-\o|\lesssim 2^{j/2}$ on the support of $\eta^\nu_j$ and the fact that $2^{j/2}|\nu-\nu'|\geq 1$. We finally consider the term $\frac{(\th-\th')^2aa'bb'}{|N_\nu-N_{\nu'}|^2}$. We have:
\begin{equation}\label{dd10}
\frac{(\th-\th')^2aa'bb'}{|N_\nu-N_{\nu'}|^2}=\frac{(\th-\th_\nu)^2aa'bb'}{|N_\nu-N_{\nu'}|^2}+2\frac{(\th-\th_\nu)(\th_\nu-\th')aa'bb'}{|N_\nu-N_{\nu'}|^2}+\frac{(\th_\nu-\th')^2aa'bb'}{|N_\nu-N_{\nu'}|^2}.
\end{equation}
The first and the last term in \eqref{dd10} are estimated like the term $\frac{(\nabla\th-\nabla\th')aa'bb'}{|N_\nu-N_{\nu'}|}$ remarking that 
\begin{equation}\label{dd11}
\normm{\frac{(\th-\th_\nu)^2ab}{|N_\nu-N_{\nu'}|^2}}_{\ll{2}}\lesssim \frac{\norm{\po\th}^2_{\ll{4}}}{(2^{j/2}|\nu-\nu'|)^2}\lesssim \norm{\po\th}^2_{\ll{2}}+\norm{\nabla\po\th}^2_{\ll{2}}\lesssim 1,
\end{equation}
and 
\begin{equation}\label{dd12}
\normm{\frac{(\th_\nu-\th')^2ab}{|N_\nu-N_{\nu'}|^2}}_{\ll{2}}\lesssim\norm{\po\th}_{\ll{4}}^2\lesssim \norm{\po\th}^2_{\ll{2}}+\norm{\nabla\po\th}^2_{\ll{2}}\lesssim 1.
\end{equation}
Finally, the second term in \eqref{dd10} is of the form \eqref{dd4} with $H_1=H_2=0$, $H_3=2^{j/2}(\th-\th_\nu)ab$ and $H_4=\frac{(\th_\nu-\th')a'b'}{|N_\nu-N_{\nu'}|}$. 

Thanks to {\bf Assumption 1} and {\bf Assumption 2} on $a, \th$ and the assumption \eqref{thregx1s} on $b$, the estimates \eqref{dd5} and \eqref{dd6} are satisfied. In particular, we have:
\begin{equation}\label{dd13}
\norm{2^{j/2}(\th-\th_\nu)ab}_{\l{\infty}{2}}\lesssim D\norm{\po\th}_{\l{\infty}{2}}\lesssim D(\norm{\po\th}_{\ll{2}}+\norm{\nabla\po\th}_{\ll{2}})\lesssim D,
\end{equation}
and 
\begin{equation}\label{dd14}
\normm{\frac{(\th_\nu-\th')a'b'}{|N_\nu-N_{\nu'}|}ab}_{\l{\infty}{2}}\lesssim D\norm{\po\th}_{\l{\infty}{2}}\lesssim D(\norm{\po\th}_{\ll{2}}+\norm{\nabla\po\th}_{\ll{2}})\lesssim D,
\end{equation}
where we have used the fact that $H^1(\s)$ embeds in $\l{\infty}{2}$ (see \cite{param1} Corollary 3.6 for a proof only using the regularity given by {\bf Assumption 1}). This concludes the proof of Lemma \ref{dd3}. \QED

\section{Proof of Lemma \ref{twodiv}}\label{sec:twodiv}

We need to compute the divergence terms involving $D_1$ and $D_2$ in \eqref{bisoa44}. We start with 
the term involving $D_1$.

\subsection{The divergence term involving $D_1$ in \eqref{bisoa44}}

Using the definition \eqref{bisoa43bb1} of $D_1$ together with the structure equation \eqref{frame} for $N$ and \eqref{ap2}, we obtain:
\begin{equation}\label{ap30}
\begin{array}{l}
\ds\textrm{div}\left(\frac{(N'-g(N,N')N)a'}{1-g(N,N')^2}D_1\right)= \ds\frac{A_1}{(1-\gn^2)(\la-\la'\frac{a}{a'}g(N,N'))}\\
\ds +\frac{A_2\la'}{(1-\gn^2)(\la-\la'\frac{a}{a'}g(N,N'))^2}+\frac{A_3}{(1-\gn^2)^2(\la-\la'\frac{a}{a'}g(N,N'))}\\
\ds +\frac{A_4\la'}{(1-\gn^2)^2(\la-\la'\frac{a}{a'}g(N,N'))^2}+\frac{A_5{\la'}^2}{(1-\gn^2)(\la-\la'\frac{a}{a'}g(N,N'))^3},
\end{array}
\end{equation}
where $A_1, A_2, A_3, A_4, A_5$ are given by:
\begin{equation}\label{ap31}
\begin{array}{ll}
\ds A_1= &\ds a'\nabla_{N-\gn N'}(abb')\trt+aa'bb'\nabla_{N'-\gn N}\trt\\
&\ds +a'\nabla^2(ab)(N,N'-\gn N)b'+\nabla_{\nabla_{N'-\gn N}N}(ab)a'b'\\
&\ds +a'\nabn(ab)\nabla_{N'-\gn N}(b')+(aa'bb'\trt+\nabn(ab)a'b')\\
&\ds \times(\trt'-\gn\trt-\nabn(\gn)+{a'}^{-1}\nabla_{N'-\gn N}(a')),
\end{array}
\end{equation}
\begin{equation}\label{ap32}
\begin{array}{l}
\ds A_2=\nabla^2a(N,N'-\gn N)abb'\gn +\nabla_{\nabla_{N'-\gn N}N}(a)abb'\gn\\
\ds +a'\nabn(a)\nabla_{N'-\gn N}({a'}^{-1}abb')\gn+\nabn(a)abb'\nabn(\gn)\\
\ds +a'\nabla_{N'-\gn N}(a^2{a'}^{-2}bb')\nabn(\gn)+a^2{a'}^{-1}bb'\nabla_{N'-\gn N}\nabn(\gn)\\
\ds +(\trt'-\gn\trt-\nabn(\gn)+{a'}^{-1}\nabla_{N'-\gn N}(a'))\\
\ds \times(a\nabn(a)bb'\gn +a^2{a'}^{-1}\nabn(\gn)bb')+(aa'bb'\trt+\nabn(ab)a'b')\\
\times (\nabla_{N'-\gn N}(a{a'}^{-1})\gn+a{a'}^{-1}\nabla_{N'-\gn N}(\gn)),
\end{array}
\end{equation}
\begin{equation}\label{ap33}
\ds A_3=2a'b'(ab\trt+\nabn(ab))\nabla_{N'-\gn N}(\gn)\gn,
\end{equation}
\begin{equation}\label{ap34}
\begin{array}{ll}
\ds A_4= & \ds 2bb'(a\nabn(a){a'}^{-1}\gn+a^2{a'}^{-2}\nabn(\gn))\\
& \ds \times\nabla_{N'-\gn N}(\gn)\gn,
\end{array}
\end{equation}
and 
\begin{equation}\label{ap35}
\begin{array}{ll}
\ds A_5= & \ds 2(\nabla_{N'-\gn N}(a{a'}^{-1})\gn+a{a'}^{-1}\nabla_{N'-\gn N}(\gn))\\
& \ds\times(a\nabn(a)bb'\gn+a^2{a'}^{-2}\nabn(\gn)bb'). 
\end{array}
\end{equation}
Note that the term $\nabla^2(ab)(N,N'-\gn N)$ appearing in \eqref{ap31} and the term $\nabla^2(a)(N,N'-\gn N)$ appearing in \eqref{ap32} contain at least one tangential derivative. In view of \eqref{ap30}-\eqref{ap35}, one easily checks that the divergence term involving $D_1$ in \eqref{bisoa44} takes the wanted form \eqref{td1} \eqref{td2} provided that we are able to control the two following terms:
\begin{equation}\label{ap36}
\ds \frac{\nabn(\gn)}{(1-\gn)^{1/2}},\, \frac{\nabla_{N-\gn N'}(\gn)}{(1-\gn)^{3/2}}.
\end{equation}
The terms in \eqref{ap36} correspond to the first and the third term of \eqref{ap9}. Thus, this control has already been proved in Appendix A.

\subsection{The divergence term involving $D_2$ in \eqref{bisoa44}}

Using the definition \eqref{bisoa43bb2} of $D_2$ together with the structure equation \eqref{frame} for $N$ and \eqref{ap2}, we obtain:
\begin{equation}\label{ap40}
\begin{array}{l}
\ds\textrm{div}\left(\frac{(N-g(N,N')N')a}{1-g(N,N')^2}D_2\right)= \ds\frac{A_1}{(1-\gn^2)(\la-\la'\frac{a}{a'}g(N,N'))}\\
\ds +\frac{A_2\la'}{(1-\gn^2)(\la-\la'\frac{a}{a'}g(N,N'))^2}+\frac{A_3}{(1-\gn^2)^2(\la-\la'\frac{a}{a'}g(N,N'))}\\
\ds +\frac{A_4\la'}{(1-\gn^2)^2(\la-\la'\frac{a}{a'}g(N,N'))^2}+\frac{A_5{\la'}^2}{(1-\gn^2)(\la-\la'\frac{a}{a'}g(N,N'))^3}, 
\end{array}
\end{equation}
where $A_1, A_2, A_3, A_4, A_5$ are given by:
\begin{equation}\label{ap41}
\begin{array}{ll}
\ds A_1= &\ds a\nabla_{N-\gn N'}(ab)\nabn(b')+a^2b\nabla^2(b')(N,N-\gn N')\\
&\ds +a^2b\nabla_{\nabla_{N-\gn N'}N}(b')\\
&\ds +(a\trt-a\gn\trt'-a\nabla_{N'}(\gn)+\nabla_{N-\gn N'}(a))ab\nabn(b'),
\end{array}
\end{equation}
\bea\label{ap42}
\ds A_2&=&-\nabla_{N-\gn N'}(a^2b)\nabn(a')\gn ab'{a'}^{-2}\\
\nn&& -a^3{a'}^{-2}\nabla^2(a')(N,N-\gn N')\gn bb'\\
\nn&& -a^3{a'}^{-2}\nabla_{\nabla_{N-\gn N'}N}(a')\gn bb'-a^3b\nabn(a')\nabla_{N-\gn N'}(b'{a'}^{-2})\gn\\
\nn&& -a^3{a'}^{-2}\nabn(a')bb'\nabla_{N-\gn N'}(\gn)-a^3{a'}^{-2}\nabn(a')\gn bb'\\
\nn&&\times(a\trt-a\gn\trt'-a\nabla_{N'}(\gn)+\nabla_{N-\gn N'}(a))+a^2b\nabn(b')\\
\nn&& \times(\nabla_{N-\gn N'}(a{a'}^{-1})\gn +a{a'}^{-1}\nabla_{N-\gn N'}(\gn)),
\eea
\begin{equation}\label{ap43}
\ds A_3=2a^2b\nabla_{N-\gn N'}(\gn)\gn\nabn(b'),
\end{equation}
\begin{equation}\label{ap44}
\ds A_4=-2a^3{a'}^{-2}bb'\nabla_{N-\gn N'}(\gn)\gn^2\nabn(a'),
\end{equation}
and 
\begin{equation}\label{ap45}
\begin{array}{ll}
\ds A_5= & \ds -2(\nabla_{N-\gn N'}(a{a'}^{-1})\gn+a{a'}^{-1}\nabla_{N-\gn N'}(\gn))\\
& \ds a^3{a'}^{-2}\nabn(a')\gn bb'.
\end{array}
\end{equation}
Note that the term $\nabla^2(b')(N,N-\gn N')$ appearing in \eqref{ap41} and the term $\nabla^2(a')(N,N-\gn N')$ appearing in \eqref{ap42} contain at least one tangential derivative. In view of \eqref{ap40}-\eqref{ap45}, one easily checks that the divergence term involving $D_2$ in \eqref{bisoa44} takes the wanted form \eqref{td1} \eqref{td2} provided that we are able to control the two terms in \eqref{ap36}. 
This control has already been proved in Appendix A. This concludes the proof of Lemma \ref{twodiv}. \QED



\begin{thebibliography}{10}

\bibitem{Ba-Ch2}
Hajer Bahouri and Jean-Yves Chemin.
\newblock \'{E}quations d'ondes quasilin\'eaires et effet dispersif.
\newblock {\em Internat. Math. Res. Notices}, (21):1141--1178, 1999.

\bibitem{Ba-Ch1}
Hajer Bahouri and Jean-Yves Chemin.
\newblock \'{E}quations d'ondes quasilin\'eaires et estimations de
  {S}trichartz.
\newblock {\em Amer. J. Math.}, 121(6):1337--1377, 1999.

\bibitem{ChBrYo}
Yvonne Choquet-Bruhat and James~W. York, Jr.
\newblock The {C}auchy problem.
\newblock In {\em General relativity and gravitation, {V}ol. 1}, pages 99--172.
  Plenum, New York, 1980.

\bibitem{CoSc}
Justin Corvino.
\newblock Scalar curvature deformation and a gluing construction for the
  {E}instein constraint equations.
\newblock {\em Comm. Math. Phys.}, 214(1):137--189, 2000.

\bibitem{Co}
Justin Corvino and Richard~M. Schoen.
\newblock On the asymptotics for the vacuum {E}instein constraint equations.
\newblock {\em J. Differential Geom.}, 73(2):185--217, 2006.

\bibitem{KR:Duke}
S.~Klainerman and I.~Rodnianski.
\newblock Improved local well-posedness for quasilinear wave equations in
  dimension three.
\newblock {\em Duke Math. J.}, 117(1):1--124, 2003.

\bibitem{Kl:2000}
Sergiu Klainerman.
\newblock P{DE} as a unified subject.
\newblock {\em Geom. Funct. Anal.}, (Special Volume, Part I):279--315, 2000.
\newblock GAFA 2000 (Tel Aviv, 1999).

\bibitem{BIL}
Sergiu Klainerman and Igor Rodnianski.
\newblock Bilinear estimates on curved space-times.
\newblock {\em J. Hyperbolic Differ. Equ.}, 2(2):279--291, 2005.

\bibitem{FLUX}
Sergiu Klainerman and Igor Rodnianski.
\newblock Causal geometry of {E}instein-vacuum spacetimes with finite curvature
  flux.
\newblock {\em Invent. Math.}, 159(3):437--529, 2005.

\bibitem{KR:Annals}
Sergiu Klainerman and Igor Rodnianski.
\newblock Rough solutions of the {E}instein-vacuum equations.
\newblock {\em Ann. of Math. (2)}, 161(3):1143--1193, 2005.

\bibitem{bil2}
Sergiu Klainerman, Igor Rodnianski, and J{\'e}r{\'e}mie Szeftel.
\newblock An ${L}^4_t{L}^4_x$ strichartz estimate for the wave equation on a
  rough background.
\newblock {\em work in progress}.

\bibitem{boundedl2}
Sergiu Klainerman, Igor Rodnianski, and J{\'e}r{\'e}mie Szeftel.
\newblock The bounded ${L}^2$ curvature conjecture.
\newblock {\em preprint, 79 p}, 2012.

\bibitem{Max}
David Maxwell.
\newblock Rough solutions of the {E}instein constraint equations.
\newblock {\em J. Reine Angew. Math.}, 590:1--29, 2006.

\bibitem{SW}
Brian Smith and Gilbert Weinstein.
\newblock Quasiconvex foliations and asymptotically flat metrics of
  non-negative scalar curvature.
\newblock {\em Comm. Anal. Geom.}, 12(3):511--551, 2004.

\bibitem{SmTa}
Hart~F. Smith and Daniel Tataru.
\newblock Sharp local well-posedness results for the nonlinear wave equation.
\newblock {\em Ann. of Math. (2)}, 162(1):291--366, 2005.

\bibitem{stein}
Elias~M. Stein.
\newblock {\em Harmonic analysis: real-variable methods, orthogonality, and
  oscillatory integrals}, volume~43 of {\em Princeton Mathematical Series}.
\newblock Princeton University Press, Princeton, NJ, 1993.
\newblock With the assistance of Timothy S. Murphy, Monographs in Harmonic
  Analysis, III.

\bibitem{param1}
J{\'e}r{\'e}mie Szeftel.
\newblock Parametrix for wave equations on a rough background {I}: Regularity
  of the phase at initial time.
\newblock {\em preprint, 145 p}, 2012.

\bibitem{param2}
J{\'e}r{\'e}mie Szeftel.
\newblock Parametrix for wave equations on a rough background {II}:
  Construction and control at initial time.
\newblock {\em preprint, 84 p}, 2012.

\bibitem{param3}
J{\'e}r{\'e}mie Szeftel.
\newblock Parametrix for wave equations on a rough background {III}: Space-time
  regularity of the phase.
\newblock {\em preprint, 276 p}, 2012.

\bibitem{param4}
J{\'e}r{\'e}mie Szeftel.
\newblock Parametrix for wave equations on a rough background {IV}: Control of
  the error term.
\newblock {\em preprint, 284 p}, 2012.

\bibitem{Ta1}
Daniel Tataru.
\newblock Strichartz estimates for operators with nonsmooth coefficients and
  the nonlinear wave equation.
\newblock {\em Amer. J. Math.}, 122(2):349--376, 2000.

\bibitem{Ta}
Daniel Tataru.
\newblock Strichartz estimates for second order hyperbolic operators with
  nonsmooth coefficients. {III}.
\newblock {\em J. Amer. Math. Soc.}, 15(2):419--442 (electronic), 2002.

\end{thebibliography}
\end{document}